\documentclass[12pt]{article}

\usepackage{amssymb}
\usepackage{latexsym}
\usepackage{amsfonts}
\newtheorem{thm}{Theorem}[section]
\newtheorem{cor}[thm]{Corollary}
\newtheorem{lem}[thm]{Lemma}
\newtheorem{pro}[thm]{Proposition}
\newtheorem{defn}[thm]{Definition}
\newtheorem{ex}[thm]{Example}
\newtheorem{remark}[thm]{Remark}
\title{The groups of Richard Thompson and complexity } 

\author{ Jean-Camille Birget
\thanks{Research supported in part by NSF grant DMS-9970471, and in
part by NSERC grant 216872-1999} 
       }
\date{}
\begin{document}
\maketitle

\begin{flushright}
{\it Dedicated with gratitude to John L.\ Rhodes on his 65th birthday.}
\end{flushright}

\begin{abstract}
We prove new results about the remarkable infinite simple groups 
introduced by Richard Thompson in the 1960s. We give a faithful 
representation in the Cuntz C$^{\star}$-algebra. For the finitely 
presented simple group $V$ we show that the word-length and the table 
size satisfy an  $n \log n$ relation. 
We show that the word problem of $V$ belongs to the parallel
complexity class AC$^1$ (a subclass of {\bf P}), whereas the 
generalized word problem of $V$ is undecidable. 

We study the distortion functions of $V$ and show that $V$ contains all 
finite direct products of finitely generated free groups as subgroups with 
{\it linear} distortion. As a consequence, up to polynomial equivalence of 
functions, the following three sets are the same:
the set of distortions of $V$, the set of Dehn functions of 
finitely presented groups, and the set of time complexity functions of 
nondeterministic Turing machines. 
\end{abstract}


\section{Introduction}

\noindent In \cite{Th0} Thompson constructed a simple finitely presented 
infinite group, one of the most remarkable groups ever found.  He denoted 
it by Pa($^{\omega}2$) and by Ft$(K)$ in \cite{Th}, by $\hat{V}$ in 
\cite{Th0}, and by ${\mathfrak C}'$ in \cite{McKTh}; in \cite{CFP} it is 
denoted by $V$, and we will follow that convention, which has been widely
adopted.
We will also use the uncountable Thompson group ${\mathcal G}_{2,1}$ 
(following the notation of \cite{ESc}). 
The proofs of the main properties of $V$ were first outlined 
in \cite{Th} and can be found in detail in \cite{CFP}, or in \cite{Hig74} 
(where it is called $G_{2,1}$, as part of an infinite
family of finitely presented simple groups). 

Thompson defined his groups as permutation groups of certain sets of 
infinite words over the alphabet $\{ 0, 1\}$. We will follow \cite{ESc}, 
and indirectly \cite{Hig74}, and define $V$ by partial bijections of the 
free monoid $\{a, b\}^*$. The advantage of this definition is that partial 
actions on finite words enable us to define algorithmic problems and their 
complexity. From now on, ``word'' will mean ``finite word''. 
 
\medskip

Our setting for the Thompson groups requires a number of elementary 
definitions and facts. Almost all of these concepts are standard  
(the literature on Thompson groups suffers from idiosyncratic terminology, 
which can usually be avoided). Since the Thompson groups are based on
partial actions, we have to choose a side for the actions.
We choose to {\it act on the left}.  The main advantage of this choice 
will turn out to be the connection between Thompson groups acting on the 
left, and prefix codes. The literature on codes greatly prefers prefix 
codes over suffix codes, and we hope that this choice improves readability.

Let $A$ be a finite alphabet. The set of all words over
$A$ (including the empty word $\varepsilon$) is denoted by $A^*$. 
Concatenation of two words $u, v \in A^*$ is denoted by $u\cdot v$ or simply 
$uv$; $A^*$ is a monoid under the concatenation operation. For two sets 
$X_1, X_2 \subseteq A^*$, we denote their concatenation by $X_1X_2$ or by  
$X_1 \cdot X_2$, defined by \  
$X_1X_2 = \{x_1x_2 \in A^* : x_1 \in X_1, x_2 \in X_2\}$. From now on we 
assume that the alphabet $A$ has a least two letters.  

A {\it right ideal} of $A^*$ is defined to be a subset $R \subseteq A^*$ such 
that \,  $R \cdot A^* \subseteq R$ \, (i.e., $R$ is closed under 
multiplication by any word in $A^*$ on the right). 

For two words $u, v \in A^*$, we say that $u$ is a {\it prefix} of $v$ iff 
$v = ux$ for some $x \in A^*$; we also write $u \geq_{\rm pref} v$ or 
$v \leq_{\rm pref} u$; this is a partial order, related to set inclusion by 
the fact that \, $v \leq_{\rm pref} u$ iff $vA^* \subseteq uA^*$.
We say that $u$ is a {\it strict prefix} of $v$ (and write $u >_{\rm pref} v$) 
iff $u \geq_{\rm pref} v$ and $u \neq v$. We say that $u$ and $v$ are 
{\it prefix-comparable} iff $v \leq_{\rm pref} u$ or $u \leq_{\rm pref} v$;
we denote this by $u \lesseqgtr_{\rm pref} v$.  
A {\it prefix code} over $A$ is defined to be a subset $C$ of $A^*$ 
such that no element of $C$ is a strict prefix of another element of $C$.  
The monograph \cite{BerstelPerrin} is an excellent reference for the material 
on prefix codes that we use here.
By definition, a {\it maximal prefix code} over an alphabet $A$ is a prefix 
code over $A$ which is not a strict subset of any other prefix code over $A$.

For a right ideal $R$ of $A^*$, a set $\Gamma \subseteq R$ is called a set of 
{\it right-ideal generators} of $R$ (or just ``generators'') iff \  
$R = \Gamma \cdot A^*$.  
One can prove (see Lemma \ref{ideals}) that any right ideal $R$ of $A^*$ 
has a unique minimal (under inclusion) set of right-ideal generators, and this 
set of generators is a prefix code.
Hence, prefix codes could be called ``right-ideal bases''. 
Moreover, since the prefix code of a right ideal is unique, right ideals of 
$A^*$ and prefix codes over $A$ are in one-to-one correspondence
(see Lemma \ref{ideals} in Appendix A1 for proofs).
 
A right ideal $R$ of $A^*$ is said to be {\it finitely generated} iff the 
prefix code corresponding to $R$ is finite.
A right ideal $R$ of $A^*$ is called {\it essential} iff $R$ has a non-empty 
intersection with every right ideal of $A^*$. (This is C$^{\star}$-algebra 
terminology.)
Note that a right ideal $R$ of $A^*$ is essential if its prefix code is 
maximal; this is also equivalent to saying that $R$ is a right ideal such 
that for every $u\in A^*$ there is $x \in A^*$ such that $ux \in R$ 
(see Lemma \ref{ideals}). Clearly, if a right ideal $J$ contains an 
essential right ideal of $A^*$, then $J$ is essential too.

\smallskip

By definition, a {\it right-ideal homomorphism} of $A^*$ is a function 
$\varphi: R_1 \to R_2$ such that $R_1$ and $R_2$ are right ideals of $A^*$, 
and such that for all 
$u \in R_1$ and all $x \in A^*$: \ $\varphi(u) \cdot x = \varphi(ux)$. 
A {\it right-ideal isomorphism} of $A^*$ is a  bijective right-ideal 
homomorphism. 

One can prove (see Lemma \ref{bijPrefCode}) that the set of all right-ideal 
homomorphisms (or isomorphisms) of $A^*$ is in one-to-one correspondence 
with the set of all functions (respectively  bijections) between prefix codes 
of $A^*$. For a right-ideal isomorphism $\varphi: P_1A^* \to P_2A^*$, where 
$P_1$ and $P_2$ are prefix codes, the restriction \ 
$\tau_{\varphi}: P_1 \to P_2$ \ is a bijection, and $\tau_{\varphi}$ 
determines $\varphi$ uniquely. Following Thompson, the restriction 
$\tau_{\varphi} : P_1 \to P_2$ of $\varphi$ will be called the 
{\it table of} \ $\varphi$, and will be used to represent 
$\varphi$ by a traditional (finite) function table. (In \cite{Hig74} and 
\cite{ESc} this was called the ``symbol of $\varphi$''). 
The maximal prefix code $P_1$ is called the {\it domain code} of $\varphi$,
and $P_2$ is called the {\it image code} or {\it range code} of $\varphi$.

By definition, an {\it extension} of a right-ideal isomorphism 
$\varphi: R_1 \to R_2$ is a right-ideal isomorphism 
$\Phi: J_1 \to J_2$ where $J_1, J_2$ are right ideals such that 
$R_1 \subseteq J_1, \ R_2 \subseteq J_2$, and $\Phi$ agrees with $\varphi$ 
on $R_1$ (i.e., $\Phi(x) = \varphi(x)$ for all $x \in R_1$). In that case we 
also call $\varphi$ a {\it restriction} of $\Phi$.  
The extension, and the restriction, are called strict iff $\varphi \neq \Phi$. 

A right-ideal isomorphism is said to be {\it maximal} iff it has no strict 
extension in $A^*$; it is called {\it extendable} otherwise.   
We denote the maximum extension of $\varphi$ by ${\sf max} \, \varphi$; \ we 
will prove in Lemma \ref{uniqueMaxExt} that the {\em maximum} extension
of an isomorphism between essential right ideals is unique (for this 
uniqueness, it is necessary that the ideals be essential).

\medskip

Most of the above concepts can be pictured, using trees. The monoid $A^*$ can 
be described by the Cayley graph of the right regular representation of $A^*$ 
relative to the generating set $A$. We will simply call this {\it the tree of} 
$A^*$.  It is an infinite tree rooted at the empty word $\varepsilon$. 
Every vertex has $|A|$ children.  
Every subset of $A^*$  is pictured as a set of vertices of this infinite tree. 
A prefix code is pictured as a set of vertices, no two of
which lie on a same directed path from the root.
A finite prefix code is maximal iff it is a prefix code that forms a ``cut'' 
in the tree (i.e., a set of vertices whose removal disconnects the root of the 
tree from all the ``ends'' of the tree). Infinite maximal prefix codes are 
harder to visualize, but the following concept is useful: 

For any prefix code $P \subset A^*$ ($P \neq \emptyset$), the 
{\it prefix tree} of $P$ is defined to be the subtree of the tree of $A^*$, 
whose vertex subset consists of all 
the prefixes of words in $P$ (and whose root is still $\varepsilon$). 
Hence, the set of leaves of this subtree is $P$. 
We have the following general facts about non-empty subsets $P \subseteq A^*$ 
(finite or infinite):

\smallskip

\noindent $\bullet$ \ 
{\it $P$ is a prefix code \ iff \ $P$ is the set of leaves of a subtree of the 
tree of $A^*$.}

\smallskip

\noindent $\bullet$ \ 
{\it A prefix code $P$ is maximal \ iff \ every non-leaf vertex of the prefix 
tree of $P$ has exactly $|A|$ children (in the prefix tree of $P$).} 

\medskip

A right ideal of $A^*$ is the same thing as an {\it order ideal} relative to 
the prefix order $\leq_{{\rm pref}}$ in $A^*$.
A (maximal) prefix code in $A^*$ is the same thing as a {\it (maximal) 
anti-chain} relative to the prefix order $\leq_{{\rm pref}}$. 
Right-ideal isomorphisms are the same thing as prefix-order isomorphisms
between prefix-order ideals. So, our discussion could also be carried out in 
partial-order terminology.  

Prefix codes are well known; see e.g. \cite{Hamming}, \cite{BerstelPerrin}.
They are not only of mathematical interest but are used in practice (e.g., in 
text compression by Huffman coding, and for error correcting codes). 

\begin{ex}
\label{exPrefCod}
 \ --- \ Some infinite maximal prefix codes
\end{ex}
\vspace{-.1in}
Infinite maximal prefix codes can be extremely complex.
Here are some examples.

\smallskip

\noindent {\bf (1)} \ For any fixed infinite sequence \
$(a_1, a_2, \ldots, a_{n-1}, a_n, \ldots) \in A^{\omega}$ \ one can build an
infinite maximal prefix code as follows. For any $a \in A$, let
$\sigma(a) \in A$ be another letter $(\neq a)$ chosen in $A$.
Consider the code \ \
$P = \{a_1 \ldots a_{n-1} \, \sigma(a_n) \ : \ n \geq 1\}.$ \   
Such infinite prefix codes have ``one infinite path-shaped end''.

\smallskip

\noindent {\bf (2)} \ 
{\it Combination of prefix codes:} Let $X = \{x_i : i \in I\}$ be a
prefix code,  and let  $(Q_i : i \in I)$ be a family of prefix codes, 
with $|X| = |I|$ (finite or infinite).
Then \ $\bigcup_{i \in I} x_i Q_i $ \ is a prefix code, which is maximal
if $X$ and each $Q_i$ are maximal.
This enables us to construct maximal prefix codes with any number of
``infinite path-shaped ends''.

\smallskip

\noindent {\bf (3)} \ An infinite maximal prefix code does not need to have 
any ``path-shaped ends''; instead, it could have any number of infinite 
``tree-shaped ends''.
For example, consider the following code over $\{a,b\}$:
 \ $P = \{a^2, b^2\}^*\cdot \{ab, ba\}$.  \  
By looking at the prefix tree of $P$ it is easy to see that $P$ is
a prefix code (all the words in $P$ are leaves of the prefix tree)
and that it is maximal (all non-leaves have two children).
Another example of a similar infinite maximal prefix code is
 \ $\{a^2, ab, b^2\}^*\cdot \{ba\}$.
See \cite{BerstelPerrin} for more examples.


\section{The Thompson groups}

Before defining the Thompson groups we prove a few facts about isomorphisms of
right ideals of $A^*$. Proposition \ref{uniqueMaxExt} and Lemmas 
\ref{ExtCritFin} and \ref{ImDom} appear in Thompson's work and in \cite{ESc}, 
with a similar content but a different formalism. Lemmas \ref{ExtCritInf}
and \ref{findMaxExt} are new.

\begin{pro}
\label{uniqueMaxExt} \  
An isomorphism between essential right ideals of $A^*$ has a 
{\em unique} maximum extension.

Equivalently, if two isomorphisms $\varphi_1, \varphi_2$ between essential 
right ideals agree on an essential right ideal then $\varphi_1$ and 
$\varphi_2$ have the same maximum extension.
\end{pro}
{\bf Proof.} \ Let $\varphi: P_1 A^* \to P_2 A^*$ be an isomorphism of 
essential right ideals, where $P_1$ and $P_2$ are maximal prefix codes.  
If $\varphi(x)$ is not defined for some $x \in A^*$ then (by Lemma 
\ref{ideals} (4)), there exists $p \in P_1$ with $x >_{\rm pref} p$. Let 
$p_x$ be the first element in the lexicographic order (assuming we have 
chosen a fixed total order for the finite alphabet $A$) such that 
$p_x \in P_1$ and $x >_{\rm pref} p_x = x \,u_x$ \, (for some $u_x \in  A^*$).
Then $p_x$ and $u_x$ are uniquely determined by $x$ and $\varphi$. 

If there is an extension $\Phi$ of $\varphi$ such that $\Phi(x)$ is defined, 
then $\varphi(p_x) = \Phi(p_x) = \Phi(x) \, u_x$. Hence, $\Phi(x)$ is uniquely
determined by $x$ and $\varphi$. 

It follows from this the union of extensions of $\varphi$ is a well 
defined extension too. Thus, we can take the union of all extensions of 
$\varphi$ to obtain the maximum extension of $\varphi$. 
 \ \ \ $\Box$

\medskip

The next two Lemmas give useful characterizations of extendability and 
maximality of right-ideal isomorphisms of essential right ideals. 
They will be used in the next section.

\begin{lem} 
\label{ExtCritFin} \  
Let $\varphi: P_1 A^* \to P_2 A^*$ be an isomorphism of essential right ideals, 
where $P_1$ and $P_2$ are {\em finite} maximal prefix codes. Then $\varphi$ is 
extendable iff there are $x_0, y_0 \in A^*$ such that for every letter 
$\alpha \in A$: \ \  $x_0 \alpha \in P_1$, $y_0 \alpha \in P_2$, and 
$\varphi(x_0 \alpha) = y_0 \alpha$. 

(If this condition holds, $\varphi$ can be extended by mapping $x_0$ to $y_0$.)
\end{lem}
{\bf Proof.} \ If $x_0 \alpha \in P_1$, $y_0 \alpha \in P_2$ and 
$\varphi(x_0 \alpha) = y_0 \alpha$ for every letter $\alpha \in A$,
then $\varphi$ can be extended by defining $\varphi(x_0)$ to be $y_0$. 
The prefix code of the domain then becomes $P_1 \cup \{x_0\} - x_0A$, 
and the prefix code of the range becomes $P_2 \cup \{y_0 \} - y_0A$.

Conversely, suppose $\varphi$ can be strictly extended to $\Phi$. 
Consider a word $x_0$ on which $\varphi$ is not defined, but on which $\Phi$ 
is defined.

\smallskip

\noindent Case 1: \ If for all $\alpha \in A$, $\varphi(x_0 \alpha)$ is 
defined, i.e., $x_0 \alpha \in P_1A^*$, then actually $x_0 \alpha \in P_1$ 
(since $\varphi$ is not defined on any strict prefix of $x_0 \alpha$). 
Also, $\varphi(x_0 \alpha) = \Phi(x_0 \alpha) =  \Phi(x_0) \alpha$. 
So we pick $y_0$ to be  $\Phi(x_0)$. Then $y_0\alpha \in P_2A^*$, for all 
$\alpha \in A$; 
but actually, $y_0 \alpha \in P_2$ (since $\varphi^{-1}$ is not defined on any 
strict prefix of $y_0 \alpha$).  Now $x_0$ and $y_0$ satisfy the properties of 
the Lemma.

\smallskip

\noindent Case 2: \ If for some $\alpha \in A$, $\varphi(x_0 \alpha)$ is not
defined, we replace $x_0$ by $x_0 \alpha$ and continue the reasoning. 
Eventually, we reach case 1, since $P_1$ is finite and maximal.
 \ \ \  $\Box$ 

\medskip

Note that the Lemma is not always true for {\it infinitely} generated 
essential right ideals. For example, let $A = \{a,b\}$, consider the 
maximal prefix code $P_1 = P_2 = \{a^n b : n \in {\mathbb N}\}$, and let 
$\varphi$ be the identity map on \ 
$P_1 A^* = \{a^n b : n \in {\mathbb N}\} \cdot \{a,b\}^*= \{a,b\}^*- \{a\}^*$.
Then $\varphi$ can obviously be extended to the identity map on $\{a,b\}^*$,
but there is no word $x$ such that $xa, xb \in P_1$.  

\medskip

The following Lemma gives a criterion for extendability in the general case.

\begin{lem} 
\label{ExtCritInf} \   
Let $\varphi: P_1 A^* \to P_2 A^*$ be an isomorphism of essential right ideals, 
where $P_1$ and $P_2$ are arbitrary maximal prefix codes. Then $\varphi$ is 
extendable iff there exists a maximal prefix code $Q \subseteq A^*$ with 
$|Q| > 1$, and there are $x_0, y_0 \in A^*$ such that for all $q \in Q:$ \ \ 
$x_0q \in P_1$, \ $y_0q \in P_2$, and $\varphi(x_0q) = y_0q$.

If $P_1$ and $P_2$ are finite then is is enough to consider finite 
codes $Q$.
\end{lem}
One sees that the general Lemma differs from the finite case by the fact that
all possible maximal prefix codes are used (instead of the alphabet $A$, 
which is a very special maximal prefix code). 

\smallskip

\noindent We will use the following notation:
For any set $L \in A^*$ and any word $x \in A^*$, we define 

\smallskip

 \ \ \ \ $\overline{x}L = \{ w \in A^* : xw \in L\}$. 

\medskip

\noindent
{\bf Proof.} \ If the condition in the Lemma holds (i.e., 
$\varphi(x_0q) = y_0q$ for all $q \in Q$), then $\varphi$ can be extended 
by defining the image of $x_0w$ to be $y_0w$ (for all $w \in A^*$). 

Conversely, suppose $\varphi$ can be extended to an isomorphism of essential
right ideals $\Phi$. Consider a word $x_0$ on which $\varphi$ is not defined, 
but on which $\Phi$ is defined, and suppose $\Phi(x_0) = y_0$. 
Define $Q$ as follows:

\smallskip

 \ \ \ \ $Q = \{ w \in A^* : x_0w \in P_1 \} = \overline{x_0}P_1$. 

\smallskip

\noindent 
By Lemma \ref{Q} in Appendix A1, $Q$ is a maximal prefix code.

\smallskip

\noindent {\it Claim:} \ \ \  
$\overline{x_0}P_1 = \overline{y_0}P_2 \ \ ( \ = \ Q)$. \ 
{\it Hence, $y_0q \in P_2$ for all $q \in Q$.} 

Indeed, $w \in \overline{x_0}P_1$ iff $x_0w \in P_1$, iff 
$\varphi(x_0w) \in P_2$. 
Moreover, $\varphi(x_0w) = \Phi(x_0w) = \Phi(x_0) \, w = y_0w$. 
Hence, $w \in \overline{x_0}P_1$ iff $y_0w \in P_2$; the latter holds iff
$w \in \overline{y_0}P_2$. This proves the Claim.

\smallskip

Finally, for any $q \in Q$, \ 
$\varphi(x_0q) = \Phi(x_0q) = \Phi(x_0) \, q = y_0q$.  \ \ \ $\Box$

\begin{lem}
\label{findMaxExt} \ Let $\varphi: P_1A^* \to P_2A^*$ be as in the previous
Lemma; then the maximum extension of $\varphi$ can be obtained as follows.
There are two maximal prefix codes
$\{x_i : i \in I\}, \{y_i : i \in I\} \subset A^*$ (for an index set
$I \subseteq \mathbb{N}$), such that

\smallskip

\noindent {\bf (a)} \ for each $i \in I$ there is a maximal prefix code $Q_i$
such that $x_iQ_i \subseteq P_1$, $y_iQ_i \subseteq P_2$, and for all
$q \in Q_i$, \ $\varphi(x_iq) = y_iq$;

\smallskip

\noindent {\bf (b)} \ the sets $x_iQ_i$ and $y_iQ_i$ are
$\subseteq$-{\em maximal}.

\smallskip

Then $\Phi$, defined by $x_i \mapsto y_i$ (for all $i \in I$), is the
maximum extension of $\varphi$.
\end{lem}

Figure 1 below gives the tree picture representing the prefix code $P_1'$
of $\Phi: P_1'A^* \to P_2'A^*$, the prefix code $P_1$ of 
$\varphi : P_1A^* \to P_2A^*$, an element $x_i \in P_1'$, and the prefix 
code $Q_i$ used to extend $\varphi$ to $\Phi$, at $x_i$.
To extend $\varphi$ to $\Phi$, we need a prefix code $Q_i = \overline{x_i}P_1$
for each $x_i \in P_1'$. In this picture, $x_i$ should be viewed as the root 
of the tree for $Q_i$.

\newpage

\unitlength=0.90mm \special{em:linewidth 0.4pt}
\linethickness{0.4pt}

\begin{picture}(70,55)
\put(70,50){\line(-1,-1){50}}
\put(70,50){\line(1,-1){50}}
\put(45,30){\makebox(0,0)[cc]{$P_1'$}}
\qbezier(50,30)(54,23)(60,28)
\qbezier(60,28)(70,35)(90,30)

\put(58.5,30.5){\makebox(0,0)[cc]{$x_i$}}
\put(60,28){\circle*{1}}     
\put(60,28){\line(-1,-1){30}}
\put(60,28){\line(1,-1){30}}
\put(25,10){\makebox(0,0)[cc]{$P_1$}}
\qbezier(30,10)(34,12)(40,8) 
\qbezier(40,8)(62,4)(80,8)
\qbezier(80,8)(100,12)(116,4)
\qbezier(39.8,7.8)(62,3.5)(80.2,7.8)
\put(60,2){\makebox(0,0)[cc]{$x_iQ_i = x_i(\overline{x_i}P_1)$}}
\end{picture}

\bigskip

\bigskip

{\sf Fig.~1: \ Extendability condition} 

\bigskip

\noindent 
{\bf Remarks.} \ {\bf (0)} In Lemma \ref{findMaxExt} we allow $Q_i$ to 
consist of just the empty word.    

\smallskip

\noindent
{\bf (1)} {\it $\subseteq$-Maximality} of the sets $x_iQ_i$ and $y_iQ_i$ 
($i \in I$) is defined as follows:
Suppose there exist a maximal prefix code $Q$ and words $x, y$ such that
$xQ \subseteq P_1$, $yQ \subseteq P_2$,  and for all $q \in Q$:
$\varphi(xq) = yq$. And suppose there exists $j \in I$
such that $x_jQ_j \subseteq xQ$ and $y_j Q_j \subseteq yQ$.
Then $x = x_j, y = y_j, Q = Q_j$.

\smallskip

\noindent  
{\bf (2)} When $\varphi$ is already maximum then the Lemma 
holds with $\{x_i : i \in I\} = P_1$, $\{y_i : i \in I\} = P_2$ and 
$Q_i = \{ \varepsilon \}$ for all $i \in I$ (where $\varepsilon$ denotes 
the empty string).   

\smallskip

\noindent
{\bf (3)} The sets $x_iQ_i$ (for $i \in I$) are two-by-two disjoint (and 
similarly for the sets $y_i Q_i$). So all the extensions, from  $x_iQ_iA^*$ to 
$x_iA^*$, as $i$ ranges over $I$, are independent (for different $i$'s), and 
can be viewed as being carried out in parallel.

\smallskip

\noindent
{\bf Proof} of Lemma \ref{findMaxExt}. \ 
If $\varphi$ is already maximal there is nothing to
prove (by Remark (2)). If $\varphi$ can be extended then, by Lemma 
\ref{ExtCritInf}, the words $x_i$, $y_i$ and the sets $Q_i$ exist, satisfying 
the claimed properties (a) and (b), and the map $\Phi$ obtained is an 
extension of $\varphi$. 
We have to show that $\Phi$ is the maximum extension of $\varphi$. 

After each set $x_iQ_i$ has been replaced by $x_i$ (and each $y_iQ_i$ by $y_i$)
in the extension process, the domain code of $\Phi$ is $\{x_i : i \in I\}$
and the image code of $\Phi$ is $\{y_i : i \in I\}$. 

The following Claim implies, by Lemma \ref{ExtCritInf}, that $\Phi$
cannot be extended.

\smallskip

\noindent {\it Claim.} 
The sets $\{x_i : i \in I\}$, $\{y_i : i \in I\}$ have no strict subset of 
the form $xP$, respectively $yP$, with $\Phi(xp) = yp$ for all $p \in P$ 
(where $P$ is a maximal prefix code with more than one element).

\smallskip

\noindent
Proof of the Claim: If, by contradiction, the Claim is false then there
is a non-trivial maximal prefix code $P$ such that \   
$xP = \{x_j : j \in J\} = \{ xp_j : j \in J\}$ \ and  
$yP = \{y_j : j \in J\} = \{yp_j : j \in J\}$, \ for some $J \subseteq I$.
Then \ $\bigcup_{j \in J} xp_jQ_j = x \bigcup_{j \in J} p_jQ_j$; moreover,   
 \ $\bigcup_{j \in J} p_jQ_j = Q$ \ is a maximal prefix code 
(by construction (2) in Example \ref{exPrefCod}), with 
$xQ \subseteq P_1$ and $yQ \subseteq P_2$. Now, $xQ$ contains \ 
$x_{j_1}Q_{j_1} \cup x_{j_2}Q_{j_2} \cup  \dots$ \ (for some 
$j_1, j_2, \ldots \in  J \subseteq I$ with $j_1 \neq j_2$), 
which contradicts $\subseteq$-maximality (assumption (b)).
 \ \ \ $\Box$

\begin{lem}
\label{ImDom} \  
Let $\varphi_1$ and $\varphi_2$ be right-ideal isomorphisms between essential 
right ideals of $A^*$. Then $\varphi_1$, $\varphi_2$, have restrictions
$\varphi_1'$, respectively $\varphi_2'$, 
such that the range of $\varphi_1'$ is equal to the domain of $\varphi_2'$.

If the domain and ranges of $\varphi_1$ and $\varphi_2$ are finitely generated
then so are the domains and ranges of $\varphi_1'$ and $\varphi_2'$. 
\end{lem}
{\bf Proof.} \ 
Let $\varphi_1: R_1 \to Q_1$ and $\varphi_2:  R_2 \to Q_2$ be right-ideal 
isomorphisms between essential right ideals. We want to show that
there is a restriction $\varphi_1'$ of $\varphi_1$, and a restriction 
$\varphi_2'$ of $\varphi_2$, such that

\smallskip

 \ \ \ \ \ \ \ \ \ \ \ \ \ 
$R_1' \ \stackrel{\varphi_1'}{\longrightarrow} \ S' \    
\stackrel{\varphi_2'}{\longrightarrow} \ Q_2'$  

\smallskip

\noindent 
where $R_1' \ (\subseteq R_1), \ Q_2' \ (\subseteq Q_2$), and 
$S' \ (\subseteq Q_1 \cap R_2)$ are essential right ideals, and 
$\varphi_1', \varphi_2'$ are right-ideal isomorphisms. 

Let $S' = Q_1 \cap R_2$, $R_1' = \varphi_1^{-1}(S')$, and
$Q_2' = \varphi_2'(S')$. Also let $\varphi_1'$ be the restriction of 
$\varphi_1$ to $R_1'$ and let $\varphi_2'$ be the restriction of $\varphi_2$
to $S'$. 

The intersection of two right ideals $Q_1, \, R_2$, is obviously a right 
ideal, and $Q_1 \cap R_2$ is essential if $Q_1, \, R_2$ are essential (Lemma
\ref{IntersectTrans}). Also, $\varphi_1^{-1}(S')$ and $\varphi_2'(S')$
are essential, by Lemma \ref{referee1}.
The other properties of the Lemma are straightforward.

If $R_1, Q_1, R_2, Q_2$ are finitely generated right ideals then 
$S' = Q_1 \cap R_2$ is finitely generated, by Lemma \ref{idealInter}.
Moreover, $R_1' = \varphi_1^{-1}(S')$ and $Q_2' = \varphi_2(S')$
are finitely generated since $\varphi_1$ and $\varphi_2$ are isomorphisms. 
 \ \ \ $\Box$

\bigskip

The following definition of the Thompson groups is very close to the 
definition of Scott \cite{ESc} (and indirectly, to the definition in 
\cite{Hig74}).
The tree representation of codes connects this definition and the definition 
by action on finite trees used in \cite{CFP}.  Although the Thompson 
groups are traditionally defined with the alphabet $\{0,1\}$, we prefer 
to use $\{a,b\}$ (the symbols ``0'' and ``1'' have too many meanings already).

\begin{defn} \  
The {\bf Thompson group} $V$ is the partial action group on 
$\{a,b\}^*$ consisting of all maximal isomorphisms between finitely 
generated essential right ideals of $\{a,b\}^*$.

\smallskip 

\noindent The {\bf Thompson group} ${\mathcal G}_{2,1}$ is the partial action
group on $\{a,b\}^*$ consisting of all maximal isomorphisms between 
essential right ideals of $\{a,b\}^*$.

\smallskip

\noindent {\bf Multiplication:} For $\varphi, \psi \in {\mathcal G}_{2,1}$ (or 
$\in V$), the product $\varphi \cdot \psi$ is \ 
{\sf max}$(\varphi \circ \psi)$ \ 
(i.e., the maximum extension of the composition of $\psi$ and $\varphi$, where
$\psi$ is applied first).
\end{defn}

We will usually write just $\varphi \psi$ for $\varphi \cdot \psi$. One can 
check easily that $\varphi^{-1} \varphi = \varphi \varphi^{-1} = {\bf 1}$ \ 
(the identity map on $\{a,b\}^*$). Maximum extension is needed for this to be 
true; without taking the maximum extension, the composite 
$\varphi^{-1} \circ \varphi$ is the restriction of {\bf 1} to the domain of 
$\varphi$, and $\varphi \circ \varphi^{-1}$ is the restriction of {\bf 1} to 
the image of $\varphi$. It is also easy to check that for any (not necessarily 
maximal) isomorphisms of essential right ideals, \ 
${\sf max}(\psi \circ \varphi) = $
${\sf max}({\sf max} \, \psi \ \circ \ {\sf max} \, \varphi)$. 
Associativity follows from this. Hence, ${\mathcal G}_{2,1}$ and $V$ are 
groups.

In connection with the definition of the Thompson groups it is natural to
introduce the following terminology: Two isomorphisms $\varphi$ and $\psi$
between essential right ideals of $A^*$ are {\it congruent} iff $\varphi$ and
$\psi$ have the same maximum extension 
(${\sf max} \, \varphi = {\sf max} \, \psi$).  

\medskip

\noindent As mentioned before, Thompson showed the following: \\  
$\bullet$ \ 
      {\it ${\mathcal G}_{2,1}$ and $V$ are simple}; \\  
$\bullet$ \ 
      {\it $V$ is finitely presented};  \\  
$\bullet$ \ 
      {\it $V$ contains all finite groups as subgroups.} 

\noindent The first two facts are not obvious at all (see \cite{Th0}, 
\cite{Hig74}, \cite{CFP}); the third fact is straightforward (but also 
remarkable). 

\bigskip

It is not obvious to construct examples of finitely generated groups which, 
like $V$, contain all finite groups. 

\medskip

\noindent {\bf Example } \ --- \   
{\it Another finitely generated group containing all finite groups:} 

\smallskip

\noindent  
Let ${\mathfrak S}_{\mathbb Z}$ be the set of all permutations of ${\mathbb Z}$
(the integers), and let ${\mathfrak S}_{{\rm fin}{\mathbb Z}}$ be the set of
all {\it finitary} permutations of ${\mathbb Z}$ (i.e., permutations that
fix all but a finite number of integers). Let \ 
$\sigma : z \in {\mathbb Z} \mapsto z+1$ \ be the right-shift function. 
Then the group \ 
$G = \langle {\mathfrak S}_{{\rm fin}{\mathbb Z}} \cup \{ \sigma \} \rangle$ \ 
(i.e., the subgroup of ${\mathfrak S}_{\mathbb Z}$ generated by 
${\mathfrak S}_{{\rm fin}{\mathbb Z}}$ and $\sigma$)
contains all finite groups. Moreover, $G$ is generated by $\sigma$ and the 
transposition $\tau_{(0,1)}$, so $G$ is finitely generated. This group has
been known for a long time; it is less well known, but easy to prove that 
 \ $\langle {\mathfrak S}_{{\rm fin}{\mathbb Z}} \cup \{ \sigma \} \rangle$ \   
{\em is a subgroup of} $V$.

To show that 
$G = \langle {\mathfrak S}_{{\rm fin}{\mathbb Z}} \cup \{ \sigma \} \rangle$ 
is a subgroup of $V$ we use a one-to-one correspondence 
between ${\mathbb Z}$ and the maximal prefix code \ 
$a^*ab \cup b^*ba \subset \{a,b\}^*$, defined by
$$z \mapsto \left\{ \begin{array}{ll} 
                  a^{-z} ab \ & \ \mbox{if $z \leq 0$} \\  
                  b^z a    \  & \ \mbox{if $z > 0$} 
             \end{array} \right.     $$
It follows that $\tau_{(0,1)}$ is represented by the following 
element  of ${\mathcal G}_{2,1}$:  \ $ab \mapsto ba, \ ba \mapsto ab$, \ and 
$\tau_{(0,1)}$ is the identity elsewhere on $a^*ab \cup b^*ba$. 
The shift $\sigma$ is represented by $a^{-z}ab \mapsto a^{-(z+1)}ab$ (for all 
$z < 0$), $ab \mapsto ba$, and $bb^za \mapsto bb^{z+1}a$ (for all $z \geq 0$). 
Here we just indicate how the maps are defined on the maximal prefix code; 
the definition on the corresponding essential right ideal follows 
automatically. 

Maximum extension of these two maps reveals that they actually belong to $V$.  
The map representing $\tau_{(0,1)}$ is easy to extend to \ 
$aa \mapsto aa, \ bb \mapsto  bb, \ ab \mapsto ba, \ ba \mapsto ab$. 
(Again, we just indicate the map on a maximal prefix code.)
The shift can be extended to \ 
$aa \mapsto a, \ ab \mapsto ba, \ b \mapsto bb$ \ 
(as defined on maximal prefix codes); we are using Lemma \ref{ExtCritInf}. 

\smallskip

Notice that the representation of the right-shift $\sigma$ above is (the 
inverse of) the generator of $V$, called ``$A$'' in \cite{CFP} (see also 
the remarks below on other Thompson groups). This provides a nice 
interpretation of the generator ``$A$''.

\bigskip

\noindent {\bf Remark on partial actions:} \ It would not be correct to say 
that ${\mathcal G}_{2,1}$ and $V$ act on $\{a,b\}^*$ by partial maps, 
since in addition to the composition of the partial maps we also take the
maximum extension. 
What we have here is a {\it partial action} (see \cite{LawsonInv}), 
as opposed to an (ordinary) action by partial maps. Here, 
with each element $g$ of the group one associates a partial 
transformation  $\tau(g)$ on some chosen set such that: \ For the identity
of $G$ we have \ $\tau(1_G) = ${\bf 1}; \  for all $g \in G$, \
$\tau(g^{-1}) = \tau(g)^{-1}$; \ and for all $g_1, g_2 \in G$, \  
$\tau(g_1) \circ \tau(g_2) \subseteq \tau(g_1g_2)$. 

On the other hand, Thompson \cite{Th} used an ordinary action, by total 
permutations, but on infinite words. But from a computational point of view, 
partially defined operations are common and easy to deal with, whereas 
infinite objects pose problems (e.g., it is not clear how one should define 
complexity of computations on infinite words); for that reason we will not 
use Thompson's original definition. Higman \cite{Hig74} defined $V$ 
by total functions on larger algebras, that contain $\{a,b\}^*$. 
But Higman's actions are uniquely determined by the partial action on 
$\{a,b\}^*$ (and are in fact the same as the partial actions in \cite{ESc}, 
which themselves are the same as ours, up to terminology). 

\medskip

\noindent {\bf Remark on other Thompson groups:} \ 
Richard Thompson defined subgroups of $V$ that are of great interest. One of
them, denoted $\hat{\mathbb P}$ in \cite{Th}, ${\mathfrak P}'$ in 
\cite{McKTh}, and $F$ in \cite{CFP}, is defined as follows: \  
$F$ is the subgroup of $V$ consisting of all maximal right-ideal 
isomorphisms of $\{a,b\}^*$ that {\it preserve the dictionary order}.

The dictionary order on $\{a,b\}^*$ (with $a < b$) is a very classical concept
and is defined as follows: 
For $x, y \in \{a,b\}^*$ we have $x \leq_{\rm dict} y$ iff $x$ is a prefix of
$y$ or there exists $p \in \{a,b\}^*$ such that $x \in pa\{a,b\}^*$ and 
$y \in pb\{a,b\}^*$ (so $p$ is the maximal common prefix of $x$ and $y$).
One can easily verify that this is a total order, compatible with 
concatenation on the right for non-prefix comparable words (i.e., 
$x \leq_{\rm dict} y$ and $x$, $y$ prefix incomparable implies 
$xw \leq_{\rm dict} yw$). In the tree picture of $\{a,b\}^*$, if $x$, $y$ are
prefix incomparable then we have: 
$x <_{\rm dict} y$ iff $x$ is in a tree branch that is more to the left than 
the tree branch containing $y$. 
We will say that a map $\varphi : \{a,b\}^* \to \{a,b\}^*$ 
preserves the dictionary order iff the following holds for all 
$x_1, x_2 \in \{a,b\}^*$: if $x_1 \leq_{\rm dict} x_2$ and if $x$, $y$ are 
prefix incomparable then $\varphi(x_1) \leq_{\rm dict} \varphi(x_2)$. 
One can prove easily that if an isomorphism between two essential right 
ideals preserves the dictionary order then its maximum extension also 
preserves the dictionary order. Thompson proved that $F$ is a finitely 
presented group, whose commutator is simple and of finite index.


\section{Word length in the Thompson group $V$}

Let $\Delta$ be a finite set of generators of $V$. 

\begin{defn} \  
 For every element $g \in V$, the word length of $g$ 
(over the generating set $\Delta$) is the length of a shortest word 
$\in (\Delta^{\pm 1})^*$ that represents $g$. 

The word length of $g$ is denoted by $|g|_{\Delta}$.  
\end{defn}
It is easy to prove that if $\Omega$ is another finite set of generators of 
$V$ then \ $|g|_{\Omega} \leq C_{1 2} \cdot |g|_{\Delta}$ \ and \ 
$|g|_{\Delta} \leq C_{2 1} \cdot |g|_{\Omega}$, where 
$C_{1 2}, C_{2 1} > 0$ depend on $\Delta$ and $\Omega$, but not on $g$.

Since the elements of $V$ are functions, there is another size 
measure for elements of $V$. In the following we will denote the 
(finite of infinite) {\it cardinality} of a set $X$ by $|X|$.

\begin{defn} \  
For a right-ideal isomorphism $\varphi: P_1A^* \to P_2A^*$, where
$P_1$ and $P_2$ are finite maximal prefix codes, the restriction
$P_1 \to P_2$ of $\varphi$ is called the {\bf table} of $\varphi$. 
(Recall that this restriction is a bijection.)

We define $\|\varphi \|$ to be $|P_1| \ (= |P_2|)$; we call this the 
{\bf table size} of $\varphi$.

For an element $g \in V$, the table size $\|g\|$ of $g$ is defined 
to be the table size of the maximally extended right-ideal isomorphism that 
represents $g$.
\end{defn}
The following lemmas will be useful when we study the table size of 
right-ideal isomorphisms.

\begin{lem}
\label{idealInter} \  
Let $P, Q, R \subseteq A^*$ be such that $PA^* \cap QA^* = RA^*$, and 
$R$ is a prefix code. Then $R \subseteq P \cup Q$.

As a consequence, the intersection of two finitely generated right
ideals is finitely generated. 
\end{lem}
{\bf Proof.} \ The Lemma has a simple and intuitive interpretation in terms
of prefix trees. We'll give a formal proof, which is almost as simple.

For any $r \in R$ there exist $p \in P, q \in Q$ and  $v, w \in A^*$
such that $r = pv = qw$. Hence $p$ and $q$ are prefix-comparable.
Let us assume $p \geq_{\rm pref} q = px$, for some $x \in A^*$ (the other 
case is handled the same way). 
Hence $q = px \in PA^* \cap QA^* = RA^*$, and $q$ is a 
prefix of $r = qw$. Since $R$ is a prefix code, $r = q$, hence $r \in Q$. 
 \ \ \ $\Box$

\begin{lem}
\label{idealIncl} \  
Let $P, Q \subseteq A^*$ be maximal prefix codes such that 
$PA^* \subseteq QA^*$.  Then $|Q| \leq |P|$.
\end{lem}
{\bf Proof.} \ For every $p \in P$ there is $q \in Q$ such that 
$p \leq_{\rm pref} q$. 
In fact, this correspondence $p \mapsto q$ is a function; indeed, if there 
are $q_1, q_2 \in Q$ such that $p \leq_{\rm pref} q_1$ and 
$p \leq_{\rm pref} q_2$ then $p = q_1x_1 = q_2x_2$ (for some 
$x_1, x_2 \in A^*$), hence $q_1$ and $q_2$ are prefix-comparable. 
This implies $q_1 = q_2$ since $Q$ is a prefix code. 

Moreover, this map is surjective. Indeed, let $q \in Q$. If $q \in PA^*$ 
then there exists a prefix of $q$ in $P$, hence an inverse. If 
$q \not\in PA^*$ then the inverse exists by Lemma \ref{ideals} (4). 

Since there is a surjective function $P \to Q$, the result follows. \ \ \ 
$\Box$

\begin{pro}
\label{composLength} \  
For any right-ideal isomorphisms $\varphi_2$ and $\varphi_1$ between
essential right ideals of $A^*$:

\smallskip

 \ \ \ \ \ \ \ \ \ \ \ \ \ \ \ \ \ \ \ \ \ \ \ \ \ \ 
 $\|{\sf max} \, \varphi_1\| \ \leq \ \|\varphi_1\|$

\noindent and

 \ \ \ \ \ \ \ \ \ \ \ \ \  \ \ \ \ \ \ \ \ \ \ \ \ \ 
 $\|\varphi_2 \cdot \varphi_1\| \ \leq \ \|\varphi_2 \circ \varphi_1\| \ 
\leq \ \|\varphi_2\| + \|\varphi_1\|$
\end{pro}
{\bf Proof.} \ The fact that \ 
$\|{\sf max} \, \varphi_1\| \ \leq \ \|\varphi_1\|$ \ 
follows directly from Lemma \ref{idealIncl}.  

\smallskip 

\noindent Let $\varphi_1: P_1A^* \to P_1'A^*$ and 
$\varphi_2: P_2A^* \to P_2'A^*$, where $P_1, P_1', P_2, P_2'$ are maximal 
prefix codes, and  $\|\varphi_1\| = |P_1| = |P_1'|$, 
$\|\varphi_2\| = |P_2| = |P_2'|$.
Then the domain of the functional composite $\varphi_2 \circ \varphi_1$ is 
a right ideal $RA^*$ where $R$ is a maximal prefix code; hence, 
$\|\varphi_2 \circ \varphi_1\| = |R|$. Moreover, \ 
$RA^* =  \varphi_1^{-1}(P_1'A^* \cap  P_2A^*)$.
By Lemma \ref{idealInter}, $P_1'A^* \cap  P_2A^* = SA^*$ for some maximal 
prefix code $S$ such that $|S| \leq |P_1'| + |P_2|$ 
$ = \|\varphi_1\| + \|\varphi_2\|$. Since $RA^*$ is the domain of 
$\varphi_2 \circ \varphi_1$, $\varphi_1$ is defined everywhere on $RA^*$; 
and since $SA^* \subseteq P_1'A^*$ (which is the domain of $\varphi_1^{-1}$),
$\varphi_1^{-1}$ is defined everywhere on $SA^*$. Thus, $\varphi_1^{-1}$ is
a bijection from $SA^*$ onto $RA^*$, hence by Lemma \ref{bijPrefCode}, 
$|S| = |R|$ $(= \|\varphi_2 \circ \varphi_1\|)$. 
 \ It follows that $|R| \leq |P_1| + |P_2|$. \ \ \ $\Box$

\begin{lem}
\label{lengthCode} \  
Let $P \subseteq A^*$ be a finite maximal prefix code. Then any word in $P$
has length at most \ $\frac{|P|-1}{|A|-1}.$ \ 
In particular, when the alphabet has 2 letters, the length is at most $|P|-1$.
\end{lem}
{\bf Proof.} \ Consider the prefix tree of $P$, which has $|P|$ leaves. Let 
the number of non-leaves be $N$. Then $N = \frac{|P|-1}{|A|-1}$, as can easily
be shown by induction on $N$. The length of a word in $P$ is equal to the 
length of the path from the root to the leaf labeled by this word; such a path
has length at most $N$. \ \ \ $\Box$

As a consequence of Lemmas \ref{composLength} and \ref{lengthCode}, we have:

\begin{cor}
\label{lengthVStable} \  
Let $\Delta$ be a fixed finite generating set of $V$. If 
$\varphi \in V - \{ {\bf 1} \}$ is described by a word 
of length $n$ over $\Delta^{\pm 1}$, then the table size satisfies \   
$\| \varphi\| \leq C_{\Delta} \, n$  \ (where \  
$C_{\Delta} = {\rm max}\{ \| \delta \| : \delta \in \Delta\}$). 

Similarly, the length of the longest word in the table of $\varphi$ 
is $\leq C_{\Delta} \, n$.
\end{cor}
We will prove that the two size measures (namely word length and table size) on 
elements of $V$ are closely related. This similar to what happens in 
the symmetric groups ${\mathfrak S}_k$, concerning the relation between $k$ 
and the word length of permutations (over a bounded number of generators, with 
bound independent of $k$). 

\begin{thm}
\label{sizes} \  
The table size and word size of an element $g \in V$ are 
related as follows:

\smallskip

\noindent {\bf (1)} \ 
There are $c_{_{\Delta}}, c'_{_{\Delta}} > 0$ (depending on the choice 
of $\Delta$) such that for all $g \in V - \{ {\bf 1} \}$:  
$$c'_{_{\Delta}} \, \|g\| \ \leq \ |g|_{\Delta} \ \leq \ c_{_{\Delta}} \, \|g\| \cdot \log_{_2} \|g\|.$$
{\bf (2)} \   For {\em almost all} $g \in V$,  
$$|g|_{\Delta} \ > \ \|g\| \cdot \log_{_{2 \, |\Delta|}} \|g\|.$$
``Almost all'' means here that in the set \   
$\{g \in V : \|g\| = n\}$, 
the subset that does {\em not} satisfy the above inequality has a proportion 
that tends to 0 exponentially fast as $n \to \infty$.
\end{thm}

\noindent Inequality (2) shows that up to big-O, the function $x \cdot \log x$ 
is the best possible upper bound in terms of $\|g\|$.

However, although inequality (2) holds for ``almost all'' $g \in V$, 
it also fails to hold for infinitely many $g \in V$; for example,
we will see in Proposition \ref{wLengthTfin<} below that for all 
$g \in F$, \ $|g|_{\Delta}  <  c \, \|g\|$.

\medskip

\noindent {\bf Proof of (2).} \ 
The proof is a counting argument. The number of maximal prefix codes 
of cardinality $n$ over the alphabet $\{a,b\}$ is the Catalan number $C_{n-1}$ 
(see the beginning of our Appendix A1). 

If we count only elements of $V$ with domain code 
$\{a^ib : i = 0, 1, \ldots , n-2\} \cup \{a^{n-1}\}$, and an arbitrary fixed 
range code of cardinality $n$, the number of elements of $V$ 
obtained is  $\geq  n (n-2) (n-2)!$. Note that the number is not $n!$ because 
we want to make sure to count only maximal right-ideal isomorphisms; therefore,
if we choose to map $a^{n-1}$ to some word $ua$, we cannot map $a^{n-2}b$ to 
$ub$, respectively $ua$; thus only $n-2$ choices exist for the image of 
$a^{n-2}b$.
Asymptotically, however, $n (n-2) (n-2)!$ and $n!$ are equivalent. Hence, 
the number of elements $g \in V$ with $\|g\| = n$ is at least \  
$C_{n-1} \, n! = \frac{1}{2n-1} \, \frac{(2(n-1))!}{(n-1)!(n-1)!} \, n!$. 
By Stirling's formula, this is equal to \  
$2^{-1/2} (4/e)^{n-1} (n-1)^{n-1} \cdot (1 + \varepsilon(n))$ \ (where \  
lim$_{n \to \infty} \, \varepsilon(n) = 0$).  

For any $\ell$, the number of words over $\Delta^{\pm 1}$ of length 
$\leq \ell$ is $\leq  \delta^{\ell +1}$, where $\delta = 2 \, |\Delta|$. 
Hence, in $V$ we have:  
The ratio of the number of elements that have word length \  
$\ell \leq n \cdot \log_{\delta} n$,  
over the number of elements that have table size $n$, is less than \   

\smallskip

$n^n \cdot 2^{1/2} (e/4)^{n-1} (n-1)^{1-n} \cdot (1 + \varepsilon_1(n)) \ = \ $
$2^{1/2} (e/4)^{n-1} n \, (\frac{n}{n-1})^{n-1} \cdot (1 + \varepsilon_1(n))$

\smallskip

$ = \ 2^{1/2}  (e/4)^{n-1} n \, e \cdot (1 + \varepsilon_2(n))$.

\smallskip

\noindent This ratio tends to $0$ exponentially fast as $n \to \infty$ 
(since $e < 4$).  \ \ \ $\Box$

\bigskip

\noindent {\bf Proof of (1).} \ We proved the first inequality of (1) already  
in Corollary \ref{lengthVStable}.
The proof of the second inequality consists of three steps. In summary:  
 
\smallskip

\noindent 
(1.1) We give a canonical factorization of any element of $V$, as
a right-ideal automorphism and two elements of $F$.  
 
\smallskip

\noindent 
(1.2) We show that all elements of $F$ have linearly bounded word length.  

\smallskip

\noindent
(1.3) We prove that the word length of right-ideal automorphisms is \ 
$\leq c \, \|g\| \cdot \log \|g\|$. 

\bigskip

\noindent {\bf (1.1) Canonical factorization.} \ 
At the end of Section 1 we already mentioned the subgroup $F$, which
consists of the elements of $V$ that preserve the dictionary order 
of $\{a,b\}^*$. 

A {\bf right-ideal automorphism} of a finitely generated essential right ideal 
$P\{a,b\}^*$ (where $P$ is a maximal prefix code) has a table whose domain 
code and range code are the same (namely $P$); the table gives a permutation 
of $P$.  

Contrary to a first impression, the set of right-ideal automorphisms (of all
essential right ideals) 
is not a group, and it is not closed under restriction nor under extension. 

\begin{pro}
\label{canFact} \  
For every $n \geq 1$ let us fix one maximal prefix code $S_n$ of cardinality
$n$. Then for every $g \in V - \{ {\bf 1} \}$ there exist unique 
elements $\alpha_g, \beta_g, \pi_g \in V$ such that \ 

\smallskip

 \ \ $g = \beta_g \pi_g \alpha_g$, \ 

\smallskip

 \ \ $\alpha_g$ and $\beta_g$ belong to $F$,  

\smallskip

 \ \ $\pi_g$ is an automorphism whose table is a permutation of $S_{\|g\|}$.  

\smallskip

\noindent Moreover, $\|\alpha_g\|, \|\beta_g\|, \|\pi_g\| \leq \|g\|$.
\end{pro}   
{\bf Proof.} \ Let $\| g \| = n$. 
Consider a maximal right-ideal isomorphism that represents $g$,
and let $\varphi: P_1 \to P_2$ be its table, where $P_1, P_2$ are maximal
prefix codes, $|P_1| = |P_2| = n$. 

We define $\alpha_g$ by mapping $P_1$ in an order-preserving way bijectively
onto $S_n$; in other words, the table of $\alpha_g$ is obtained by taking
the elements of $P_1$ in increasing dictionary order as the domain, and by 
taking the elements of $S_n$ in increasing dictionary order as the 
corresponding range.
Similarly, $\beta_g$ is defined by mapping $S_n$ in an order-preserving 
way bijectively onto  $P_2$. Uniqueness of $\alpha_g$ and $\beta_g$ 
follows from the fact that once the domain and image codes of elements
of $F$ are specified, the elements are uniquely determined.
Finally, we simply let \ $\pi_g =  \beta_g^{-1} \varphi \alpha_g^{-1}$; 
hence, $\pi_g$ is also uniquely determined.
 \ \ \ $\Box$

\bigskip 

The idea of a factorization of the above type appears in proofs of Thompson's
\cite{Th}, where he uses the family of maximal prefix codes \ 
$S_n = \{a^ib : 0 \leq i \leq n-2\} \cup \{a^{n-1}\}$. This family has the 
following nice property (which Thompson does not mention or use, however): 
If one only considers automorphisms whose domain (and range) code is of the
form $S_n = \{a^ib : 0 \leq i \leq n-2\} \cup \{a^{n-1}\}$ (for 
$n \in \mathbb{N}$), then this particular set of automorphisms (as $n$ ranges
over all integers $> 1$) is a subgroup of 
$V$, and this set of automorphisms is closed under 
extension and restriction. The only other family of maximal prefix codes with 
this property is $\{b^ia : 0 \leq i \leq n-2\} \cup \{b^{n-1}\}$. 
The reason is that for those two families of prefix codes there is only one 
place in a domain (and range) code, where extension of some automorphisms is 
possible (namely just above the deepest point in the prefix tree, where 
$\{a^{n-1}, a^{n-2}b\}$ will be replaced by $\{a^{n-2}\}$); restriction 
(subject to the constraint that the code of the restriction should belong to 
this particular family of prefix codes) is also only possible at 
one place (again, at the deepest point in the prefix tree). 

The above two families of prefix codes have a disadvantage for us: The depth 
of the prefix tree of such an $S_n$ is $n-1$; as a consequence we would get 
$\|g\|^2$ in our theorem, instead of $\|g\| \cdot \log \|g\|$. So, we will use 
the following family of codes, that have logarithmic depth.  

As the domain code and range code of the table of $\pi_g$ we pick a set 
$S_n \subset \{a,b\}^{k-1} \cup \{a,b\}^k$, where $|S_n| = n = \|g\|$, 
and $k = \lceil \log_{_2} n \rceil$ (equivalently, $2^{k-1} < n \leq 2^k$).
When $n = 2^k$, $S_n = \{a,b\}^k$. When $n < 2^k$, we choose 
$S_n$ as in the figure below, representing the prefix tree of the 
maximal prefix code $S_n$. The horizontal lines indicate the leaves of 
the prefix tree (i.e., the elements of $S_n$); the higher one of the two 
lines pictures the vertices at depth $k-1$ (i.e., the elements of $S_n$ 
in $\{a,b\}^{k-1}$, of which there are $2^k-n$), the lower horizontal line 
pictures the vertices at depth $k$ (i.e., the elements of $S_n$ in 
$\{a,b\}^k$, of which there are $2n-2^k$). 

\bigskip

\unitlength=0.90mm \special{em:linewidth 0.4pt}
\linethickness{0.4pt}

\begin{picture}(70,55)
\put(70,50){\line(-1,-1){44}}
\put(70,50){\line(1,-1){39}}

\put(26,6){\line(1,0){50}}  
\put(26,5.8){\line(1,0){50}}  

\put(76,11){\line(1,0){33}}  
\put(76,10.8){\line(1,0){33}}  

\put(77,0){\makebox(0,0)[cc]{$S_n$}}
\end{picture}

\bigskip

{\sf Fig.~2: The maximal prefix code $S_n$}

\bigskip 

\bigskip

\noindent {\bf (1.2) Word length in} $F${\bf .} \ The subgroup 
$F$ is generated by the following two elements of $F$ (Thompson \cite{Th0}, 
\cite{Th}):

\[ \sigma \ = \ \left[ \begin{array}{ccc}
         \ a^2 \ & \ ab \  & \ b \\
           \ a \ & \ ba  \ & \ b^2 
\end{array}        \right]  \]

\[ \theta \ = \ \left[ \begin{array}{cccc}
         \ a \ & \ ba^2 \ & \ bab \  & \ b^2 \\
         \ a \ & \ ba   \ & \ b^2a \ & \ b^3
\end{array}        \right]  \]

\medskip

\noindent In \cite{CFP}, our $\sigma$ is called $A^{-1}$ (we changed the 
notation because we use $A$ to denote the alphabet); in Section 2 we saw 
that $\sigma$ can be interpreted as the shift operator on $\mathbb Z$. 
In \cite{CFP}, our $\theta$ is called $B^{-1}$. 

\begin{pro}
\label{wLengthTfin<} \ 
For every $g \in F$ we have \ \ 
$|g|_{ \{\sigma,\theta\}}  \ < \  4 \, \|g\|$.  

In words: \ $F$ has linearly bounded word length. 
\end{pro}
{\bf Proof.} \  Let $g \in F$ be represented by a table 
$\varphi: R \to S$, \ let $n+1 = \|g\| = \|\varphi\|$, and let \ 
$X_i = \sigma^{i-1} \theta^{-1} \sigma^{-i+1}$ for all $i \geq 1$, and
$X_0 = \sigma^{-1}$. 
Cannon, Floyd and Parry \cite{CFP} (Theorem 2.5, page 223) prove that 

\medskip 

$g \ = \ $
$X_0^{b_0}X_1^{b_1} \ldots X_n^{b_n}X_n^{-a_n} \ldots X_1^{-a_1}X_0^{-a_0}$

\medskip

\noindent where $b_{\ell}$ $(0 \leq \ell \leq n)$ is the length of the longest 
path in the prefix tree of $S$, subject to the following conditions:  

$\bullet$ the path consists only of left-edges,   
 
$\bullet$ the start vertex of the path is leaf number $\ell$ (the leaves are
numbered from 0 through $n$),  

$\bullet$ the end vertex of the path does not have a label in 
 $b^* \ (\subset \{a,b\}^*)$.  

\smallskip  

\noindent Similarly, one defines $a_{\ell}$ $(0 \leq \ell \leq n)$ for the 
prefix tree of $R$.  (One observes that for the right-most leaf 
$a_n = b_n =0$, so the above expression could be simplified; but that doesn't 
matter.)  By replacing each $X_i$ we obtain  
 
\medskip 

$g \ = \ $
$\sigma^{-b_0} \theta^{-b_1} \sigma \theta^{-b_2} \sigma \theta^{-b_3}$
$ \ldots \ \sigma \theta^{-b_{n-1}} \sigma \ \theta^{-b_n + a_n} \ $
$\sigma \theta^{a_{n-1}} \sigma \ \ldots \ $
$\theta^{a_3} \sigma \theta^{a_2} \sigma \theta^{a_1} \sigma^{a_0}$

\medskip

\noindent hence, \ \ $|g|_{ \{\sigma,\theta\} } \ < \ $
$\sum_{i=0}^n b_i  + n + \sum_{i=0}^n a_i  + n$. 

\medskip

\noindent By the definition of $b_{\ell}$,  \  
$\sum_{i=0}^n b_i$ \ is less than the number of left-edges in the  
prefix tree of $S$. In a prefix tree over the alphabet $\{a,b\}$ there is an 
equal number of left-edges and right-edges (since every vertex has 0 or 2 
children). Moreover, the total number $e$ of edges in a prefix tree with 
$n+1$ leaves satisfies \ $e = 2n$ \ (since such a tree has $n$ interior
vertices, and each interior vertex corresponds to two edges, and vice versa). 
Therefore, \ 
$\sum_{i=0}^n b_i < n$ \ (and similarly, \ $\sum_{i=0}^n a_i < n$). 
The result follows. \ \ \ $\Box$

\bigskip

\noindent {\bf (1.3) Word length of right-ideal automorphisms.} \ 
Let us prove the claimed bound on the word length of all right-ideal 
automorphisms.
Let $\pi: S \to S$ be the table of any automorphism, where $S$ is any finite 
maximal prefix code, and where $\pi$ is a permutation of $S$ (so here $S$ is 
not necessarily of the form $S_n$). It is 
well known that every permutation of a finite set $\{1, 2, \ldots, n\}$ can be 
expressed as the composition of $\leq 3n$ transpositions of the form $(1 | i)$
with $i \in \{1, 2, \ldots, n\}$. Indeed, we can first take disjoint cycles; 
and for a cycle we have \ $(x_1|x_2|x_3| \ldots |x_{r-1}|x_r) = $
$(x_1|x_r)(x_1|x_{r-1}) \ldots (x_1|x_3)(x_1|x_2)$; finally, for a 
transposition, $(i|j) = (1|i)(1|j)(1|i)$. Recall that all our functions and
permutations are applied on the left of the argument.

We will write the automorphism $\pi$ as the product of $\leq 3 \, \|\pi\|$ 
transpositions of the form $(a^k|w)$
with $a^k, w \in S$; in particular, $w \not\in a^*\ (\subset \{a,b\}^*)$, and
$k \geq 1$. 

The transposition $(a^k|w)$ is defined as follows.
Let $j$ $(0 \leq j < k)$ be such that $a^jb$ is  a prefix of $w$; such a $j$
exists (and is unique) since $w \not\in a^*$ and since $w$ is not 
prefix-comparable with $a^k$. So $w$ can be written as $w = a^jbv$, for 
some $v \in \{a,b\}^*$.  
Then $(a^k|w)$ is defined by the table 

\smallskip

\[ (a^k|w) \ = \ \left[ \begin{array}{cccc}
 \ a^k \ & \ w  \  & \ a^ib \       & \ a^jb p\ell \      \\
 \   w \ & \ a^k \ & \ a^ib \       & \ a^jb p\ell \      
\end{array}        \right]  \]

\vspace{-.1in}

\makebox[3.4in]{} $_{0\leq i < k}$ \ \ $_{p > v}$

\vspace{-.05in}

\makebox[3.5in]{} $_{i \neq j}$ \ \ \ $_{p\ell \not\geq v}$

\bigskip

\noindent Here the range of $i$ is $0 \leq i \leq k-1$ and $i \neq j$. The word
$p$ ranges over all strict prefixes of $v$; the notation $x > y$ is short for 
$x >_{\rm pref} y$ and means that $x$ is a strict prefix of $y$ (as defined in
the Introduction); $\ell \in \{a,b\}$ is such 
that $p\ell$ is {\it not} a prefix of $v$. For every strict prefix $p$ there 
will be exactly one letter $\ell$ such that $p\ell$ is not a prefix of $v$. 

In our canonical factorization $g = \beta_g \pi_g \alpha_g$, the automorphism
$\pi_g$ has a table which is a permutation of the maximal prefix code $S_n$,
where $n = \|g\|$. We saw that all words in $S_n$ have length \ 
$\leq \lceil \log_{_2} n \rceil$.  Also, $\| \pi_g \| = |S_n| = n$.
So, when we factor $\pi_g$ as $\leq 3 \, n$ transpositions of the form
$(a^k|w)$, the parameters $k$ and $w$ satisfy \ \ 
$k, |w| \ \leq \ \lceil \log_{_2} n \rceil$.

\smallskip

Therefore the next Lemma will complete the proof of Theorem \ref{sizes}.

\begin{lem} 
\label{wLengthTaut} \ 
Every transposition $(a^k|w)$ has word length  $\leq c \cdot (k + |w|)$
 $\leq 2c \, \lceil \log_{_2} n \rceil$ 
over some finite set of generators of $V$ (for some constant $c>0$).
\end{lem}
{\bf Proof.} \ We will use the following {\bf generators} of $V$: 
 \ $\sigma$ and $\theta$ \ 
(the generators of $F$ used before), and 

\[ \gamma_1 \ = \ \left[ \begin{array}{cccc}
 \ a^2 \ & \ aba \ & \ ab^2 \ & \ b  \      \\
 \ a^2 \ & \ ba  \ & \ ab   \ & \ b^2 \
\end{array}        \right]  \]

\[ \gamma_2 \ = \ \left[ \begin{array}{cccc}
 \ a^2 \ & \ aba \ & \ ab^2 \ & \ b  \      \\
 \ a^2 \ & \ ab  \ & \ ba   \ & \ b^2 \
\end{array}        \right]  \]

\[ \delta  \ = \ \left[ \begin{array}{cccc}
 \ a^3 \ & \ ab \ & \ a^2b \ & \ b  \      \\
 \ a^2 \ & \ ab  \ & \ ba   \ & \ b^2 \
\end{array}        \right]  \]

\[ (a^2|ab)  \ = \ \left[ \begin{array}{ccc}
 \ a^2 \ & \ ab  \ & \ b  \      \\
 \ ab  \ & \ a^2 \ & \ b  \   
\end{array}        \right]  \]

\[ (a^2|a ba)  \ = \ \left[ \begin{array}{cccc}
 \ a^2 \ & \ aba \ & \ ab^2 \ & \ b \  \\
 \ aba \ & \ a^2 \ & \ ab^2 \ & \ b \
\end{array}        \right]  \]

\[ (ab | b) \ = \ \left[ \begin{array}{ccc}
 \ a^2 \ & \ ab \ & \ b  \      \\
 \ a^2 \ & \  b \ & \ ab \ 
\end{array}        \right]  \]

\[ (a| b) \ = \ \left[ \begin{array}{cc}
 \ a \ & \ b  \      \\
 \ b \ & \ a \
\end{array}        \right]  \]

\noindent {\sc Case 1}: \ The transposition $(a^k|w)$ is such that 
$w \in b \{a,b\}^*$  \ (i.e., $w$ starts with $b$).  

\smallskip

Recall that \ $k, |w| \in [1, \lceil \log_{_2} n \rceil ]$.
We will eliminate Case 1 by showing:
{\it If $(a^k|w)$ is conjugated by at most $k + |w|$ generators, a 
transposition of the form $(a^K|az)$ is obtained, where \ 
$K + |az| \leq k + |w|+1$.  
} 
Indeed, we have:

\medskip 

(1.1)  \ \ \ \ \ \ \ \ \ \ \ 
$\sigma^{-1} \cdot (a^k|b^h av) \cdot \sigma =  (a^{k+1}|b^{h-1} av)$ \ \ 
 for all $h \geq 2$, $v \in  \{a,b\}^*$, $k \geq 1$; 

\medskip

(1.2) \ \ \ \ \ \ \ \ \ \ \
$\sigma^{-1} \cdot (a^k|bav) \cdot \sigma =  (a^{k+1}|abv)$ \ \
 for all $v \in  \{a,b\}^*$, $k \geq 1$;

\medskip

(1.3) \ \ \ \ \ \ \ \ \ \ \ 
$\sigma^{-1} \cdot (a^k|b^h) \cdot \sigma = (a^{k+1}|b^{h-1})$ \ \ 
 for all $h \geq 2$, $k \geq 1$;

\medskip

(1.4) \ \ \ \ \ \ \ \ \ \ \
$(ab|b) \cdot (a^k|b) \cdot (ab|b) = (a^k| ab)$ \ \ 
 for all $k \geq 2$. 

\smallskip

 \ \ \ \ \ \ \    \ \ \ \ \ \ \ \ \ \ \
(When $k = 1$, $(a|b)$ is a generator.)

\medskip 

\noindent For the detailed verification of (1.1), (1.2), (1.3), and (1.4),
see Appendix A2.

If the transposition $(a^k|w)$ is in  situation (1.1), i.e., $w = b^h av$ with
$h \geq 2$, then after $h-1$ conjugations we are in situation (1.2).
In situation (1.2), one conjugation gets us out of Case 1.
If the transposition $(a^k|w)$ is in  situation (1.3), i.e., $w = b^h$ with
$h \geq 2$, then after $h-1$ conjugations we are in situation (1.4).
In situation (1.4), one conjugation gets us out of Case 1.
Thus we have eliminated Case 1.

\medskip

\noindent {\sc Case 2}: \ The transposition $(a^k|w)$ is such that
$w \in a \{a,b\}^*$  \ (i.e., $w$ starts with $a$). 

Since $a^k$ and $w$ belong 
to the same maximal prefix code they are not prefix-comparable; therefore, 
$w$ is of the form $w = a^j b^h v$ where $k > j \geq 1$ (hence $k \geq 2$), 
$h \geq 1$, $v \in \{a,b\}^* - b\{a,b\}^*$ \ ($v$ does not start with $b$).

\smallskip

\noindent We have: \  

\medskip

\noindent (2.1) \ \ \ \ \ \ \ \ \ \ \ 
If \ $(k>) \ j \geq 2$, $h \geq1$, \ then \ \ \ \ 
$\sigma \cdot (a^k|a^j b^h v) \cdot \sigma^{-1} = (a^{k-1} |a^{j-1} b^h v)$.    

\medskip

\noindent (2.2) \ \ \ \ \ \ \ \ \ \ \   
If \ $j =1$, $k \geq 3$, $h \geq1$, \ then \ \ \ \
$\delta \cdot (a^k|a b^h v) \cdot \delta^{-1} = (a^{k-1} |ab^h v)$.  

\medskip

\noindent (2.3) \ \ \ \ \ \ \ \ \ \ \   
If \ $j =1$, $k =2$, $h \geq 2$, \ then \ \ \ \
$\gamma_1 \cdot (a^2|a b^h v) \cdot \gamma_1^{-1} = (a^2|a b^{h-1} v)$.

\medskip

\noindent (2.4) \ \ \ \ \ \ \ \ \ \ \   
If \ $j =1$, $k =2$, $h =1$, $|v| \geq 2$ (so $v = au$ for some 
$u \in  \{a,b\} \, \{a,b\}^*$), \ then \ \ \ \ 

\smallskip

 \ \ \ \ \ \ \ \ \ \ \ \ \ \ \ \ \ \ \ 
   $\gamma_2 \cdot (a^2|a b au) \cdot \gamma_2^{-1} = (a^2|a b u)$. 

\medskip

\noindent (2.5) \ \ \ \ \ \ \ \ \ \ \   
If \ $j =1$, $k =2$, $h =1$, $|v| \leq 1 $, \ then \ 
$(a^2|a ba)$ \ and \ $(a^2|a b)$ \ are generators.
 
\medskip

\noindent The detailed verification of (2.1) -- (2.4) appears in Appendix A2.

If the transposition $(a^k|w)$ is in  situation (2.1), then after $j-1$ steps
we reach situation (2.2); in each step, both $k$ and $|w|$ decrease. 
If the transposition $(a^k|w)$ is in situation (2.2), then after $k-2$ steps 
(during which $k$ keeps decreasing while $w$ remains unchanged) we reach 
situation (2.3).
If the transposition $(a^2|w)$ is in situations (2.3) or (2.4), we remain in 
situations (2.3) or (2.4) as long as $w$ is of the form 
$a b^h v$ (with $h \geq 2$) or $a b au$; at each step, $w$ becomes shorter.
Eventually we reach case (2.5).  

This completes the proof of Lemma \ref{wLengthTaut}, and hence of Theorem
\ref{sizes}.   \ \ \ $\Box$


\section{The word problem of $V$ }

It is fairly obvious from the representation of $V$ by partial 
functions that the word problem is decidable. In this section we show that 
the word problem of $V$ has low complexity.

Since $V$ is finitely generated, we can consider the word problem 
of this group; let $\Delta$ be a finite set of generators of $V$. 
As an input for the word problem we consider a string over the alphabet
$\Delta^{\pm 1}$; we ask whether the product of the elements in 
that string is equal to {\bf 1} (the identity).  

We can conclude immediately from Corollary \ref{lengthVStable} that 
{\em the word problem of} $V$ {\em is in} {\bf P} (deterministic 
polynomial time).
Indeed, a product of $n$ generators has table size $\leq C_{\Delta} n$.
Moreover, every entry in the table is a word of length $\leq C_{\Delta} n$, 
so the total memory space occupied by the table is $O(n^2)$.
To solve the word problem we just compose the partial functions, and we check 
whether the product is a ``subidentity'' (i.e., $x$ is mapped to $x$ for 
every $x$ in the finite maximal prefix code of the product).  
Let us mention also the interesting fact that $V$ has polynomially
bounded isoperimetric function (Guba \cite{Guba}).

One can prove much more detailed and stronger complexity results. 
For this we will use the {\it parallel complexity classes} AC$^0$ and AC$^1$, 
which are subclasses of {\bf P}. In short, AC$^k$ (for $k \in \mathbb{N}$) 
consists of those problems whose output can be computed by
acyclic boolean circuits (where the gates have unbounded finite fan-in), of 
depth $O((\log n)^k)$, and of size polynomial in $n$ (where $n$ is the input 
length).  See \cite{Handb}, \cite{Vollmer} and \cite{Wegener} for details.

It is well known that the prefix relation between two words can be checked
by an AC$^0$ circuit; moreover, given $s$ and $t =sz \in A^*$, an AC$^0$ 
circuit can output $z$. 

\medskip

Let us define the problems precisely.

\medskip

\noindent 
(1) {\it Composition problem} \\  
{\sc Input}: Two isomorphisms between essential right ideals 
$\varphi: P_1A^* \to Q_1A^*$, and $\psi: P_2A^* \to Q_2A^*$, where 
$P_1, Q_1, P_2, Q_2$ are maximal prefix codes;
$\varphi$ is described by its finite table \ 
$\{ (p, \varphi(p)) : p \in P_1 \} \subseteq P_1 \times Q_1$, \  
and similarly for $\psi$. \\ 
{\sc Output}: The composite $\psi \circ \varphi$, described by a 
finite table.

\medskip

\noindent
(2) {\it Word problem of} $V$ \\
{\sc Input}: A sequence $(\psi_1, \ldots, \psi_n)$ (where $n \in \mathbb{N}$
is also variable) of isomorphisms between essential right ideals \ 
$\psi_i : P_iA^* \to Q_iA^*$, where $P_i, Q_i$ are maximal prefix codes
($i = 1, \ldots, n$);  
each $\psi_i$ is described by its finite table \ 
$\{ (p, \psi_i(p)) : p \in P_i \} \subseteq P_i \times Q_i$. \\  
{\sc Output}: ``Yes'' if \   
$\psi_1 \cdot \ldots \cdot \psi_n = {\bf 1}$ \ (i.e., the product in 
$V$ is the identity map on $A^*$); ``no'' otherwise.

Note that this definition of the word problem is consistent with the usual
definition. Since $V$ is finitely generated the usual definition of 
the word problem can be applied to any finite set of generators. 
By Corollary \ref{lengthVStable} the length of the description of 
$(\psi_1, \ldots, \psi_n)$ by finite tables is linearly bounded by the length 
of $(\psi_1, \ldots, \psi_n)$ as words over a fixed finite set of generators
of $V$. Moreover, the description of the input of the word problem
by tables is not significantly more compact than a description by a word over
generators, by Theorem \ref{sizes}.

\begin{thm} \  
{\bf (1)} \ The composition problem (for two isomorphisms between essential 
  right ideals) is in {\rm AC}$^0$. \\  
{\bf (2)} \ The word problem of $V$ is in {\rm AC}$^1$. 
\end{thm}
{\bf Proof.} \ {\bf (1)} \  Let \ $\varphi = \{(x_i,y_i) : i=1, \ldots, n\}$
and \ $\psi = \{(u_j, v_j) : j = i, \ldots, m\}$.  
For an AC$^0$ algorithm, we consider all $nm$ indiviual composites 
$(x_i, y_i) \circ (u_j, v_j)$ 
in parallel. Each composite $(x_i, y_i) \circ (u_j, v_j)$ can be computed by 
a constant-depth circuit $C_{i,j}$ by checking the prefix-relation between 
$x_i$ and $v_j$; if $x_i$ and $v_j$ are not prefix-comparable, the circuit 
$C_{i,j}$ outputs 0; if $x_i = v_j z \leq_{\rm pref} v_j$ (for some 
$z \in A^*$), the circuit $C_{i,j}$ outputs $(u_j z, y_i)$; 
if $x_i \geq_{\rm pref} v_j = x_i z$ (for some $z \in A^*$), the circuit 
$C_{i,j}$ outputs $(u_j, y_i z)$. The size of $C_{i,j}$ is \ 
$O({\sf max}\{|x_i|, |v_j|\})$, hence $O(n)$ by Lemma \ref{lengthCode}. 

Finally, the overall circuit consists of all the $C_{i,j}$ 
($1 \leq i \leq n, 1 \leq j \leq m$), and an additional layer which masks 
outputs that are 0, and just outputs the set of non-zero outputs of the 
circuits $C_{i,j}$. The overall size of the circuit is \

\smallskip
   
 \ \ \ \ \  $O(n \cdot m \ \cdot $
${\sf max}\{|x_i|,|y_i|,|v_j|,|u_j|: 1\leq i\leq n,1\leq j\leq m\})$, 

\smallskip

\noindent  which corresponds to big-O of the cube of the input length. 

\medskip

\noindent
{\bf (2)} \ We are given a sequence $(\psi_1, \ldots, \psi_n)$ of isomorphisms
between essential right ideals; let $N$ be the total length of this input
(when each $\psi_i$ is described by a finite table).
We compose the $\psi_i$ two-by-two in parallel ($\psi_1$ with $\psi_2$, 
$\psi_3$ with $\psi_4$, etc.); then we compose the
resulting $\lceil n/2 \rceil$ functions two-by-two, then the resulting
$\lceil n/4 \rceil$, etc.; we obtain a ``composition tree'' of depth
$\lceil \log_{_2} n \rceil$, with $n-1$ nodes.
By Proposition \ref{composLength} and Lemma \ref{lengthCode}, the intermediate
and the final composites have finite tables of space $O(N^2)$ each.
By {\bf (1)} above, each node of the composition tree can be implemented by a 
constant-depth circuit of size $O((N^2)^3)$. 
So the total size of the circuit is $O(n N^6)$, and the depth is $O(\log n)$.

Finally, once we have a finite table for the composite
$\psi_1 \cdot \ldots \cdot \psi_n$
we can check easily (in constant depth) whether this function is a partial
identity.     \ \ \ $\Box$

\bigskip

Since the word problem of $V$ is in AC$^1$ it follows 
(see \cite{Handb}) that it is also in DSpace($(\log n)^2$). 

\bigskip

We saw already that the word problem of $V$ is in {\bf P}.  The next 
proposition gives more precise bounds on the time complexity; the bounds for 
deterministic and nondeterministic time are still rather crude. 
Here, ``coNTime'', ``NTime'', and ``DTime'' refer to the usual time complexity 
classes (co-nondeterministic, nondeterministic, and deterministic
time complexity, respectively; see \cite{Handb}). 

\begin{pro} 
\label{wrdprobTime} \   
The word problem of $V$ is in {\rm coNTime}$(n)$, in {\rm NTime}$(n^2)$, 
and in {\rm DTime}$(n^3)$. 
\end{pro}
{\bf Proof.} \ Let $\Delta$ be a finite generating set of $V$, let
$w \in (\Delta^{\pm 1})^*$, let $n = |w|$, and let $\varphi_w$ be the 
right-ideal isomorphism obtained by composing the generators as they appear 
in $w$, without taking maximum extensions. Obviously, $w = 1$ in $V$ 
iff $\varphi_w$ is a partial identity.  
Recall (Corollary \ref{lengthVStable}) that   
$\| \varphi_w \| \leq C_{\Delta} \, n$, and that the length of the longest 
word in the table of $\varphi_w$ is $\leq C_{\Delta} \, n$. 
Hence the sum of the lengths of the words in the table of $\varphi_w$ is 
$\leq C_{\Delta}^2 \, n^2$. 

\smallskip

\noindent {\bf (1)} The word problem of $V$ is in coNTime$(n)$ iff 
the negation of the word problem is in NTime$(n)$. To prove the latter we 
consider a nondeterministic two-tape Turing machine, one tape containing the
input $w$; on the other tape the machine guesses a word $z \in \{a,b\}^*$
such that $z$ belongs to the domain code of $\varphi_w$ (so 
$|z| \leq C_{\Delta} \, n$), and $\varphi_w(z) \neq z$. 
Such a word $z$ exists iff $w \neq 1$ in $V$. Moreover, 
$|z| \leq C_{\Delta} \, n$, hence it takes only linear time to write $z$ on
the second tape. 

On tape 1 the machine also makes a copy of $z$.

Next, the machine reads $w$ from right to left, and applies the generators to 
the word $(\in \{a,b\}^*)$ on tape 2 (initially that word is $z$). 
Each application of a generator changes a prefix of bounded length of the word 
on tape 2. (This bound is the length of the longest word in the table of 
any generator $\in \Delta$.) Hence, applying one generator takes only a 
bounded amount of time, and the total time to apply all the generators in $w$
is linear in $|w|$; hence the time 
is $\leq cn$ for some constant $c$.  The final content of tape 2 will be 
$\varphi_w(z)$. 
  
Finally, the machine compares the original copy of $z$ (saved on tape 1) with 
$\varphi_w(z)$; this takes linear time. It accepts if these two words are 
different.

\smallskip

\noindent {\bf (2)} \ To prove that the word problem of $V$ is in 
NTime$(n^2)$ we consider a nondeterministic Turing machine which guesses the
table of a partial identity map, and then verifies that the guess is indeed 
the table of $\varphi_w$. The machine accepts $w$ if the verification 
succeeds. Such an accepting computation exists iff $w = 1$ in $V$. 

More precisely, the Turing machine guesses words $x_i \in \{a,b\}^*$, 
each of length $|x_i| < \|\varphi_w \| \ (\leq C_{\Delta} \, n)$, in such a way 
that $\{x_i : 1 \leq i \leq \|\varphi_w \| \}$ is a maximal prefix code. 
This is done by means on a {\it nondeterministic depth-first search} in the 
tree of $\{a,b\}^*$. This search starts at the root of the tree and proceeds
like an ordinary depth-first search, except that we guess where the leaves of
the prefix tree of the code $\{x_i : 1 \leq i \leq \|\varphi_w \| \}$ are; 
each time we reach such a leaf we write down the corresponding word $x_i$. 
Since the set of the guessed $x_i$ consists of the leaves of a subtree of 
$\{a,b\}^*$, it is a prefix code; since the depth-first search will not end 
until it has visited all the leaves of a tree, this prefix code is maximal. 
The running time of the search is proportional to the number of leaves, which 
is $\|\varphi_w \|$ (in an accepting computation). 
The time to write down all the $x_i$ is 
$< \|\varphi_w \|^2$ (since $|x_i| < \|\varphi_w \|$ in case of acceptance).
Thus, we can guess the table of a partial identity of table size 
$\|\varphi_w\|$ in time $O(n^2)$.
 
Next, the Turing machine verifies, deterministically, that 
$\varphi_w(x_i) = x_i$ for each $i$. 
This is done in the same way as in (1), and takes linear time for each $x_i$. 
Hence the total time of a successful verification is $O(n^2)$.   

\smallskip

\noindent {\bf (3)} \ Let us first look at the deterministic time complexity 
of composition (without maximum extension) of two elements 
$\varphi, \psi \in V$, described by their tables. Let  
$\varphi = \{(s_i,t_i) : i=1, \ldots, \ell\}$, hence by Lemma 
\ref{lengthCode},  \ $|t_i|, |s_i| < \|\varphi \| = \ell$. 
Let  $\psi = \{(u_j,v_j) : j=1, \ldots, m\}$,  hence \      
     $|v_j|, |u_j| < \|\psi \| = m$. 
Suppose the two tables are written on two different tapes of a deterministic 
Turing machine. The machine can then write the set \ 
$\{(u_j,v_j) \circ (s_i,t_i) : 1 \leq i \leq \ell, \ 1 \leq j \leq m\}$ \ on a 
third tape and reduce each term $(u_j,v_j) \circ (s_i,t_i)$ to either {\bf 0}
(in which case it is removed from the set), or to a new term of the from 
$(x_k, y_k)$. This takes time \ 
$\leq c_1 \, \|\varphi \|^2 \, \|\psi \|^2$ for some constant $c_1$.  
  
Finally, we can compute the table of \ 
$\varphi_w = d_n \ldots d_k \ldots d_2 d_1$ (with $d_k \in \Delta^{\pm 1}$) 
by composing from right to left. Recall that $C_{\Delta}$ be the maximum table 
size of any generator in $\Delta$.
Inductively, assume the table of $d_k \ldots d_2 d_1$ is $\psi$, with \ 
$\| \psi \| \leq C_{\Delta} k$. 
Then the table for $d_{k+1} \psi$ can be computed in time  
$\leq c_1 \, \|d_{k+1}\|^2 \, \|\psi \|^2$
$ \leq c_1 \, C_{\Delta}^2 \, C_{\Delta}^2 \, k^2$.
So the total time to compute the table of $w$ is \ 
$\leq c_1 \, C_{\Delta}^4 \, \sum_{i=1}^n k^2 = O(n^3)$.
             \ \ \ $\Box$


\section{Generalized word problems of $V$ }

Generalized word problems of a group $G$ ask about membership of elements of
$G$ in specified subgroups; the word problem is the special case obtained by
taking the trivial subgroup $\{1\}$. Let $G$ be a group with finite 
generating set $A$. If $S \subseteq G$ and $w \in (A^{\pm 1})^*$, we 
write $w \in_G S$ iff the element of $G$ represented by $w$ is in $S$.  
If $x_1, \ldots, x_n \in (A^{\pm 1})^*$ then 
$\langle x_1, \ldots, x_n \rangle_G$ denotes the subgroup of $G$ 
generated by the elements of $G$ represented by the words 
$\{x_1, \ldots, x_n \}$.
If $S \subseteq G$ then $\langle S \rangle_G$ denotes the subgroup of $G$ 
generated by $S$. 

\smallskip

By definition, the {\bf generalized word problem} of a fixed group $G$ with a
fixed finite generating set $A$, is specified as follows: \\   
{\sc Input}: A finite set \ $X \subset (A^{\pm 1})^*$, and an additional 
``test word'' $y \in (A^{\pm 1})^*$. \\  
{\sc Question}: \ $y \in_{_G} \langle X \rangle_G$ \   
(i.e., does the element of $G$ represented by the word $y$ belong to the   
subgroup of $G$ generated by the elements represented by the words in $X$)? 

\smallskip

\noindent We consider two {\bf special forms of the generalized word problem}.

\smallskip

\noindent {\bf (1)} The generalized word problem of a fixed group $G$, with a 
fixed finite generating set $A$ and a {\it fixed subgroup generating set} \  
$X = \{x_1, \ldots, x_k \} \subset (A^{\pm 1})^*$, is specified as follows: \\
{\sc Input}: A word $y \in (A^{\pm 1})^*$. \\  
{\sc Question}: \ $y \in_{_G} \langle x_1, \ldots, x_k \rangle_G$ ?

\smallskip

\noindent {\bf (2)} The generalized word problem of a fixed group $G$ with a 
fixed finite generating set $A$ and a {\it fixed ``test word''} 
$y \in (A^{\pm 1})^*$, is specified as follows: \\
{\sc Input}: A finite set \ $X \subset (A^{\pm 1})^*$. \\  
{\sc Question}: \ $y \in_{_G} \langle X \rangle_G$ ?

\bigskip

We will use the next few facts to obtain undecidability results for the 
generalized word problem of $V$ and its special versions.

Graham Higman proved that  {\it Thompson's group $V$ contains a 
subgroup isomorphic to  ${\rm FG}_2$ (the free group on two generators).} 
A proof appears in the Appendix of \cite{EScConjug}.

Thompson \cite{Th} proved that if $G$ is a subgroup of $V$ then 
$V$ also contains a subgroup isomorphic to the direct product 
$G \times G$. Hence, $V$ contains all finite direct powers of $G$.
Thompson \cite{Th} (and also Higman \cite{Hig74}) mention that this can be 
generalized to any finite direct product and to countably infinite 
direct sums. The {\it direct sum} of a family of groups 
$(G_i : i \in I)$ is defined to be the subgroup of the direct product 
$\prod_{i\in I} G_i$ consisting of the sequences that have only a finite 
number of non-identity components.  

As a consequence of these results of Thompson and Higman we obtain: \  
{\em $V$ contains a subgroup isomorphic to} ${\rm FG}_2 \times {\rm FG}_2$.

\medskip

We will also need to talk about ``uniform word problems'', and related 
problems for Turing machines. 

\smallskip

\noindent $\bullet$ \ The {\it uniform word problem} for groups is specified 
as follows. \\
{\sc Input:}  A finite presentation $\langle B, R \rangle$ of a group, and 
a word $w \in (B^{\pm 1})^*$ \ (where letters in the alphabet $B$ are 
encoded over some fixed finite alphabet). \\
{\sc Question:} $w =_{\langle B, R \rangle} 1$\ ? \  I.e., is  the element of
the group $\langle B, R \rangle$ represented by $w$ the identity?

\smallskip

\noindent $\bullet$ \ 
The {\it uniform word problem for a fixed word} $w_0 \in (B_0^{\pm 1})^*$
(where $B_0$ is a fixed finite alphabet), is specified as follows. \\
{\sc Input:}  A finite presentation $\langle B, R \rangle$, with
$B_0 \subseteq B$ \ (where letters in the alphabet $B$ are encoded over
some fixed finite alphabet). \\
{\sc Question:} $w_0 =_{\langle B, R \rangle} 1$ ?

\smallskip

\noindent $\bullet$ \  
The {\it acceptance problem} is specified as follows. \\
{\sc Input:}  A Turing machine $M$ and a word $w$ over the tape alphabet of
$M$ (where letters and states of $M$ are encoded over some fixed finite
alphabet). \\
{\sc Question:} Does $M$ accept $w$ ?

\smallskip

\noindent $\bullet$ \  
The {\it fixed-word acceptance problem} for a fixed word $w_0 \in B_0^*$
(where $B_0$ is a fixed finite alphabet), is specified as follows. \\
{\sc Input:}  A Turing machine $M$ whose tape alphabet contains $B_0$
(such that letters and states of $M$ are encoded over some fixed finite
alphabet). \\
{\sc Question:} Does $M$ accept $w_0$ ?

\begin{lem}
\label{unifWP} \  
There are infinitely many words $w_0$
such that the uniform word problem with fixed word $w_0$ is undecidable.
\end{lem}
{\bf Proof.} \ By Rice's theorem, the acceptance problem for any fixed word 
$w_0$ (and variable Turing machine $M$) is undecidable. 
By Boone's Lemma (see Lemma 12.7 in \cite{Rotman}), there is a computable 
function with linear time complexity which maps $(w, M)$
(where $M$ is a Turing machine and $w$ is a word over the tape alphabet
of $M$) to a finite presentation $G_M = \langle B_M, R_M \rangle$ and a
word $\rho(w) \in (B_M^{\pm 1})^*$ such that:
 \ $M$ accepts $w$ \ iff \ $\rho(w) =_{G_M} 1$.

Hence, for any fixed word $\rho(w_0)$, the fixed-word uniform word problem 
of $G_M$ is undecidable.
 \ \ \ $\Box$

\begin{pro} \ 
For $V$ (with any fixed finite generating set $\Delta$), the 
generalized word problem is undecidable.

The first special form of the generalized word problem is undecidable for some
fixed subgroups.

The second special form of the generalized word problem is undecidable for 
some fixed non-empty test words.
\end{pro}

\noindent {\bf Proof.} \ We saw that $V$ contains 
${\rm FG}_2 \times {\rm FG}_2$. Therefore we can apply Mikhailova's theorem 
\cite{Mikh58} (see Theorem IV.4.3 in \cite{LyndonSchupp}), which states that 
the generalized word problem, as well as its first special form (for some 
fixed finite subgroup generator sets $X$), are undecidable for 
${\rm FG}_2 \times {\rm FG}_2$.

 Moreover, it follows from the proof of Mikhailova's theorem 
(see Lemma IV.4.2 in \cite{LyndonSchupp}) that the second special form of the 
generalized word problem (for some fixed test words) of 
${\rm FG}_2 \times {\rm FG}_2$ is also undecidable. Indeed, this proof gives a 
reduction of the uniform word problem ``is $w =_H 1$ ?'' 
(where $H$ is a variable finitely presented 
group, and $w$ is a variable word) to the generalized word problem of 
${\rm FG}_2 \times {\rm FG}_2$,  with some 
subgroup generating set $L_H$ (where $L_H$ is finite). 
We saw in Lemma \ref{unifWP} that the uniform word problem for a fixed 
word (but variable finite presentations) is undecidable.

Hence, the generalized word problem and its two special forms are undecidable 
for any group containing ${\rm FG}_2 \times {\rm FG}_2$, and in particular 
for $V$.  \ \ \ $\Box$

\bigskip

\noindent {\bf Other decision problems}

\smallskip

Graham Higman (Theorem 9.3 in \cite{Hig74}) proved that the
{\bf conjugacy problem} and the {\bf order problem} of $V$ are decidable.
See also \cite{EScSurvey}, \cite{EScConjug}.

{\bf Open problem:} Is the {\bf generation problem} of $V$ decidable? 
(The generation problem is specified as follows: The input is a finite set 
$\Gamma$ of elements of $V$, given by their tables; the question is whether
$\Gamma$ generates all of $V$.)


\section{Distortion}

We just saw that in some (in fact, infinitely many) cases, the first special 
form of the generalized word problem of $V$ (for a fixed subgroup)
is undecidable. This leads to the question: What can the complexity of
this problem be when the problem is decidable? 
 
In relation to the first special form of the generalized word problem 
(for a fixed subgroup), the ``Cayley graph distortion function'' will play 
the role of inherent (group-theoretic) complexity. We will see that it is 
closely connected with nondeterministic time complexity. 
Let $G$ be a fixed group with fixed finite generating set $A$, and let
$X \subset (A^{\pm 1})^*$ be a fixed finite set, generating a subgroup 
$H = \langle X\rangle_G$. If $y \in (A^{\pm 1})^*$
is such that $y \in_{_G} \langle X\rangle_G$,  then $y$ has two lengths,
namely, one over the alphabet $A^{\pm 1}$ and one over $X^{\pm 1}$. 
How the two lengths are related is an important question, first addressed
by Gromov \cite{Gromov}.

\begin{defn} \ 
{\bf (1)} For a group $G$ with finite generating set $A$, and an element 
$g \in G$, the {\em $A$-length} of $g$, denoted by $|g|_A$ is the minimum 
length of any word $y \in (A^{\pm 1})^*$ that represents $g$. 

In the Cayley graph $\Gamma(G,A)$ of $G$ with generating set $A$,  
$|g|_A$ is the distance from the root (i.e., the identity element of $G$)
of the element $g$.

\smallskip

{\bf (2)} Let $X$ be a finite subset of $(A^{\pm 1})^*$ and let 
$H = \langle X\rangle_G$ be the subgroup of $G$ generated by $X$.
If $y \in_{_G} \langle X\rangle_G$, we define the {\em $X$-length}, 
denoted by $|y|_X$, to be the length of a shortest sequence 
$(x_{i_1}, x_{i_2}, \ldots, x_{i_{\ell}})$ of elements of $X^{\pm 1}$ 
such that \ $y =_{_G} x_{i_1} x_{i_2} \ldots x_{i_{\ell}}$. 

In the Cayley graph $\Gamma(G,A)$, the subgroup $H$ with generating set $X$ 
generates a path-subgraph, whose vertex set is $H$, and 
whose edges correspond to $X$-labeled paths between vertices in 
$H$. This path-subgraph is the image of a path-embedding of the Cayley graph 
$\Gamma(H, X)$ into $\Gamma(G,A)$. 
(If $X \subseteq A$, then $\Gamma(H, X)$ is a subgraph of $\Gamma(G,A)$.)  
\end{defn}
The Cayley graph distortion function (we'll call it simply ``distortion
function'') describes the relation between the two lengths of an element 
$y$ of the subgroup $\langle X\rangle_G$ of $G$. 

\begin{defn} \  
Let us fix a group $G$ with finite generating set $A$, and a finite set
$X \subset (A^{\pm 1})^*$.  A non-decreasing function 
$f: \mathbb{N} \to \mathbb{N}$ is a {\bf distortion function} 
of the subgroup $H = \langle X\rangle_G$ in $G$ iff for all elements 
$g \in \langle X\rangle_G$ we have: \ \ \ $|g|_X \leq f(|g|_A)$.
\end{defn}
Note that according to our definition, distortion functions are 
non-decreasing (i.e., if $m \leq n$ then $f(m) \leq f(n)$). 
The minimum distortion function for a given $G$, $A$, and $X$ (as above)
is called ``the'' distortion function.

\medskip

Two functions $f_1$ and $f_2: \mathbb{N} \to \mathbb{N}$ are said to be 
{\em equivalent} (or more precisely, {\it linearly equivalent}) iff 
there exist positive constants $c_0$, $C_{1 2}$, $c_{1 2}$, $C_{2 1}$, 
$c_{2 1}$ such that for all $n \geq c_0$: \ \  
$f_1(n) \leq C_{1 2} \, f_2(c_{1 2} \, n)$ \ and \  
             $f_2(n) \leq C_{2 1} \, f_1(c_{2 1} \, n)$.  
 \ In big-O notation, $f_1$ and $f_2$ are linearly equivalent \ iff \  
$f_1(n) = O(f_2(O(n)))$ and $f_2 = O(f_1(O(n)))$.

\medskip

The concept of distortion was formally introduced by Gromov \cite{Gromov}. 
However, Gromov's definition used an additional factor $n$; so, Gromov's 
distortion is constant when ours is linear. We follow Ol'shanskii and Sapir
\cite{OlSap}, whose theorem connecting distortion to Dehn functions (stated
below) makes Gromov's version of the definition look less well motivated.

It is easy to see the following: When one changes the generating set $A$ of 
$G$ to another finite generating set (of the same group $G$), and one changes
the set $X$ that generates the subgroup $H = \langle X\rangle_G$ to another
finite set (that generates the same subgroup $H$), then the minimum
distortion does not change (up to linear equivalence); see e.g., \cite{Farb}. 
Hence for finitely generated groups the distortion depends only on the 
groups $G$ and $H$ $(\subseteq G)$ (up to linear equivalence). 

The distortion cannot be less than linear, since (up to finite change of
generators) one can increase $A$ so that $X \subseteq A$. When the distortion 
is linear one also says that ``there is no distortion'', or that $H$ is 
``isometrically embedded'' in $G$.  

In \cite{Bi}, \cite{Olsh} and \cite{BORS} it was proved that the
{\it Higman embedding theorem} for semigroups, respectively groups, can be 
strengthened in such a way that  the distortion is linear. 
Another important result about distortion functions is the following theorem 
of Ol'shanskii and Sapir \cite{OlSap}: {\it The set of distortion functions 
of finitely generated subgroups of ${\rm FG}_2 \times {\rm FG}_2$ coincides 
(up to linear equivalence) with the set of all Dehn functions of finitely
presented groups.}  
Guba and Sapir \cite{GubaSapir} proved that for any integer $d \geq 2$ the 
Thompson group $F$ (in the notation of \cite{CFP}) has a subgroup with 
distortion $\geq n^d$ (up to linear equivalence). 

\begin{thm}
\label{distTfin}
 \ The set of distortion functions of the finitely generated subgroups of 
the Thompson group $V$ contains (up to linear equivalence) the set 
of all Dehn functions of finitely presented groups.
\end{thm}
{\bf Proof.} \ This is an immediate consequence of the Ol'shanskii-Sapir 
theorem and Theorem \ref{linearDist}. \ \ \ $\Box$

\begin{thm}  
\label{linearDist} \  
{\bf (1)} \ Thompson's group $V$ contains a subgroup isomorphic to the
free group ${\rm FG}_2$ with linear distortion.  \\     
{\bf (2)} \ The group $V$ also contains a subgroup isomorphic to
${\rm FG}_2 \times {\rm FG}_2$ with linear distortion.
\end{thm}
{\bf Proof of (1).} \ We start out from Graham Higman's result (mentioned 
earlier), that Thompson's group $V$ contains a subgroup isomorphic to 
${\rm FG}_2$. In the Appendix of \cite{EScConjug} it is shown that the 
following two elements of $V$ generate a free group:

\[ \alpha \ = \ \left[ \begin{array}{cccccc}
         \ a \ & \ b^3 \  & \ b^2ab \ & \ b^2a^3 \ & \ b^2a^2b \ & \ ba \\
      \ b^4a \ & \ b^3a \ & \ b^2a \  & \  a \     & \  ba  \    & \  b^5
\end{array}        \right]  \]
\[ \beta \ = \ \left[ \begin{array}{cccccc}
         \ b  \ & \ ab \  & \ a^2b^2 \ & \ a^2ba^2 \ & \ a^2bab \ & \ a^3 \\
      \ a^3ba \ & \ a^3b^2 \ & \ a^2b \ & \ b \      & \ ab  \ & \ a^4
\end{array}        \right]  \]

\noindent Our main goal now is to show that the free subgroup generated by 
$\{ \alpha, \beta \}$ has {\it linear distortion} in $V$.

Let $\mu$ be any element of $\langle \alpha, \beta \rangle$, with 
$\mu \neq {\bf 1}$. If we write $\mu$ in such a way that there are no 
cancellations then $\mu$ has one of the following expressions:  

\smallskip

{\sc Case} $\alpha \alpha$: \ \ \ \ \
$\mu = \ \ \ \ \ \alpha^{h_m} \beta^{k_{m-1}} \ldots 
        \beta^{k_i} \alpha^{h_i} \ldots \alpha^{h_2} \beta^{k_1} \alpha^{h_1}$, 

\smallskip 

{\sc Case} $\beta \alpha$: \ \ \ \ \ 
$\mu = \beta^{k_m} \alpha^{h_m} \beta^{k_{m-1}} \ldots 
      \beta^{k_i} \alpha^{h_i} \ldots \alpha^{h_2} \beta^{k_1} \alpha^{h_1}$, 

\smallskip 

{\sc Case} $\alpha \beta$: \ \ \ \ \  
$\mu = \ \ \ \ \ \alpha^{h_m} \beta^{k_{m-1}} \ldots 
      \beta^{k_i} \alpha^{h_i} \ldots \alpha^{h_2} \beta^{k_1}$,

\smallskip 

{\sc Case} $\beta \beta$: \ \ \ \ \
$\mu = \beta^{k_m} \alpha^{h_m} \beta^{k_{m-1}} \ldots 
       \beta^{k_i} \alpha^{h_i} \ldots \alpha^{h_2} \beta^{k_1}$.

\smallskip

\noindent with \   
$h_m, k_{m-1}, h_{m-1}, \ldots, k_2, h_2, k_1 \in \mathbb{Z} -\{0\}$, \  
$k_m \in \mathbb{Z} -\{0\}$ in cases $\beta \alpha$ and $\beta \beta$, \ 
$k_m = 0$ in cases $\alpha \alpha$ and $\alpha \beta$, \ 
$h_1 \in \mathbb{Z} -\{0\}$ in cases $\alpha \alpha$ and $\alpha \beta$, \
$h_1 = 0$ in cases $\alpha \beta$ and $\beta \beta$. 

Since $\langle \alpha, \beta \rangle$ is a free group and $\mu$ is reduced,  
the $\{\alpha, \beta\}$-length of $\mu$ is \ 
$|\mu|_{\{\alpha, \beta\}} = |k_m| + |h_m| + \ldots + |k_1| + |h_1|$.

In order to prove that the subgroup $\langle \alpha, \beta \rangle$ of 
$V$ has linear distortion, we need to show that the minimum length 
of $\mu$ over some fixed finite generating set $\Delta$ of $V$ is \ 
$|\mu|_{\Delta} \geq c \cdot |\mu|_{\{\alpha, \beta\}}$ \ (for some constant 
$c >0$ depending only on the chosen set of generators $\Delta$ of $V$).

We will use the following method. We will show that for any $\mu$ there is 
a word $y \in \{a,b\}^*$ such that 

\medskip

\noindent {\sc Case} $(* \alpha)$. \ 
If $\mu$ is in cases $\alpha \alpha$ or $\beta \alpha$, then $y$ satisfies:

\smallskip

$\mu(a) = ya$, 

\smallskip

$\mu(ba) = yb$, 

\smallskip

$|y| > |\mu|_{\{\alpha, \beta\}}$. 

\medskip

\noindent {\sc Case} $(* \beta)$. \ 
If $\mu$ is in cases $\alpha \beta$ or $\beta \beta$, then $y$ satisfies:

\medskip

$\mu(b) = ya$

\smallskip

$\mu(ab) = yb$, 

\smallskip

$|y| > |\mu|_{\{\alpha, \beta\}}$.

\medskip

\noindent
The existence of $y$, as above, implies  that the table of the {\em maximum
extension} of $\mu$ contains an entry of length $> |\mu|_{\{\alpha, \beta\}}$.
Indeed, $ya$ and $yb$ are such entries (in the range code of 
${\sf max} \, \mu$); by Lemma \ref{ExtCritFin}, $ya$ and $yb$ remain in the 
table when $\mu$ is maximally extended, since $ya$ and $yb$ are the images of 
$a$, resp.\ $ba$ (or $b$ resp.\ $ab$). 

Now, by Lemma \ref{lengthCode}, the table size of the maximum extension 
of $\mu$ satisfies \   $\|{\sf max} \, \mu\| > |\mu|_{\{\alpha, \beta\}}$.
Also, by Proposition \ref{composLength}, \ 
$\|{\sf max} \, \mu\| \leq \| \mu\|$. 

By Corollary \ref{lengthVStable}, \  
$(\|{\sf max} \, \mu\| \leq)$ $\| \mu\| \leq C_{\Delta} \, |\mu|_{\Delta}$, 
therefore by the above, \ 
$|\mu|_{\{\alpha, \beta\}} < C_{\Delta} \, |\mu|_{\Delta}$. 

So the proof of Theorem \ref{linearDist} will be complete once we prove the 
following two claims, which show that the appropriate $y$ exists.

\medskip

\noindent {\bf Claim } $(*\alpha)$.  \ {\it 
Let $\mu$ be as above, according to cases $\alpha \alpha$ or $\beta \alpha$. 
Then,

\smallskip

$\mu(a) =  (w_m v_m^{|k_m|-1}) \ t_m u_m^{|h_m|-1} w_{m-1} \ldots $
             $ w_i v_i^{|k_i|-1} t_i u_i^{|h_i|-1} \ldots $
                         $ w_1 v_1^{|k_1|-1} t_1 \ u_1^{|h_1|-1} w_0$

\medskip

\noindent where: 

\smallskip

$ w_0 = \left\{ \begin{array}{ll}
                ba  & \mbox{if \ \ $h_1> 0$} \\
                a^2 & \mbox{if \ \ $h_1 < 0$} \end{array}  \right. $
 \ \ \ (Note that $w_0$ ends in $a$.) 

\medskip

\noindent For \ $m-1 \geq i \geq 1$: \ \ \ \ 
$ w_i = \left\{ \begin{array}{ll}
                ba^3  & \mbox{if \ \ $h_{i+1} > 0$, $k_i > 0$} \\ 
                ba^2b & \mbox{if \ \ $h_{i+1} > 0$, $k_i < 0$} \\ 
                a^4   & \mbox{if \ \ $h_{i+1} < 0$, $k_i > 0$} \\
                a^3b  & \mbox{if \ \ $h_{i+1} < 0$, $k_i < 0$}
                \end{array}  \right. $   

\medskip

\noindent In case $\beta \alpha$: \ \ \ \ 
$ w_m = \left\{ \begin{array}{ll}
                a^3  & \mbox{if \ \ $k_m > 0$} \\ 
                a^2b & \mbox{if \ \ $k_m < 0$} \end{array}  \right. $ 

\medskip

\noindent For \ $m-1 \geq i \geq 1$,  or $i = m$ in case $\beta \alpha$: 
 \ \ $ t_i = \left\{ \begin{array}{ll}
                bab^2  & \mbox{if \ \ $h_i > 0$, $k_i > 0$} \\
                a^2b^2 & \mbox{if \ \ $h_i > 0$, $k_i < 0$} \\
                baba   & \mbox{if \ \ $h_i < 0$, $k_i > 0$} \\
                a^2ba  & \mbox{if \ \ $h_i < 0$, $k_i < 0$}
                \end{array}  \right. $ 

\medskip

\noindent In case $\alpha \alpha$: \ \ \ \ 
The factor $(w_m v_m^{|k_m|-1})$ is absent, and \ 
 $t_m = \left\{ \begin{array}{ll}
                b^3  & \mbox{if \ \ $h_m > 0$} \\
                b^2a & \mbox{if \ \ $h_m < 0$} \end{array}  \right. $ 

\medskip
 
\noindent For $m \geq i \geq 1$: \ \ \ \ $ u_i = \left\{ \begin{array}{ll}
                a & \mbox{if \ \ $h_i > 0$} \\
                b & \mbox{if \ \ $h_i < 0$}
                \end{array}  \right. $   \ \ \ \ \  \ \ \  
$ v_i = \left\{ \begin{array}{ll}
                a & \mbox{if \ \ $k_i > 0$} \\
                b & \mbox{if \ \ $k_i < 0$}
                \end{array}  \right. $ 

\medskip
 
\noindent We also have: 
 
\smallskip 

$\mu(ba) = (w_m v_m^{|k_m|-1}) \ t_m u_m^{|h_m|-1} w_{m-1} \ldots $
             $ w_i v_i^{|k_i|-1} t_i u_i^{|h_i|-1} \ldots $
                         $ w_1 v_1^{|k_1|-1} t_1 u_1^{|h_1|-1} W_0$

\smallskip

\noindent where:

\medskip
 
$W_0 = \left\{ \begin{array}{ll}
                b^2  & \mbox{if \ \ $h_1> 0$} \\
                ab & \mbox{if \ \ $h_1 < 0$} \end{array}  \right. $
 \ \ \ (Note that $W_0$ ends in $b$.)

\medskip

\noindent and all other $w_i$, $v_i$, $t_i$, and $u_i$ are the same as 
for $\mu(a)$.
} 

\medskip

\noindent  {\bf Proof of Claim} $(* \alpha)$. \ The proof goes by induction 
on the number of exponents $k_m, h_m, \ldots, k_1, h_1$. 
For the details, see Appendix A3.

\bigskip
 
\noindent {\bf Claim } $(*\beta)$.  \ {\it
Let $\mu$ be as above, according to cases $\alpha \beta$ or $\beta \beta$.
Then,

\smallskip

$\mu(b) =  (w_m v_m^{|k_m|-1}) \ t_m u_m^{|h_m|-1} w_{m-1} \ldots $
             $ w_i v_i^{|k_i|-1} t_i u_i^{|h_i|-1} \ldots $
                         $ w_1 v_1^{|k_1|-1} t_1$

\smallskip

\noindent where:

\smallskip

\noindent $t_1 = \left\{ \begin{array}{ll}
                  ba  & \mbox{if \ \ $k_1> 0$} \\
                  a^2 & \mbox{if \ \ $k_1 < 0$} \end{array}  \right. $
 \ \ \ (Note that $t_1$ ends in $a$.)

\medskip

\noindent and all other $w_i$, $v_i$, $t_i$, and $u_i$ are the same as
for $\mu(a)$ in Claim $(*\alpha)$.

\medskip

\noindent We also have:

\smallskip

$\mu(ab) = (w_m v_m^{|k_m|-1}) \ t_m u_m^{|h_m|-1} w_{m-1} \ldots $
             $ w_i v_i^{|k_i|-1} t_i u_i^{|h_i|-1} \ldots $
                         $ w_1 v_1^{|k_1|-1} T_1$

\smallskip

\noindent where:

\medskip

$T_1 = \left\{ \begin{array}{ll}
                  b^2  & \mbox{if \ \ $k_1> 0$} \\
                  ab   & \mbox{if \ \ $k_1 < 0$} \end{array}  \right. $ 
 \ \ \ (Note that $T_1$ ends in $b$.)

\medskip

\noindent and all other $w_i$, $v_i$, $t_i$, and $u_i$ are the same as
for $\mu(b)$.
} 

\medskip

\noindent {\sc Proof of Claim} $(* \beta)$. \ The proof goes by induction on
the number of exponents $k_m, h_m, \ldots, k_1$ in $\mu$.  
For details, see Appendix A3.

This completes the proof of part {\bf (1)} of Theorem \ref{linearDist}.

\bigskip

{\bf Proof of (2).} \  
For each element $\varphi \in V$ we consider two elements
$\varphi_a, \varphi_b \in V$ defined as follows:

\medskip

 \ \ \ $\varphi_a(ax) =  a \, \varphi(x)$,  \ \ for all $x$
          in the domain of $\varphi$ ;

\smallskip

 \ \ \ $\varphi_a(b) = b$.

\medskip

 \ \ \ $\varphi_b(a) = a$ ;

\smallskip

 \ \ \ $\varphi_b(bx) = b \, \varphi(x)$, \ \ for all $x$
            in the domain of $\varphi$.

\medskip

\noindent Then from the generators $\alpha, \beta$ of FG$_2$ seen in 
the proof of {\bf (1)}, we obtain a set 
$\{\alpha_a, \beta_a, \alpha_b, \beta_b\}$ which generates a subgroup of
$V$ isomorphic to ${\rm FG}_2 \times {\rm FG}_2$. The proof that the 
distortion of $\langle \alpha_a, \beta_a, \alpha_b, \beta_b \rangle$ in $V$
is linear is very similar to the proof for the distortion of 
$\langle \alpha, \beta \rangle$ in $V$. 

Any element of $\langle \alpha_a, \beta_a, \alpha_b, \beta_b \rangle$ can 
be put in the form $\varphi_a \psi_b$, for some $\varphi, \psi \in V$. 
Looking at cases
$(* \alpha)$ and $(* \beta)$ for both $\varphi$ and $\psi$ (in the proof of
part {\bf (1)}), we obtain four cases. In each case we find that there are 
words $y, z \in \{a,b\}^*$  such that

\medskip

\noindent {\sc Case} $(* \alpha, * \alpha)$ \ \ \ \ \  

$\varphi_a \psi_b(aa) = a \, \varphi(a) = aya$, \ \   
$\varphi_a \psi_b(aba) = a \, \varphi(ba) = ayb$, \ \ and 
$|y| > |\varphi|_{\{\alpha,\beta\}}$;

$\varphi_a \psi_b(ba) = b \, \psi(a) = bza$, \ \  
$\varphi_a \psi_b(bba) = b \, \psi(ba) = bzb$, \ \ and 
$|z| > |\psi|_{\{\alpha,\beta\}}$. 

\medskip

\noindent {\sc Case} $(* \alpha, * \beta)$ \ \ \ \ \  

$\varphi_a \psi_b(aa) = a \, \varphi(a) = aya$, \ \   
$\varphi_a \psi_b(aba) = a \, \varphi(ba) = ayb$, \ \ and 
$|y| > |\varphi|_{\{\alpha,\beta\}}$;

$\varphi_a \psi_b(bb) = b \, \psi(b) = bza$, \ \  
$\varphi_a \psi_b(bab) = b \, \psi(ab) = bzb$, \ \  and 
$|z| > |\psi|_{\{\alpha,\beta\}}$. 

\medskip

\noindent {\sc Case} $(* \beta, * \alpha)$ \ \ \ \ \  

$\varphi_a \psi_b(ab) = a \, \varphi(b) = aya$, \ \  
$\varphi_a \psi_b(aab) = a \, \varphi(ab) = ayb$, \ \ and 
$|y| > |\varphi|_{\{\alpha,\beta\}}$;

$\varphi_a \psi_b(ba) = b \, \psi(a) = bza$, \ \  
$\varphi_a \psi_b(bba) = b \, \psi(ba) = bzb$, \ \ and 
$|z| > |\psi|_{\{\alpha,\beta\}}$. 

\medskip

\noindent {\sc Case} $(* \beta, * \beta)$ \ \ \ \ \  

$\varphi_a \psi_b(ab) = a \, \varphi(b) =  aya$, \ \  
$\varphi_a \psi_b(aab) = a \, \varphi(ab) = ayb$, \ \ and 
$|y| > |\varphi|_{\{\alpha,\beta\}}$;

$\varphi_a \psi_b(bb) = a \, \psi(b) = aza$, \ \  
$\varphi_a \psi_b(bab) = b \, \psi(ab) = bzb$, \ \ and 
$|z| > |\psi|_{\{\alpha,\beta\}}$. 

\medskip

\noindent Note that by Lemma 7 (page 300) of \cite{ESc}, there is an 
isomorphism between $\langle \alpha,\beta \rangle$ and 
$\langle \alpha_a, \beta_a \rangle$, mapping generators to generators.
Hence \     
$|\varphi|_{\{\alpha,\beta\}} = |\varphi_a|_{\{\alpha_a,\beta_a\}}$, and 
 \ $|\psi|_{\{\alpha,\beta\}} = |\psi_b|_{\{\alpha_b,\beta_b\}}$.
Also, the generators $\alpha_b,\beta_b$ cannot occur in a shortest word 
over $\{\alpha_a,\beta_a,\alpha_b,\beta_b\}$ representing $\varphi_a$; thus 
 \ $|\varphi|_{\{\alpha,\beta\}} = $
$|\varphi_a|_{\{\alpha_a,\beta_a,\alpha_b,\beta_b\}}$.
Similarly, \ $|\psi|_{\{\alpha,\beta\}} = $
$|\psi_b|_{\{\alpha_a,\beta_a,\alpha_b,\beta_b\}}$.

\smallskip

The existence of $y$ and $z$ as above implies (as in part {\bf (1)}), that 
the table of the {\em maximum extension} of $\varphi_a \psi_b$ contains an 
entry of length 
$> |\varphi_a|_{\{\alpha_a,\beta_a,\alpha_b,\beta_b\}}$, and an entry of 
length
$> |\psi_b|_{\{\alpha_a,\beta_a,\alpha_b,\beta_b\}}$; hence, by Lemma 
\ref{lengthCode}, the table size of the maximum extension of 
$\varphi_a \psi_b$ satisfies  

\medskip
  
$\|{\sf max} \, \varphi_a \psi_b\| \ > \ $
${\sf max}\{|\varphi_a|_{\{\alpha_a,\beta_a,\alpha_b,\beta_b\}}, $
$|\psi_b|_{\{\alpha_a,\beta_a,\alpha_b,\beta_b\}} \}$
$ \ > \ \frac{1}{2} $
$(|\varphi_a|_{\{\alpha_a,\beta_a,\alpha_b,\beta_b\}} $ 
$ + |\psi_b|_{\{\alpha_a,\beta_a,\alpha_b,\beta_b\}})$.

\medskip

\noindent Moreover, by Corollary \ref{lengthVStable}, \
$\|{\sf max} \, \, \varphi_a \psi_b\| \leq \|\varphi_a \psi_b\|$
$\leq C_{\Delta} \, |\varphi_a \psi_b|_{\Delta}$.
Therefore, \ 

\smallskip

$|\varphi_a \psi_b|_{\Delta} \ > \ c \, $
$(|\varphi_a|_{\{\alpha_a,\beta_a,\alpha_b,\beta_b\}}$
$  + |\psi_b|_{\{\alpha_a,\beta_a,\alpha_b,\beta_b\}})$, 

\smallskip

\noindent for some constant $c > 0$. Obviously, \ 
$|\varphi_a|_{\{\alpha_a,\beta_a,\alpha_b,\beta_b\}} $ 
$  + |\psi_b|_{\{\alpha_a,\beta_a,\alpha_b,\beta_b\}} \ \geq \ $
$|\varphi_a \psi_b|_{\{\alpha_a,\beta_a,\alpha_b,\beta_b\}}$.
Thus, 

\smallskip

$|\varphi_a \psi_b|_{\Delta} \ > \ c \, $
$|\varphi_a \psi_b|_{\{\alpha_a,\beta_a,\alpha_b,\beta_b\}}$,

\smallskip

\noindent  which completes the proof of part {\bf (2)} of
Theorem \ref{linearDist}.
 \ \ \ $\Box$

\bigskip

The following is of independent interest.

\begin{pro} 
\label{distProd} \  
 If $G$ is a finitely generated subgroup of $V$ with a distortion function 
$\delta$ then $V$ also contains a subgroup isomorphic 
to the direct product $G \times G$, with distortion function linearly 
equivalent to $n \mapsto \delta(n \log n)$.
\end{pro}
{\bf Proof.} \   Recall that for each element 
$\varphi \in G \subseteq V$ we considered the two elements 
$\varphi_a, \varphi_b \in V$ defined as follows:

\medskip

 \ \ \ $\varphi_a(ax) =  a \, \varphi(x)$,  \ \ for all $x$
          in the domain of $\varphi$ ;

\smallskip

 \ \ \ $\varphi_a(b) = b$.

\medskip

 \ \ \ $\varphi_b(a) = a$ ;

\smallskip

 \ \ \ $\varphi_b(bx) = b \, \varphi(x)$, \ \ for all $x$
            in the domain of $\varphi$.

\medskip

Let $\Delta$ be a finite generating set of $V$, and let $\Gamma$ 
be a finite generating set of $G$. 
For $\varphi \in G$ let $n = |\varphi|_{\Delta}$. Then 
$|\varphi|_{\Gamma} \leq \delta(n)$. Let \  
$\Delta_a = \{\theta_a : \theta \in \Delta\}$, and 
$\Delta_b = \{\theta_b : \theta \in \Delta\}$. Then 
$\Delta_a \cup \Delta_b$ is a finite generating set of a subgroup isomorphic 
to $V \times V$.
Similarly,  $G \times G$ is generated by $\Gamma_a \cup \Gamma_b$. 
Moreover, for any $(\varphi_a, \psi_b) \in G \times G$ we have:
$|(\varphi_a, \psi_b)|_{\Gamma_a \cup \Gamma_b} = $
$|\varphi_a \cdot \psi_b|_{\Gamma_a \cup \Gamma_b} = $
$|\varphi_a|_{\Gamma_a} + |\psi_b|_{\Gamma_b}$
$\leq \delta(|\varphi_a|_{\Delta_a}) + \delta(|\psi_b|_{\Delta_b})$
$\leq 2 \cdot \delta(|\varphi_a \psi_b|_{\Delta_a \cup \Delta_b})$. 
The latter inequality uses the fact that distortion functions are 
non-decreasing, by definition.
Moreover we have

\smallskip

$|\varphi_a \psi_b|_{\Delta_a \cup \Delta_b} \ \leq \ $
$|\varphi_a|_{\Delta_a \cup \Delta_b} + |\psi_b|_{\Delta_a \cup \Delta_b}$
$ \ = \ |\varphi|_{\Delta} + |\psi|_{\Delta}$.

\smallskip

\noindent The last equality holds by Lemma 7 (page 300) of \cite{ESc}, 
as we observed in the proof of part (2) of Theorem \ref{linearDist}.  
Next, by Theorem \ref{sizes}, 

\smallskip

$|\varphi|_{\Delta} + |\psi|_{\Delta} \ \leq \ C_{\Delta} \, $
$(\|\varphi\| \, \log_2 \|\varphi\| + \|\psi\| \, \log_2 \|\psi\|)$ 

\smallskip

$ \ \ \ \ \ \leq  \ C_{\Delta} \, (\|\varphi\| + \|\psi\|) \, $
$\log_2 (\|\varphi\| + \|\psi\|)$.

\medskip

\noindent {\sc Claim:} \ For all $\varphi, \psi \in V$, \ \    
$\|\varphi_a \psi_b\| = \|\varphi\| + \|\psi\|$.

\smallskip

\noindent Proof of the Claim: Suppose the tables for $\varphi, \psi$, in 
maximally extended form are 

$ \varphi = \left[ \begin{array}{ccc}
        x_1 & \ldots & x_m \\
        y_1 & \ldots & y_m
\end{array}        \right], \ \   
\psi = \left[ \begin{array}{ccc}
        u_1 & \ldots & u_n \\
        v_1 & \ldots & v_n
\end{array}        \right]. $  \\  
Then \ 
$ \varphi_a = \left[ \begin{array}{ccc c}
        ax_1 & \ldots & ax_m & b \\
        ay_1 & \ldots & ay_m & b 
\end{array}        \right], \ \  
\psi_b = \left[ \begin{array}{c ccc}
      a & bu_1 & \ldots & bu_n \\
      a & bv_1 & \ldots & bv_n
\end{array}        \right]$, \ and  

\smallskip
  
$ \varphi_a \psi_b = \left[ \begin{array}{ccc ccc}
        ax_1 & \ldots & ax_m & \ bu_1 & \ldots & bu_n \\
        ay_1 & \ldots & ay_m & \ bv_1 & \ldots & bv_n
\end{array}        \right]. $   

\smallskip

\noindent  
It is easy to see that the latter table is in maximally extended from.
Indeed, no pair $ax_i, bu_j$ can lead to extension; and since the original
tables of $\varphi$ and $\psi$ were maximally extended already, no extension
can happen among pairs $ax_i, ax_j$ or pairs $bu_i, bu_j$. 
This proves the Claim.

\medskip

\noindent By the Claim, and the inequalities proved just before the Claim, 
we obtain:

\smallskip

$|\varphi_a \psi_b|_{\Delta_a \cup \Delta_b} \ \leq \ C_{\Delta} \, $
$\|\varphi_a \psi_b\| \, \log_2 \|\varphi_a \psi_b\|$

\smallskip

\noindent By the first inequality of Theorem \ref{sizes} this implies \  
$|\varphi_a \psi_b|_{\Delta_a \cup \Delta_b} \ \leq \ c \, $
$|\varphi_a \psi_b|_{\Delta} \, \log_2 |\varphi_a \psi_b|_{\Delta}$

\noindent for some constant $c > 0$. We already proved \ 
$|(\varphi_a, \psi_b)|_{\Gamma_a \cup \Gamma_b} \ \leq \ $
$2 \cdot \delta(|\varphi_a \psi_b|_{\Delta_a \cup \Delta_b})$.
Therefore we have \  

\smallskip

$|(\varphi_a, \psi_b)|_{\Gamma_a \cup \Gamma_b} \ \leq \ 2 \  $
$\delta(C\, |\varphi_a\psi_b|_{\Delta}\, \log_2 |\varphi_a\psi_b|_{\Delta})$
 
\smallskip

\noindent for some constant $C > 0$. This proves the Proposition.   
  \ \ \ $\Box$

\bigskip

We will need the following definitions:  
\begin{defn} 
\label{polynEquivDef} \  
Two functions $f_1$ and $f_2: \mathbb{N} \to \mathbb{N}$ are
{\bf polynomially equivalent} iff there exist positive constants $c_0$, 
$C_{1 2}$, $c_{1 2}$, $C_{2 1}$, $c_{2 1}$, $d_1$, $d_2$, $d_3$, $d_4$ such 
that for all $n \geq c_0$:   

\smallskip

  \ \ \ \ $f_1(n) \leq C_{1 2} \cdot (f_2(c_{1 2} \, n^{d_1})^{d_3})$ \ and \
             $f_2(n) \leq C_{2 1} \cdot (f_1(c_{2 1} \, n^{d_2})^{d_4})$.
\smallskip

\noindent In big-O notation this means \ 
$f_1(n) = O(f_2(O(n^{O(1)}))^{O(1)})$ \  
and \ $f_2(n) = O(f_1(O(n^{O(1)}))^{O(1)})$.
\end{defn}

\begin{defn}
A function $f : \mathbb{N} \to \mathbb{N}$ is {\bf superadditive} iff
for all $n, m$: \ \ $f(n + m) \geq f(n) + f(m)$.
\end{defn}

\begin{pro}
\label{superAdd} \
Let $t: \mathbb{N} \to \mathbb{N}$ be a non-decreasing function such that
$t(n) \geq n$ for all $n$. Then there is a superadditive function
$T: \mathbb{N} \to \mathbb{N}$ such that for all $n$,

\smallskip

 \ \ \ \ \ \ \ $t(n) \leq T(n) \leq n \cdot t(n)$.

\smallskip

\noindent Hence, every non-decreasing function which is larger than the
identity function is polynomially equivalent to a non-decreasing 
{\em superadditive} function.
\end{pro}
{\bf Proof.} \ Given the function $t$ we define the desired superadditive
function $T$ by

\smallskip

 $T(n) = {\sf max} \{ \ \sum_{i=1}^k t(n_i) \ : \ k \geq 1, \ $
 $(n_1, \ldots, n_k) \in (\mathbb{N}-\{0\})^k, \ \ \sum_{i=1}^k n_i = n \}$.

\smallskip

\noindent In other words, $(n_1, \ldots, n_k)$ is any partition of $n$; we
define $T(n)$ to be $\sum_{i=1}^k t(n_i)$, maximized over all partitions 
of $n$.
It is straightforward to verify that $T$ is superadditive and that \  
$t(n) \leq T(n) \leq n \cdot t(n)$.  \ \ \ $\Box$

\bigskip

Theorem \ref{distTfin} gives us a large subset of the set of distortion 
functions of $V$, namely the set of all Dehn functions. Moreover, 
we know from \cite{SBR} that for any time complexity function of a 
nondeterministic Turing machine, its fourth power is linearly equivalent 
to a Dehn function (if it is superadditive). By Proposition \ref{superAdd}, 
we can always assume 
that our functions are superadditive (up to polynomial equivalence).

Corollary \ref{threeEq} below will give a kind of converse to 
Theorem \ref{distTfin}. 

\begin{lem}  \  
If a finitely generated subgroup $H$ of a finitely generated group $G$ has 
distortion function $\delta$, and if the word problem of $G$ has 
nondeterministic time complexity $T(\cdot)$, then the generalized word problem 
of $H$ in $G$ has nondeterministic time complexity bounded by a function 
linearly equivalent to $T(\delta(\cdot))$.

Moreover, if $T$ is the time complexity of a nondeterministic Turing 
machine and $\delta$ is the distortion of a subgroup of $G$ then 
$T(\delta(\cdot))$ is linearly equivalent to the time complexity of a 
nondeterministic Turing machine.  
\end{lem}
{\bf Proof.} \ Suppose $G$ has a finite generating set $A$ and $H$ has a finite
generating set $\Delta \subset (A^{\pm 1})^*$.   
Let $w \in (A^{\pm 1})^*$ of length $|w| = n$ be an input to the generalized
word problem. We guess a word $z \in (\Delta^{\pm 1})^*$ of length 
$|z| \leq \delta(n)$ such that $w = z$ in $G$. Such a $z$ exists iff the 
answer to the generalized word problem is ``yes''. It takes time $O(\delta(n))$
to guess $z$.

To check correctness of the guessed $z$ we solve the word problem  
``Is $w = z$ in G?'', in nondeterministic time linearly equivalent to 
$T(n + \delta(n))$, which is linearly equivalent to $T(\delta(\cdot))$.
 \ \ \ $\Box$  

\medskip

If we let $G$ be $V$, the above and Proposition \ref{wrdprobTime} 
give us a nondeterministic Turing machine with time complexity linearly
equivalent to $\delta(\cdot)^2$, for solving the generalized word problem 
of $H$ in $V$. Therefore we have:

\begin{pro} \  
Every distortion of $V$ is linearly equivalent to the square-root 
of a nondeterministic time complexity function.
\end{pro}

\noindent In summary, we have the following corollary.

\begin{cor}
\label{threeEq}
 \ The following classes of functions are polynomially equivalent:  

\smallskip

\noindent
$\bullet$ \ Time complexity functions of nondeterministic Turing machines. \\ 
$\bullet$ \ Dehn functions of finitely presented groups. \\  
$\bullet$ \ Distortions in the Thompson group $V$. 
\end{cor} 
The polynomial equivalences between the above three classes of functions 
actually have ``uniform degree''; this means that the degrees $d_1$, $d_2$,
$d_3$, $d_4$ that appear in the polynomial equivalences between functions 
(in Definition \ref{polynEquivDef}) are the same for all pairs of functions.

\bigskip

\noindent {\bf Questions}: \ {\bf (1)} 
Of course (by Ol'shanskii and Sapir's theorem \cite{OlSap}, and by the results
of \cite{SBR}), Corollary \ref{threeEq} also holds for 
${\rm FG}_2 \times {\rm FG}_2$. For what other groups does the Corollary hold? 

Some linear groups are candidates; it is well known that SL$_4(\mathbb{Z})$ has 
${\rm FG}_2 \times {\rm FG}_2$ as a subgroup (see pp.~41-42 in \cite{Miller}, 
and \cite{Miller92}), but one would also need to study the distortion of  
linear groups as subgroups of other linear groups.   

\smallskip

\noindent
{\bf (2)} Are the distortions of $V$ linearly (rather than just 
polynomially) equivalent to the Dehn functions of finitely presented groups? 
(By \cite{OlSap}, ${\rm FG}_2 \times {\rm FG}_2$ has this stronger property.)


\section{Representation of the Thompson groups in algebras}

We will show that the Thompson groups are subgroups of Cuntz 
C$^\star$-algebras. Those algebras can be defined as the completion of 
quotient algebras of the polycyclic monoid. 

The {\it polycyclic monoid} on a generating set $A$ is an inverse monoid 
with zero $\bf{0}$, defined by an inverse monoid presentation with 
relations   

\smallskip

$\{\alpha^{-1}\alpha = {\bf 1} : \alpha \in A \} \ \cup $
$\{\alpha^{-1}\beta =  {\bf 0} : $
$\alpha, \beta \in A$ with $\alpha \neq \beta \}$.

\smallskip

\noindent
It follows from this presentation that every element of $PC(A)$,
other than {\bf 1} and {\bf 0},
is of the form $yx^{-1}$ with $y, x \in A^*$. On the other hand,
$u^{-1}v = {\bf 0}$ if $u, v \in A^*$ are not prefix-comparable; and
if $u \geq_{\rm pref} v = uz$, then $u^{-1}v = z \in A^*$; \ if
$v \geq_{\rm pref} u = vw$, then $u^{-1}v = w^{-1} \in (A^{-1})^*$.

For the definition of ``inverse monoid'' and more information on
these monoids, see e.g.\ the monograph \cite{LawsonInv}.
Polycyclic monoids were introduced in \cite{NivatPerrot} (and were
re-invented in \cite{Cuntz}).
The papers \cite{MeakinSapir}, \cite{HinesLawson} and \cite{Lawson2}
give interesting applications of polycyclic monoids.

\bigskip

We assume that the alphabet is $A = \{a, b\}$, and in order to represent
$V$ we first consider the monoid algebra of $PC(a,b)$ over any field
$\mathbb{K}$

\medskip

$\mathbb{K}[PC(a,b)] \ = \  
\{ \sum_{i=1}^n \kappa_i \,y_i x_i^{-1} \ : \ n \in \mathbb{N}, \
y_i x_i^{-1} \in PC(a,b) \ {\rm and} \ \kappa_i \in \mathbb{K} \
{\rm for \ all} \ i = 1, \ldots, n \}$

\smallskip

\hspace{1.1in}  $\cup \ \ \{ {\bf 0} \}$

\bigskip

\noindent Before embedding $V$ in a Cuntz algebra, we will represent $V$ 
as a subgroup of the multiplicative part of a quotient algebra of 
$\mathbb{K}[PC(a,b)]$.
Let us look at an example before going into details.

\bigskip

\noindent {\bf Example:} Consider the two elements $C, B$ of $V$
given by tables
\[ B \ = \ \left[ \begin{array}{cccc}
        a \ & \ ba   \ & \ b^2a \ & \ b^3 \\
        a \ & \ ba^2 \ & \ bab \  & \ b^2
\end{array}        \right], \ \ \ \ \
C \ = \ \left[ \begin{array}{ccccc}
        a^2 \ & \ ab \ & \ bab  \ & \ ba^2 \ & \ b^2 \\
        a   \ & \ ba \ & \ b^3a \ & \ b^2a \ & \ b^4
\end{array}        \right]. \]
They will be represented by elements of $\mathbb{K}[PC(a,b)]$ as follows:

\smallskip

$B$ is represented by \ \
$aa^{-1} + ba^2 a^{-1}b^{-1} + baba^{-1}b^{-2} + b^2b^{-3}$,

\smallskip

$C$ is represented by \ \
$aa^{-2} + bab^{-1}a^{-1} + b^3ab^{-1}a^{-1}b^{-1} + b^2aa^{-2}b^{-1}$
$ + b^4b^{-2}$.

\smallskip

\noindent
One observes that the composite $C \circ B$ is then represented by the
product of the corresponding elements of $\mathbb{K}[PC(a,b)]$:

\smallskip

$(aa^{-2} + bab^{-1}a^{-1} + b^3ab^{-1}a^{-1}b^{-1} + b^2aa^{-2}b^{-1}$
$ + b^4b^{-2})$ $\cdot$
$(aa^{-1} + ba^2 a^{-1}b^{-1} + baba^{-1}b^{-2} + b^2b^{-3}) \ $

$ = \ aa^{-2} + bab^{-1}a^{-1} + b^3aa^{-1}b^{-2} + b^2aa^{-1}b^{-1} $
  $+ b^4b^{-3}$ ,

\smallskip

\noindent where we applied the relations of $PC(a,b)$ and omitted the terms
that are {\bf 0}.

The maximum extension of $C \circ B$ turns out to be \
$A = \left[ \begin{array}{ccc}
         a^2 \ & \ ab \ & \ b \\
         a   \ & \ ba \ & \ b^2
\end{array}        \right] $.
The $\mathbb{K}[PC(a,b)]$-representation of $A$ can be obtained from the
$\mathbb{K}[PC(a,b)]$-representation of $C \circ B$ by repeatedly applying
the relation \ $aa^{-1} + bb^{-1} = {\bf 1}$, as follows:

\medskip

$aa^{-2} + bab^{-1}a^{-1} + b^3aa^{-1}b^{-2} + b^2aa^{-1}b^{-1} + b^4b^{-3}$

\smallskip

$ \ = \ aa^{-2}+ bab^{-1}a^{-1}+ b^3(aa^{-1}+bb^{-1})b^{-2}+ b^2aa^{-1}b^{-1}$

\smallskip

$ \ = \ aa^{-2} + bab^{-1}a^{-1} + b^3b^{-2} + b^2aa^{-1}b^{-1}$

\smallskip

$ \ = \ aa^{-2} + bab^{-1}a^{-1} + b^2(bb^{-1}+aa^{-1})b^{-1}$

\smallskip

$ \ = \ aa^{-2} + bab^{-1}a^{-1} + b^2b^{-1}$.

\medskip

\noindent This ends the Example. We will now formalize this for
$V$, then generalize the representation to ${\mathcal G}_{2,1}$, and
finally prove properties.

\medskip

For an algebra $\cal A$ and a set $X \subseteq {\cal A}$ we write
$\langle \langle X \rangle \rangle$ for the {\it ideal} generated by $X$ in
$\cal A$.  We now take the following quotient algebra:
\begin{center}
${\cal C}_V \ =\ \mathbb{K}[PC(a,b)]/ {\bf I}_V$
\end{center}
where \
\begin{center}
${\bf I}_V \ = \ $
$\langle \langle a \, a^{-1} + b \, b^{-1} - \bf{1} \rangle \rangle$.
\end{center}
We can interpret the ideal {\bf I}$_V$ as a {\it term rewrite system}
with the two rules, \
$a \, a^{-1} + b \, b^{-1} \to {\bf 1}$ \ and \
${\bf 1} \to a \, a^{-1} + b \, b^{-1}$.
The ideal ${\bf I}_V$ and the corresponding rewrite system
interpretation are inspired from Lemma \ref{ExtCritFin}, which tells us
how to find maximum extensions of right-ideal isomorphisms.

\smallskip

Note that in the quotient algebra ${\cal C}_V$ we have \
$\sum_{q \in Q} q \, q^{-1} = {\bf 1}$ \ for any finite maximal prefix code
$Q$. Hence, in the above example, the sum representing $C \circ B$ could have
been maximally extended in one step as follows:

\smallskip

$aa^{-2} + bab^{-1}a^{-1} + b^3aa^{-1}b^{-2} + b^2aa^{-1}b^{-1} + b^4b^{-3}$

$ \ = \ aa^{-2}+ bab^{-1}a^{-1}+ b^2(baa^{-1}b^{-1}+aa^{-1}+b^2b^{-2})b^{-1}$

$ \ = \ aa^{-2}+ bab^{-1}a^{-1}+ b^2b^{-1}$

\smallskip

\noindent where $\{ ba, a, b^2\}$ is a maximal prefix code.

\bigskip

Representing ${\mathcal G}_{2,1}$ by algebras is more complicated.
For a right-ideal
isomorphism $\varphi$ between essential right ideals we want the
representation \ $\sum_{x \in P} \varphi(x) \, x^{-1}$, where
$P \subset A^*$ is the (possibly infinite) domain code of $\varphi$.
Then for any $w \in PA^*$ we have:

\medskip

 \ \ \ \ \ \ \ \ \ \ \ $\varphi(w) \ = $
 $ \ \sum_{x \in P} \varphi(x) \, x^{-1} \ w \ \ \cap \ \ A^*$.

\medskip

\noindent
The sum \ $\sum_{x \in P} \varphi(x) \, x^{-1}$ \ belongs to an algebra
consisting of possibly infinite sums over the monoid $PC(a,b)$ over any field
$\mathbb{K}$; but we need to restrict the infinite sums in order to get the
desired algebraic properties. The following property of sums will
guarantee that our algebra is closed under multiplication, but further
restrictions will be needed.

\begin{defn} \
We call a relation $S \subset A^* \times A^*$ {\em finite-to-finite}
iff for any $x \in A^*$ there are only finitely many $y \in A^*$
such that $(x,y) \in S$, and for any $y \in A^*$ there are only finitely
many $x \in A^*$ such that $(x,y) \in S$.
\end{defn}
As a preliminary step we start out with the following algebra.

\bigskip

${\cal B}_{\infty} \ = \
\{ \sum_{i \in I} \kappa_i \,y_i x_i^{-1} \ : \ I \subseteq \mathbb{N},
 \ \ y_i x_i^{-1} \in PC(a,b) \ {\rm and} \ \kappa_i \in \mathbb{K}
 \ {\rm for \ all} \ i \in I,   $

\smallskip

 \makebox[1.8in]{}  and the relation
    $\{ (y_i, x_i) : i \in I\}$ is finite-to-finite$\}$
 \ $\cup \ \ \{ \bf{0} \}$

\bigskip

\noindent It is obvious that ${\cal B}_{\infty}$ is closed under addition.
In Lemma \ref{BclosedMult} we'll prove that ${\cal B}_{\infty}$ is closed under
multiplication.

\bigskip

We define the set of {\bf unary sums} to be the subset ${\cal U}_{\infty}$
of ${\cal B}_{\infty}$ consisting of sums
of the form $\sum_{i\in I} y_i x_i^{-1}$ with the following properties:

\smallskip

\noindent $\bullet$ \ All the coefficients in the sum are equal to {\bf 1}.

\smallskip

\noindent $\bullet$ \ The sets $\{x_i: i\in I\}$ and $\{y_i: i\in I\}$ are
{\it maximal} prefix codes, such that the indexing $i \in I \mapsto x_i$
and the indexing $i \in I \mapsto y_i$ are bijective functions.

\smallskip

The set ${\cal U}_V \ (\subset \mathbb{K}[PC(a,b)])$ is defined in a
similar way, by taking the index sets $I$ to be finite in the above 
definition.

\smallskip

We also define the sets of {\bf partial unary sums},
${\cal U}_{\infty}^{\rm part}$ and ${\cal U}_V^{\rm part}$, by
just requiring $\{x_i: i\in I\}$ and $\{y_i: i\in I\}$ to be prefix codes
(not necessarily maximal),
while the rest of the definitions is kept unchanged.

\medskip

We will prove later (Lemma \ref{Umult}) that the sets ${\cal U}_{\infty}$,
${\cal U}_V$, ${\cal U}_{\infty}^{\rm part}$ and
${\cal U}_V^{\rm part}$ are closed under multiplication, so they are monoids.

\medskip

We define the algebra ${\cal A}_{\infty}$ as the subalgebra of
${\cal B}_{\infty}$ generated as an algebra by ${\cal U}_{\infty}^{\rm part}$:
$${\cal A}_{\infty} \ = \ \langle {\cal U}_{\infty}^{\rm part} \rangle.$$
Note that $\mathbb{K}[PC(a,b)]$ is the subalgebra
of ${\cal B}_{\infty}$ generated by ${\cal U}_V^{\rm part}$.
Since ${\cal U}_{\infty}^{\rm part}$ is closed under multiplication (as we
will prove in Lemma \ref{Umult}), every element of ${\cal A}_{\infty}$ is a
{\it linear combination} of elements of ${\cal U}_{\infty}^{\rm part}$;
in other words, ${\cal A}_{\infty}$ is the {\it monoid algebra} of the
monoid ${\cal U}_{\infty}^{\rm part}$.

The algebra ${\cal A}_{\infty}$ can also be characterized as follows
(as will be proved in Lemma \ref{charactA}):
${\cal A}_{\infty}$ consists of the elements
$\sum_{i \in I} \kappa_i \,y_i x_i^{-1}$ of ${\cal B}_{\infty}$ that have
the following three properties.

\smallskip

\noindent (1) The relation $S = \{ (x_i, y_i) : i \in I\}$ is
{\it bounded finite-to-finite}. This means that there exists $n_0$ such that
for every $x_{i_0} \in \{x_i : i \in I\}$, \
 $|S(x_{i_0})| = |\{ y_j : (x_{i_0}, y_j) \in S\}| \leq n_0$, and
for every $y_{j_0} \in \{y_i : i \in I\}$, \
$|S^{-1}(y_{j_0})| = |\{ x_i : (x_i, y_{j_0}) \in S\}| \leq n_0$.
(i.e., there is a bound on the cardinalities
of all the sets $S(x_i)$ and $S^{-1}(y_i)$ as $i$ ranges over $I$).

\smallskip

\noindent (2) In $\{x_i : i \in I\}$ and in $\{y_i : i \in I\}$, all
$>_{\rm pref}$-chains have bounded length.

\smallskip

\noindent (3) The set $\{ \kappa_i : i \in I \}$ is finite.

\bigskip

\noindent
Finally, in order to embed ${\mathcal G}_{2,1}$ we consider the quotient
algebra

\begin{center}
${\cal C}_{\infty} \ = \ {\cal A}_{\infty} / {\bf I}_{\infty}$
\end{center}

\noindent where the ideal ${\bf I}_{\infty}$ is defined by

\bigskip

${\bf I}_{\infty} \ = \ \langle \langle \ \sum_{i \in I}$
$y_i \, (\sum_{q \in Q_i} q \, q^{-1} - {\bf 1}) \, x_i^{-1}$
$ \ : \ \ I \subseteq \mathbb{N}$, \

\smallskip

 \makebox[1in]{} $\{x_i : i \in I\}$ and
 $\{y_i : i \in I\}$ are maximal prefix codes,

\smallskip

 \makebox[1in]{} $i \in I \mapsto x_i \in \{x_i : i \in I\}$ \ and
    \ $i \in I \mapsto y_i \in \{y_i : i \in I\}$ are bijections,

\smallskip

 \makebox[1in]{} and $Q_i$ is a maximal prefix code of $A^*$
                 (for every $i \in I$)
$\rangle \rangle$.

\bigskip

Note that $\bigcup_{i \in I} x_iQ_i$ and $\bigcup_{i \in I} y_iQ_i$ are maximal
prefix codes, by construction (2) in Example \ref{exPrefCod}; hence
$\sum_{i \in I} \sum_{q \in Q_i} y_i \, q \, q^{-1} x_i^{-1} \in $
$U_{\infty}^{\rm part}$, and thus
${\bf I}_{\infty} \subset {\cal A}_{\infty}$.
We can interpret the ideal ${\bf I}_{\infty}$ as a {\it generalized
term rewrite system} with the rules \ $\sum_{q \in Q} q \, q^{-1} \to {\bf 1}$
 \ and \ ${\bf 1} \to \sum_{q \in Q} q \, q^{-1}$, \ where $Q$ ranges over all
maximal prefix codes of $A^*$. This is a {\it generalized} rewrite system, in
the sense that these rules are applied to infinite sums, and we allow
infinitely many rules to be applied ``in parallel''; i.e., infinitely many
rules
of the form \ $\sum_{q \in Q_i} q \, q^{-1} \rightleftarrows {\bf 1}$ \
$(i \in I)$ \  are applied in infinitely many non-overlapping
locations in the infinite sum.

The ideals ${\bf I}_V$ and ${\bf I}_{\infty}$ (and the corresponding
rewrite system interpretation) are inspired by Lemmas \ref{ExtCritFin},
\ref{ExtCritInf} and \ref{findMaxExt} (especially Remark 3), which tell us
how to find maximum extensions of right-ideal isomorphisms.

The congruence class of \ $\sum_{i \in I} \kappa_i \,y_i x_i^{-1}$ \ is \
$\sum_{i \in I} \kappa_i \,y_i x_i^{-1} + {\bf I}_V$, \ respectively
 \ $\sum_{i \in I} \kappa_i \,y_i x_i^{-1} + {\bf I}_{\infty}$.
Two elements of $\mathbb{K}[PC(a,b)]$ or ${\cal A}_{\infty}$ that belong to 
the same congruence class are called {\it congruent}.

\medskip

We will obtain the following representations of $V$ and ${\mathcal G}_{2,1}$
as subgroups of the multiplicative part of ${\cal C}_V$,
respectively ${\cal C}_{\infty}$:

\begin{thm}
\label{thmThompson} \   
The Thompson group $V$ is isomorphic to the subgroup 
${\cal U}_V / {\bf I}_V$ of the multiplicative part of the algebra 
${\cal C}_V \ = \ \mathbb{K}[PC(a,b)] / {\bf I}_V$.

The Thompson group ${\mathcal G}_{2,1}$ is isomorphic to the subgroup 
${\cal U}_{\infty} / {\bf I}_{\infty}$ of the multiplicative part of the 
algebra ${\cal C}_{\infty} \ = \ {\cal A}_{\infty} / {\bf I}_{\infty}$.
\end{thm}
The fact that only the elements {\bf 0} and {\bf 1} of the field are used in
this representation is reminiscent of the regular representation by
permutation matrices.

\medskip

All the algebras above are {\it $\star$-algebras}, by defining \
$(\sum_{i \in I} \kappa_i \, y_ix_i^{-1})^{\star} $
$= \sum_{i \in I} \kappa_i^{\star} \, x_iy_i^{-1}$ \ and taking
the field $\mathbb{K}$ to be $\mathbb{C}$ (the complex numbers, where $\star$
denotes conjugation). Then for all algebra elements $s, t$
and scalars $\kappa \in \mathbb{C}$ :
 \ $(s^{\star})^{\star} = s, \ \ \ $
$(\kappa s + t)^{\star} = \kappa^{\star} s^{\star} + t^{\star}, \ \ \ $
$(st)^{\star} = t^{\star}s^{\star}$.

\begin{remark} -- Connection with the Cuntz algebras \footnote{I owe the 
observation of the connection between the above algebra ${\cal C}_V$ and 
the Cuntz algebras to John Meakin \cite{Meakin}.}

\end{remark}
\vspace{-.15in}
The use of the polycyclic monoid and the relations \
$\sum_{q \in Q} q \, q^{-1} = {\bf 1}$ \
means that the Cauchy completion of our algebra ${\cal C}_V$ (defined
for any alphabet $A$) is the {\bf Cuntz} C$^{\star}$-algebra
${\cal O}_{|A|}$. These C$^{\star}$-algebras were first introduced by Dixmier
\cite{Dixmier}, then studied by Cuntz \cite{Cuntz} who proved many remarkable
properties, many of which are reminiscent of the properties of $V$
itself (e.g., that ${\cal O}_{|A|}$ is a simple algebra). See also
\cite{KDavidson},\cite{APaterson} and \cite{Renault}, where connections
between C$^{\star}$-algebras and inverse monoids and, in particular, the 
relation between the Cuntz algebras and the polycyclic monoid, are exposited.

In summary, from this and Theorem \ref{thmThompson} we obtain: 

\begin{cor}
The Thompson group $V$ is a subgroup of the multiplicative part of the 
Cuntz-Dixmier {\rm C}$^{\star}$-algebra ${\cal O}_2$. 
More generally, the Thompson-Higman group $G_{n,1}$  
$(n \in \mathbb{N}, n \geq 2)$ is a 
subgroup of the multiplicative part of the Cuntz {\rm C}$^{\star}$-algebra 
${\cal O}_n$. 
\end{cor}

\medskip

\noindent We outline the proof of Theorem \ref{thmThompson}. First we list
some lemmas. Other lemmas that play an indirect role, as well as the proofs 
of all the lemmas, are given in Appendix A4. 

\medskip

\noindent $\bullet$ {\it Lemma} \ref{corresp}. \
There is a one-to-one correspondence between (1) the set of all isomorphisms
between (essential) right ideals of $\{a,b\}^*$, and (2) the set
${\cal U}_{\infty}^{\rm part}$ (respectively ${\cal U}_{\infty}$)
$$\Sigma: \ \ \varphi \ \ \longmapsto \ \ 
\sum_{x \in {\rm DomC(\varphi)}} \varphi(x) \, x^{-1} \ \in
{\cal U}_{\infty}^{\rm part}.$$
Its inverse is
$$\Phi: \ \ \sum_{i\in I}y_ix_i^{-1}\in {\cal U}_{\infty}^{\rm part}
  \ \ \longmapsto \ \  
(\varphi : \{x_i:i\in I\} A^* \to \{y_i:i \in I\} A^*)$$
\noindent where $\varphi$ is defined by \
$\varphi(x) = \sum_{i\in I} y_i x_i^{-1} \, x  \ \cap \ A^*$.

Similarly, there is a one-to-one correspondence between (1) the set of all
isomorphisms between finitely generated (essential) right ideals of
$\{a,b\}^*$, and (2) the set ${\cal U}_V^{\rm part}$
(respectively ${\cal U}_V$).

\medskip

\noindent $\bullet$ {\it Lemma} \ref{Umult}. \
The sets ${\cal U}_{\infty}$, ${\cal U}_V$,
${\cal U}_{\infty}^{\rm part}$ and ${\cal U}_V^{\rm part}$ are
closed under multiplication.
We have the following formula for the multiplication in
${\cal U}_{\infty}^{\rm part}$:
\[ \sum_{j \in J} y_j \, x_j^{-1} \cdot \sum_{i \in I} v_i \, u_i^{-1}
 \ = \ \]
\[ \sum_{j \ \in \ {\rm dom}f \ \cap \ {\rm im}g} \hspace{-.3in}
                    y_j \, u_{f(j)}^{-1}
 \ +  \
\sum_{j \ \in \ {\rm dom}f \ - \ {\rm im}g} \hspace{-.3in}
                     y_j \, z_j^{-1} \, u_{f(j)}^{-1}
 \ + \
\sum_{i \ \in \ {\rm dom}g \ - \ {\rm im}f} \hspace{-.3in}
                       y_{g(i)} \, t_i \, u_i^{-1} \]
where $f: J \to I$ and $g: I \to J$ are partial functions defined by
$x_j \leq_{\rm pref} v_{f(j)}$ and $x_j = v_{f(j)} z_j$; similarly, 
$x_{g(i)} \geq_{\rm pref} v_i$ and $v_i = x_{g(i)} t_i$. 

\medskip

\noindent $\bullet$ {\it Lemma} \ref{homU} \
The one-to-one correspondence $\Phi$ of Lemma \ref{corresp} is a
{\em homomorphism}, i.e., for all $\sigma_2$, $\sigma_1$
$ \in {\cal U}_{\infty}^{\rm part} :$ \ \
$\Phi(\sigma_2 \cdot \sigma_1) = \Phi(\sigma_2) \circ \Phi(\sigma_1)$.

\medskip

\noindent $\bullet$ {\it Lemma} \ref{respectEquiv} \
The one-to-one correspondence $\Sigma$ of Lemma \ref{corresp}
respects the congruence relations on the set ${\cal U}_{\infty}$ (induced
by {\bf I}$_{\infty}$) and on the set of all isomorphisms between essential
right ideals. In other words, two isomorphisms between essential right 
ideals, $\varphi_1$ and $\varphi_2$ are congruent (i.e., they have the same 
maximum extension) \ iff \ $\Sigma(\varphi_1)$ and $\Sigma(\varphi_2)$ are 
congruent (relative to the ideal ${\bf I}_{\infty}$).
A similar fact holds for ${\cal U}_V$.

\bigskip

\noindent
{\bf Proof of Theorem \ref{thmThompson}.} We prove the theorem for
${\mathcal G}_{2,1}$; for $V$ the proof is similar.
By Lemma \ref{corresp} there is a one-to-one correspondence $\Sigma$ from
the set of all isomorphisms between essential right ideals, to the set
${\cal U}_{\infty}$.
By Lemma \ref{respectEquiv}, this one-to-one correspondence respects the
congruence relation of the set ${\cal U}_{\infty}$ and the congruence relation
of the set of all isomorphisms between essential right ideals. Therefore 
$\Sigma$ determines a one-to-one correspondence between ${\mathcal G}_{2,1}$ 
and ${\cal U}_{\infty} / {\bf I}_{\infty}$.

Moreover, ${\mathcal G}_{2,1}$ and ${\cal U}_{\infty} / {\bf I}_{\infty}$
are isomorphic as groups:
Any subidentity isomorphism of the form $pw \in PA^* \to pw \in PA^*$ (with
$p$ ranging over a maximal prefix code $P$, $w$ ranging over $A^*$)
is mapped by $\Sigma$ to $\sum_{p\in P} p \, p^{-1}$, which is congruent to
{\bf 1} modulo ${\bf I}_{\infty}$.
Also, for any isomorphism $\varphi: P_1A^* \to P_2A^*$ between essential right
ideals, the one-to-one correspondence $\Sigma$ maps the inverse $\varphi^{-1}$
to \ $\sum_{p_2\in P_2} \varphi^{-1}(p_2) \, p_2$
$= \sum_{p_1\in P_1} p_1 \, \varphi(p_1)^{-1}$, which is congruent to the
inverse of $\Sigma(\varphi)$. Indeed,
$\sum_{i\in I} y_i x_i^{-1} \cdot \sum_{i\in I} x_i y_i^{-1} = $
$\sum_{i\in I} y_i  y_i^{-1}$, and
$\sum_{i\in I} y_i  y_i^{-1} + {\bf I}_{\infty} = {\bf 1} + {\bf I}_{\infty}$.
Finally, the product $\Sigma(\psi) \, \Sigma(\varphi)$ is
congruent to $\Sigma(\psi \varphi)$, as we saw in Lemma \ref{homU}. 
 \ \ \ $\Box$


\newpage 


\section{Appendix A1} 

In this appendix we present basic information about prefix codes and right 
ideals of free monoids.

\medskip

\noindent {\bf A combinatorial observation:} \  
The number of maximal prefix codes of cardinality $n$ over the alphabet 
$\{ a, b\}$ is $C_{n-1}$, where \ $C_n = \frac{1}{n+1}$$2n \choose n$ \   
(the classical Catalan number). This is proved by counting binary trees
(see e.g.~\cite{SedgFlaj}, Chapter 5). 
Asymptotically, \  
$C_n = \frac{4^n}{\sqrt{\pi n^3}} \, (1 + \varepsilon(n))$, \  
with \ lim$_{n\to\infty} \ \varepsilon(n)= 0$.

\smallskip

So, instead of defining the elements of the Thompson group $V$ as 
``maximum'' bijections between maximal prefix codes over the alphabet 
$\{ a, b\}$, one could define $V$ as bijections between other 
combinatorial objects, parameterized by a positive integer $n$, such that the 
number of objects of ``size'' $n$ is the Catalan number (see
e.g.~\cite{SedgFlaj}, Section 5.11).  

\medskip

The next lemmas give some elementary properties of prefix codes.
\begin{lem} 
\label{ideals} 
{\bf (1)} \ $R$ is a right ideal of $A^*$ \ iff \  there exists a 
prefix code $P$ over $A$ such that $R = PA^*$.   

\smallskip

\noindent {\bf (1')} \ For a right ideal $R$ the prefix code $P$ such that 
$R = PA^*$ is unique.  

\smallskip

\noindent 
{\bf (2)} \ $R$ is an essential right ideal of $A^*$ \ iff \ $R$ is a right 
ideal such that 

\smallskip

 \ \ \ \ \ \ $(\forall u \in A^*)(\exists x \in A^*) \ ux \in R$  

\smallskip

\noindent 
(i.e., in the terminology of \cite{ESc}, $R$ is ``inescapable'').  

\smallskip

\noindent   
{\bf (3)} \ $R$ is an essential right ideal of $A^*$ \ iff \ $R$ is a 
right ideal whose corresponding prefix code is maximal.  

\smallskip

\noindent 
{\bf (4)} \ $P$ is a maximal prefix code of $A^*$ \ iff \ $P$ is a prefix 
code and 

\smallskip

 \ \ \ \ \ \ $(\forall u \in A^* - PA^*)(\exists x \in A^*) \ ux \in P$. 
\end{lem}
{\bf Proof of (1):} \ 
$[ \Leftarrow ]$ \  Obviously, $PA^*$ is a right ideal. 

\smallskip

\noindent
$[ \Rightarrow ]$ \  We claim that for any right ideal $R \subseteq A^*$ we
have:

\smallskip 

 \ \ \ \ \ \ \ \ \ \ \ \ $P = R - R \, A$ \ \ is a prefix code such that 
 \ $R = PA^*$.

\smallskip

\noindent
Obviously, $PA^* \subseteq RA^* \subseteq R$, since $R$ is a right ideal.
Conversely, let us show that $R \subseteq PA^*$. For any $r \in R$, let $p$ 
be the shortest prefix of $r$ that belongs to $R$. Since $r$ itself is in 
$R$, $p$ exists, and by definition, $p \in R$. Also, $r = px$ for some 
$x \in A^*$, since $p$ is a prefix of $r$. 
Finally, $p \notin RA$ (otherwise, we would have $p = r'a$ for some 
$r' \in R$, $a \in A$, which would imply that $p$ is not the shortest prefix 
that $r$ has in $R$). Thus, $p \in P$. Since $r = px$, we have
$R \subseteq PA^*$. 

To show that $P$ is a prefix code, let $p, p' \in P$ and suppose $p'$ is a 
prefix of $p$: \ $p = p'x$, for some $x \in A^*$. 
If $x$ is not empty then $p \in RA$, contradicting the assumption 
$p \in P \ (= R - RA)$. Thus, the words in $P$ that are prefixes of each other
are equal to each other.  \ \ 

\medskip

\noindent {\bf Proof of (1'):} \ If $P_1A^* = P_2A^*$ for two prefix  
codes $P_1, P_2$, then for every $p_1 \in P_1$ there exists $p_2 \in P_2$
such that $p_1 = p_2 x$ (for some $x \in A^*$). Also, there is $p_1' \in P_1$
such that $p_2 = p_1' y$ (for some $y \in A^*$). Hence $p_1 = p_1'xy$, which 
implies $x = y = \varepsilon$ (the empty word), since $P_1$ is a prefix code.
Thus, $p_1 = p_2 \in P_2$. Therefore, $P_1 \subseteq P_2$. Similarly, 
$P_2 \subseteq P_1$, so $P_1 = P_2$. 

\medskip

\noindent
{\bf Proof of (2).} \ If $R$ is essential then $R$ intersects any right ideal, 
in particular $uA^*$. Thus, \ $R \cap uA^* \neq \emptyset$, hence $ux \in R$ 
for some $x \in A^*$. 

Conversely, consider any right ideal $PA^*$. For any $u \in P$, \ 
$ux \in R$ for some $x \in A^*$, hence $R$ intersects $uA^*$, hence $R$ 
intersects $PA^*$.  

\medskip

\noindent
{\bf Proof of (3).} \ Let $P$ be the prefix code corresponding to the right 
ideal $R$. By definition, $P$ is a maximal prefix code iff \ 
 $(\forall u \in A^* - P): P \cup \{u\}$ is not a prefix code. 

If $PA^*$ is essential then (by (2)) for any $u \in A^* - P$, \   
$uA^* \cap PA^* \neq \emptyset$,  hence \  
$(\exists p \in P)(\exists x, y \in A^*) \ ux = py$. 
Hence, $u \leq_{\rm pref} p$ or $p \leq_{\rm pref} u$. Therefore, since $u$ 
is prefix-comparable to some word in $P$, $P \cup \{u\}$ is not a prefix 
code. Thus, $P$ is maximal.

Conversely, if $P$ is maximal and $u \in A^*$, let us show that 
$ux \in PA^*$ for some $x \in A^*$.
If $u \in P$ this is obviously true (taking $x = \varepsilon$). 
If $u \notin P$, \ $P \cup \{u\}$ is not a prefix code, by maximality
of $P$, hence $u$ either has a prefix $p$ for some $p \in P$, hence 
$u = py  \in PA^*$ (for some $y \in A^*$), or $u$ is a prefix of some 
$p \in P$; then $ux = p \in PA^*$ (for some $x \in A^*$). 

\medskip

\noindent
{\bf Proof of (4).} \ Let $P$ be a maximal prefix code. By (1) -- (3) of this 
Lemma,  for any $u \in A^* - PA^*$ there is $x \in A^*$ such that $ux \in PA^*$.
Then, replacing $x$ by the shortest prefix $x'$ of $x$ for which $ux' \in PA^*$,
we have $ux' \in P$.

The converse is immediate form (2).  \ \ \ $\Box$

\begin{lem}
\label{bijPrefCode} \  
Let $\varphi: P_1A^* \to P_2A^*$ be a right-ideal isomorphism, where $P_1, P_2$
are prefix codes. Then $\varphi$ maps $P_1$ bijectively onto $P_2$.
\end{lem}      
{\bf Proof.} \ Since $\varphi$ is injective on $P_1A^*$ it is also injective 
on $P_1$. Let us show that $\varphi(P_1) \subseteq P_2$. 
If $p_1 \in P_1$ then $\varphi(p_1) = p_2 w$ for some 
$p_2 \in P_2, \ w \in A^*$.
Also, $p_1 = \varphi^{-1}(p_2 w) = \varphi^{-1}(p_2) \, w = p_1' vw$, for 
some $v \in A^*$. Since $P_1$ is a prefix code it follows that $p_1 = p_1'$ 
and hence $w= \varepsilon$. Thus, $\varphi(p_1) = p_2 \in P_2$.
A similar reasoning, applied to $\varphi^{-1}$, implies that 
$\varphi^{-1}(P_1) \subseteq P_2$, hence $P_1 \subseteq \varphi(P_2)$.  
  \ \ \ $\Box$

\begin{lem}
\label{IntersectTrans} \  
The intersection of two essential right ideals $R_1, R_2$ of $A^*$ is an 
essential right ideal of $A^*$. 
\end{lem}
{\bf Proof.} \ Let $R_1, R_2$ be two essential right ideals of $A^*$.
The intersection $R_1 \cap R_2$ is obviously a right ideal. 
By Lemma \ref{ideals} (2), for all $u \in A^*$ there is $x \in A^*$ such that 
$ux \in R_1$ (since $R_1$ is essential). Moreover, applying Lemma \ref{ideals} 
(2) to $R_2$: for $ux$ there is $y \in A^*$ such that $uxy \in R_2$. 
Since $R_1$ 
is a right ideal and since $ux \in R_1$, we also have $uxy \in R_1$. 
Thus, for all $u \in A^*$ there is $z = xy \in A^*$ such that 
$uxy \in R_1 \cap R_2$, which implies that $R_1 \cap R_2$ is essential.   
 \ \ \ $\Box$

\begin{lem}
\label{referee1} \ 
Let $\varphi: R \to Q$ be an isomorphism between essential right ideals, and
let $\varphi': R' \to Q'$ be a restriction of $\varphi$ to right subideals
$R' \subset R$, $Q' \subset Q$, with $\varphi'(R') = Q'$. Then we have:
 \ $R'$ is essential iff $Q'$ is essential. 
\end{lem}
{\bf Proof.} \ Assume $Q'$ is essential. Let $u$ be any word over $A^*$. Since
$R$ is essential, there exists $x \in A^*$ such that $ux \in R$; hence 
$\varphi(ux) \in Q$. Since $Q'$ is essential, there exists $y \in A^*$ such 
that $\varphi(ux) \, y \in Q'$; hence $\varphi^{-1}(\varphi(ux) \, y) \in R'$.
Moreover, $\varphi^{-1}(\varphi(ux) \, y) = ux \, \varphi^{-1}(y)$. Thus,
every word $u \in A^*$ is the prefix of a word in $R'$, which implies that
$R'$ is essential. 

The proof in the other direction is symmetric to this proof, since 
$\varphi$ is an isomorphism.
 \ \ \ $\Box$

\begin{lem}
\label{Q} \  
Assume $Q$ is a prefix code of $A^*$, and let $x\in A^*$. Then 
$\overline{x}Q = \{y\in A^* : xy \in Q\}$ is either the empty set or a prefix 
code. If $Q$ is a maximal prefix code of $A^*$ then $\overline{x}Q$ is 
maximal too.
\end{lem}
{\bf Proof.} \ If $y_1, y_2 \in \overline{x}Q$ are prefix-comparable then 
$xy_1, xy_2$ will also be prefix-comparable, contradicting the fact that 
$Q$ is a prefix code. 

To show maximality of $\overline{x}Q$ if $Q$ is maximal, consider any word 
$z \in A^*$; we want to show that $z$ is
prefix-comparable with some element of $\overline{x}Q$. 
Since $Q$ is a maximal prefix code, $xz$ is prefix-comparable with some 
$q \in Q$. \\  
$\bullet$ \ If $xz$ is a prefix of $q$ then $xzt = q$ for some $t \in A^*$, 
so $z$ is a prefix of $zt \in \overline{x}Q$.  \\   
$\bullet$ \ If $q$ is a prefix of $xz$ then no other element of $Q$ is a 
prefix of $xza$. Either (case 1), $q$ is a prefix of $x$ or (case 2), 
$q = xp$ for some $p \in \overline{x}Q$ such that  $p$ is a prefix of $z$; hence
in case 2, $z$ is prefix-comparable to an element of $\overline{x}Q$. In case 1,
$\overline{x}Q = \emptyset$; indeed, $q$ is a prefix of $x$ (so 
$\overline{x}\{q\} = \emptyset$), and no other element of $Q$ is 
prefix-comparable with $x$ ($Q$ being a prefix code), hence 
$\overline{x}\{q'\} = \emptyset$ for all $q' \in Q-\{q\}$. 
 \ \ \ $\Box$

\newpage


\section{Appendix A2}

In this appendix we give details of the proof of Lemma \ref{wLengthTaut}.
If $x, y \in \{a,b\}^*$, we abbreviate \ $x >_{{\rm pref}} y$ \ to \ $x > y$;
recall that this means that $x\{a,b\}^*$ strictly contains $y\{a,b\}^*$, i.e., 
$x$ is a strict prefix of $y$. Similarly, $x \not\geq y$
means that $x$ is not a prefix of $y$.

\bigskip

\noindent {\bf Fact A2.1} \ 

\smallskip

(1.1)  \ \ \ \ \ \ \ \ \ \ \
$\sigma^{-1} \cdot (a^k|b^h av) \cdot \sigma =  (a^{k+1}|b^{h-1} av)$ \ \
 for all $h \geq 2$, $v \in  \{a,b\}^*$, $k \geq 1$;

\medskip

(1.2) \ \ \ \ \ \ \ \ \ \ \
$\sigma^{-1} \cdot (a^k|bav) \cdot \sigma =  (a^{k+1}|abv)$ \ \
 for all $v \in  \{a,b\}^*$, $k \geq 1$;

\medskip

(1.3) \ \ \ \ \ \ \ \ \ \ \
$\sigma^{-1} \cdot (a^k|b^h) \cdot \sigma = (a^{k+1}|b^{h-1})$ \ \
 for all $h \geq 2$, $k \geq 1$;

\medskip

(1.4) \ \ \ \ \ \ \ \ \ \ \
$(ab|b) \cdot (a^k|b) \cdot (ab|b) = (a^k| ab)$ \ \
 for all $k \geq 2$.

\bigskip

\noindent {\bf Proof.} \ 
Verification of (1.1): \ When $h \geq 2$, $v \in  \{a,b\}^*$, 
$k \geq 1$, then 

\medskip

$\sigma^{-1} \cdot (a^k|b^h av) \cdot \sigma $

\medskip

$= \ \sigma^{-1} \cdot
\left[ \begin{array}{cccccc}
 \ a^k    \ & \ b^h av \ & \ a^ib \ & \ b^ja \ & \ b^{h+1} \ & \ b^h au\ell \  \\
 \ b^h av \ & \ a^k    \ & \ a^ib \ & \ b^ja \ & \ b^{h+1} \ & \ b^h au\ell \
\end{array}        \right]                                 \cdot  $
$\left[ \begin{array}{ccc}
   \ a^2 \ & \ ab \ & \ b    \\
   \ a   \ & \ ba \ & \ b^2
\end{array}        \right] $ 

 \hspace{1.85in} $_{1\leq i < k}$ \ $_{1\leq j < h}$
 \ \ \ \ \ \ \ \ \ \ $_{u > v, \ u\ell \not\geq v}$

\bigskip

$= \ \left[ \begin{array}{ccc}
         \ a   \ & \ ba \ & \ b^2 \\
         \ a^2 \ & \ ab \ & \ b
\end{array}        \right]              \cdot $
$\left[ \begin{array}{ccccccc}
 \ a^{k+1} \ & \ b^{h-1}av \ & \ a^ib \ & \ ab \ & \ b^ja \ & \ b^h \ & \ b^{h-1} au\ell \
\\
 \ b^h av \ & \ a^k     \ & \ a^{i-1}b\ & \ ba \ & \ b^{j+1}a\ & \ b^{h+1} \ & \ b^h au\ell
\
\end{array}        \right]   $

\makebox[3.1in]{} $_{2 \leq i \leq k}$ \ \ \ \ \ \ \ \ \ \ \ $_{1\leq j < h-1}$

\bigskip

$= \ \left[ \begin{array}{ccccccc}
 \ a^{k+1} \ & \ b^{h-1}av\ & \ a^ib \ & \ ab \ & \ b^ja \ & \ b^h \ & \ b^{h-1} au\ell \
\\
 \ b^{h-1}av\ & \ a^{k+1} \ & \ a^ib \ & \ ab \ & \ b^ja \ & \ b^h \ & \ b^{h-1} au\ell \
\end{array}        \right]   $

 \makebox[1.75in]{} $_{2 \leq i \leq k}$ \ \ \ \ \ \ \ \ $_{1\leq j < h-1}$

\bigskip

$= \ (a^{k+1} | b^{h-1}av)$.

\bigskip

\noindent 
$\bullet$ Verification of (1.2): \ When $v \in  \{a,b\}^*$, $k \geq 1$, 
then

\medskip

$\sigma^{-1} \cdot (a^k|b av) \cdot \sigma $

\medskip

$= \ \sigma^{-1} \cdot
\left[ \begin{array}{ccccc}
 \ a^k    \ & \ b av \ & \ a^ib \ & \ b^2 \ & \ b au\ell \  \\
 \ b av \  & \ a^k   \ & \ a^ib \ & \ b^2 \ & \ b au\ell \
\end{array}        \right]                       \cdot  $
$\left[ \begin{array}{ccc}
   \ a^2 \ & \ ab \ & \ b    \\
   \ a   \ & \ ba \ & \ b^2
\end{array}        \right] $

 \makebox[1.65in]{} $_{1\leq i < k}$ \ \ \ \ \ \ \  
             $_{u > v, \ u\ell \not\geq v}$

\bigskip

$= \ \left[ \begin{array}{ccc}
         \ a   \ & \ ba \ & \ b^2 \\
         \ a^2 \ & \ ab \ & \ b
\end{array}        \right]              \cdot $
$\left[ \begin{array}{ccccc}
 \ a^{k+1} \ & \ ab v \ & \ a^ib \ & \ b \ & \ ab u\ell \         \\
 \ b av \ & \ a^k    \ & \ a^{i-1}b\ & \ b^2 \ & \ ba u\ell \
\end{array}        \right]   $

 \makebox[2.85in]{} $_{2 \leq i \leq k}$

\bigskip

$= \ \left[ \begin{array}{ccccc}
 \ a^{k+1} \ & \ abv\ & \ a^ib \ & \ b \ & \ ab u\ell \ \\
 \ ab v\ & \ a^{k+1} \ & \ a^ib \ & \ b \ & \ ab u\ell \
\end{array}        \right]   $

 \makebox[1.45in]{} $_{2 \leq i \leq k}$

\bigskip

$= \ (a^{k+1} | abv)$.

\bigskip

\noindent $\bullet$ Verification of (1.3): \ When $h \geq 2$, $k \geq 1$,
then

\medskip

$\sigma^{-1} \cdot (a^k|b^h) \cdot \sigma$

\medskip

$=  \ \sigma^{-1} \cdot
\left[ \begin{array}{cccc}
 \ a^k \ & \ b^h \ & \ a^ib \ & \ b^ja \   \\
 \ b^h \ & \ a^k \ & \ a^ib \ & \ b^ja \
\end{array}        \right]                   \cdot  $
$\left[ \begin{array}{ccc}
   \ a^2 \ & \ ab \ & \ b    \\
   \ a   \ & \ ba \ & \ b^2
\end{array}        \right] $

 \makebox[1.5in]{} $_{1\leq i < k}$ \ $_{1 \leq j \leq h-1}$

\bigskip

$= \ \left[ \begin{array}{ccc}
         \ a   \ & \ ba \ & \ b^2 \\
         \ a^2 \ & \ ab \ & \ b
\end{array}        \right]              \cdot $
$\left[ \begin{array}{ccccc}
 \ a^{k+1} \ & \ b^{h-1} \ & \ a^ib \ & \ ab \ & \ b^ja \         \\
 \ b^h    \ & \ a^k    \ & \ a^{i-1}b\ & \ ba \ & \ b^{j+1}a \
\end{array}        \right]   $

 \makebox[2.9in]{} $_{2 \leq i \leq k}$ \ \ \ \ \ \ \ \ \ \ $_{1 \leq j \leq h-2}$

\bigskip

$= \ \left[ \begin{array}{ccccc}
 \ a^{k+1} \ & \ b^{h-1} \ & \ a^ib \ & \ ab \ & \ b^ja \ \\
 \ b^{h-1} \ & \ a^{k+1} \ & \ a^ib \ & \ ab \ & \ b^ja \
\end{array}        \right]   $

 \makebox[1.45in]{} $_{2 \leq i \leq k}$ \ \ \ \ \ \ $_{1 \leq j \leq h-2}$

\bigskip

$= \ (a^{k+1} | b^{h-1})$.

\bigskip

\noindent $\bullet$ Verification of (1.4): \ When $k \geq 2$, then 

\medskip

$(ab|b) \cdot (a^k|b) \cdot (ab|b)$

\medskip

$=  \ (ab|b) \cdot
\left[ \begin{array}{ccc}
 \ a^k \ & \ b   \ & \ a^ib \    \\
 \ b   \ & \ a^k \ & \ a^ib \
\end{array}        \right]               \cdot  $
$\left[ \begin{array}{ccc}
   \ a^2 \ & \ ab \ & \ b    \\
   \ a^2 \ & \ b \ & \ ab
\end{array}        \right] $

 \makebox[1.6in]{} $_{1\leq i \leq k-1}$

\bigskip

$= \ \left[ \begin{array}{ccc}
   \ a^2 \ & \ ab \ & \ b    \\
   \ a^2 \ & \ b \ & \ ab
\end{array}        \right]              \cdot $
$\left[ \begin{array}{cccc}
 \ a^k \ & \ ab  \ & \ a^ib \ & \ b \       \\
 \ b   \ & \ a^k \ & \ a^ib\ & \ ab \
\end{array}        \right]   $

 \makebox[2.6in]{} $_{2 \leq i \leq k-1}$

\bigskip

$= \ \left[ \begin{array}{cccc}
 \ a^k \ & \ ab  \ & \ a^ib \ & \ b \  \\
 \ ab  \ & \ a^k \ & \ a^ib \ & \ b \
\end{array}        \right]   $

 \makebox[1.1in]{} $_{2 \leq i \leq k-1}$

\bigskip

$ = (a^k| ab)$.

\bigskip

\noindent This proves Fact A2.1.  \ \ \ $\Box$


\bigskip

\bigskip

\noindent {\bf Fact A2.2} \

\smallskip

\noindent (2.1) \ \ \ \ \ \ \ \ \ \ \
If \ $(k>) \ j \geq 2$, $h \geq1$, \ then \ \ \ \
$\sigma \cdot (a^k|a^j b^h v) \cdot \sigma^{-1} = (a^{k-1} |a^{j-1} b^h v)$.

\medskip

\noindent (2.2) \ \ \ \ \ \ \ \ \ \ \
If \ $j =1$, $k \geq 3$, $h \geq1$, \ then \ \ \ \
$\delta \cdot (a^k|a b^h v) \cdot \delta^{-1} = (a^{k-1} |ab^h v)$.

\medskip

\noindent (2.3) \ \ \ \ \ \ \ \ \ \ \
If \ $j =1$, $k =2$, $h \geq 2$, \ then \ \ \ \
$\gamma_1 \cdot (a^2|a b^h v) \cdot \gamma_1^{-1} = (a^2|a b^{h-1} v)$.

\medskip

\noindent (2.4) \ \ \ \ \ \ \ \ \ \ \
If \ $j =1$, $k =2$, $h =1$, $|v| \geq 2$ (so $v = au$ for some
$u \in  \{a,b\} \, \{a,b\}^*$), 

\smallskip

 \hspace{2.4in} then \ \ \ \ 
   $\gamma_2 \cdot (a^2|a b au) \cdot \gamma_2^{-1} = (a^2|a b u)$.

\bigskip

\noindent {\bf Proof.} \ 
Verification of (2.1): \ If \ $(k>) \ j \geq 2$, $h \geq 1$, \ then 

\medskip

$\sigma \cdot (a^k|a^j b^h v) \cdot \sigma^{-1} $ \ \ \ \ \ \ \ \ \ \ \ \ \
 \ \ \ \ \ \ \ \ \  
  {\footnotesize (The column ``$a^jb^ra$'' is absent if $h = 1$.)}
                                                                   
\medskip

$= \ \sigma \cdot
\left[ \begin{array}{ccccc}
 \ a^k      \ & \ a^j b^h v \ & \ a^ib \ & \ a^jb^ra \ & \ a^jb^h p\ell \  \\
 \ a^j b^h v \ & \ a^k    \ & \ a^ib \ & \ a^jb^ra \ & \ a^jb^h p\ell \
\end{array}        \right]                                 \cdot  $
$\left[ \begin{array}{ccc}
         a   \ & \ ba \ & \ b^2 \\
         a^2 \ & \ ab \ & \ b
\end{array}        \right] $

 \makebox[1.85in]{} $_{0\leq i \leq k-1}$ \ $_{1\leq r < h}$
 \ \ \   $_{p > v, \ p\ell \not\geq v}$

 \makebox[1.85in]{} $_{i\neq j}$ 

\bigskip

$= \ \left[ \begin{array}{ccc}
   a^2 \ & \ ab \ & \ b    \\
   a   \ & \ ba \ & \ b^2
\end{array}        \right]   \cdot $
$ \left[ \begin{array}{ccccccc}
a^{k-1}\ &\ a^{j-1}b^hv\ &\ a^ib\ &\ a^{j-1}b^ra\ &\ a^{j-1}b^hp\ell\ &\ ba\ &\ bb \\ 
a^jb^hv\ & \ a^k \ &\ a^{i+1}b\ &\ a^jb^ra\ &\ a^jb^h p\ell\ &\ ab\ &\ b 
\end{array}        \right]       $

 \makebox[3.1in]{} $_{0\leq i \leq k-2}$ \ \ $_{1\leq r < h}$

 \makebox[3.2in]{} $_{i\neq j-1}$ 

\bigskip

\bigskip

$= \ \left[ \begin{array}{ccccccc}
a^{k-1}\ &\ a^{j-1}b^hv\ &\ a^ib\ &\ a^{j-1}b^ra\ &\ a^{j-1}b^hp\ell\ &\ ba\ &\ bb \\ 
a^{j-1}b^hv\ &\ a^{k-1}\ &\ a^ib\ &\ a^{j-1}b^ra\ &\ a^{j-1}b^hp\ell\ &\ ba\ &\ bb 
\end{array}        \right]       $

 \makebox[1.75in]{} $_{0\leq i \leq k-2}$ \ \ $_{1\leq r < h}$

 \makebox[1.8in]{} $_{i\neq j-1}$ 

\bigskip

$= \ (a^{k-1} |a^{j-1} b^h v)$.

\bigskip

\noindent $\bullet$ Verification of (2.2):  \ 
If \ $j =1$, $k \geq 3$, $h \geq 1$, \ then   

\medskip

$\delta \cdot (a^k|a b^h v) \cdot \delta^{-1}$
 \ \ \ \ \ \ \ \ \ \ \ \ \
 \ \ \ \ \ \ \ \ \
  {\footnotesize (The column ``$ab^ra$'' is absent if $h = 1$.)}

\medskip

$= \ \delta \cdot
\left[ \begin{array}{cccccc}
 \ a^k \ & \ ab^h v \ & \ a^ib \ & \ ab^ra \ & \ ab^h p\ell \ & \ b \ \\
 \ ab^hv \ & \ a^k  \ & \ a^ib \ & \ ab^ra \ & \ ab^h p\ell \ & \ b \  
\end{array}        \right]                                 \cdot  $
$\left[ \begin{array}{cccc}
   a^2 \ & \ ba \ & \ ab \ & \ b^2    \\
   a^3 \ & \ a^2b \ & \ ab \ & \ b
\end{array}        \right] $

 \makebox[1.5in]{} $_{2 \leq i \leq k-1}$ \ $_{1\leq r < h}$
 \ \ \   $_{p > v, \ p\ell \not\geq v}$

\bigskip

$= \ \left[ \begin{array}{cccc}
         a^3 \ & \ a^2b \ & \ ab \ & \ b  \\  
         a^2 \ & \ ba \ & \ ab \ & \ b^2    
\end{array}        \right]              \cdot $
$ \left[ \begin{array}{ccccccc}
a^{k-1}\ &\ ab^hv\ &\ a^ib\ &\ ab^ra\ &\ ab^hp\ell\ &\ ba\ &\ bb \\
ab^hv\ & \ a^k \ &\ a^{i+1}b\ &\ ab^ra\ &\ ab^h p\ell\ &\ a^2b\ &\ b
\end{array}        \right]       $

 \makebox[3.2in]{} $_{2\leq i \leq k-2}$ \ \ $_{1\leq r < h}$

\bigskip

$= \ \left[ \begin{array}{ccccccc}
a^{k-1}\ &\ ab^hv\ &\ a^ib\ &\ ab^ra\ &\ ab^hp\ell\ &\ ba\ &\ bb \\
ab^hv\ &\ a^{k-1}\ &\ a^ib\ &\ ab^ra\ &\ ab^hp\ell\ &\ ba\ &\ bb
\end{array}        \right]       $

 \makebox[1.3in]{} $_{2\leq i \leq k-2}$ \ $_{1\leq r < h}$

\medskip

$ = \ (a^{k-1} |ab^h v)$.

\bigskip

\noindent $\bullet$ Verification of (2.3): \
If \ $j =1$, $k =2$, $h \geq 2$, \ then \ \ \ \

\medskip

$\gamma_1 \cdot (a^2|a b^h v) \cdot \gamma_1^{-1} $.
 \ \ \ \ \ \ \ \ \ \ \ \ \ \ \ \ \ \
  {\footnotesize (The column ``$ab^h p\ell$'' is absent if $v$
                  is empty.)}

\medskip

$= \ \gamma_1 \cdot
\left[ \begin{array}{ccccc}
 a^2 \ & \ ab^h v \ & \ ab^ra \ & \ ab^h p\ell \ & \ b  \\
 ab^hv \ & \ a^2  \ & \ ab^ra \ & \ ab^h p\ell \ & \ b 
\end{array}        \right]                                 \cdot  $
$\left[ \begin{array}{cccc}
   a^2 \ & \ ba \ & \ ab \ & \ b^2    \\
   a^2 \ & \ aba \ & \ ab^2 \ & \ b
\end{array}        \right] $

 \makebox[1.7in]{} $_{1\leq r < h}$ \ \ $_{p > v, \ p\ell \not\geq v}$

\bigskip

$= \ \left[ \begin{array}{cccc}
 a^2 \ & \ aba \ & \ ab^2 \ & \ b        \\
 a^2 \ & \ ba  \ & \ ab   \ & \ b^2 
\end{array}        \right]              \cdot $
$ \left[ \begin{array}{cccccc}
a^2   \ &\ ab^{h-1}v \ &\ ab^ra\ &\ ab^hp\ell\ &\ ba\ &\ bb \\
ab^hv \ & \ a^2 \ & \ ab^{r+1}a\ &\ ab^{h+1}p\ell\ &\ aba\ &\ b
\end{array}        \right]       $

 \makebox[3.5in]{}  $_{1\leq r \leq h-2}$

\bigskip

$= \ \left[ \begin{array}{cccccc}
a^2 \ &\ ab^{h-1}v\ &\ ab^ra\ &\ ab^hp\ell\ &\ ba\ &\ bb \\
ab^{h-1}v\ &\ a^2 \ &\ ab^ra\ &\ ab^hp\ell\ &\ ba\ &\ bb
\end{array}        \right]       $

 \makebox[1.6in]{} $_{1\leq r \leq h-2}$

\medskip

$= \ (a^2|a b^{h-1} v).$

\bigskip

\noindent $\bullet$ Verification of (2.4): \  
If \ $j =1$, $k =2$, $h =1$, $|v| \geq 2$ (so $v = au$ for some 
$u \in \{a,b\} \, \{a,b\}^*$), \ then \ 
$\gamma_2 \cdot (a^2|a bau) \cdot \gamma_2^{-1} = (a^2|a b u)$.

\bigskip

Case (2.4.$a$): \ $v = aaw$, for some $w \in \{a,b\}^*$.

\medskip

$\gamma_2 \cdot (a^2|a b a^2 w) \cdot \gamma_2^{-1}$ 
 \ \ \ \ \ \ \ \ \ \ \ \ \ \ \ \ \ \ \ \ \ \ \ \ \ \ \ \ \
  {\footnotesize (The column ``$aba^2p\ell$'' is absent if $w$
                  is empty.)}

\medskip

$= \ \gamma_2 \cdot
\left[ \begin{array}{cccccc}
 a^2 \ & \ aba^2w \ & \ abab \ & \ ab^2 \ & \ aba^2p\ell \ & \ b  \\
 aba^2w \ & \ a^2 \ & \ abab \ & \ ab^2 \ & \ aba^2p\ell \ & \ b
\end{array}        \right]                                 \cdot  $
$\left[ \begin{array}{cccc}
   a^2 \ & \ ab \ & \ ba \ & \ b^2    \\
   a^2 \ & \ aba \ & \ ab^2 \ & \ b
\end{array}        \right] $

 \makebox[2.9in]{}  $_{p > w, \ p\ell \not\geq w}$

\bigskip

$= \ \left[ \begin{array}{cccc}
 a^2 \ & \ aba \ & \ ab^2 \ & \ b        \\
 a^2 \ & \ ab  \ & \ ba   \ & \ b^2
\end{array}        \right]              \cdot $
$ \left[ \begin{array}{cccccc}
a^2   \ &\ abaw \ &\ ab^2\ &\ abap\ell\ &  \ ba\ &\ bb \\
aba^2w \ & \ a^2 \ & \ abab\ &\ aba^2p\ell\ &\ ab^2\ &\ b
\end{array}        \right]       $

\medskip

$= \ \left[ \begin{array}{cccccc}
a^2 \ &\ abaw\ &\ ab^2\ &\ abap\ell\ &\ ba\ &\ bb \\
abaw \ &\ a^2 \ &\ ab^2\ &\ abap\ell\ &\ ba\ &\ bb
\end{array}        \right]       $

\medskip

$= \ (a^2|a ba w)$.

\bigskip

Case (2.4.$b$): \ $v = abw$, for some $w \in \{a,b\}^*$.

\medskip

$\gamma_2 \cdot (a^2|a b abw) \cdot \gamma_2^{-1}$
 \ \ \ \ \ \ \ \ \ \ \ \ \ \ \ \ \ \ \ \ \ \ \ \ \ \ \ \ \
  {\footnotesize (The column ``$ababp\ell$'' is absent if $w$
                  is empty.)}

\medskip

$= \ \gamma_2 \cdot
\left[ \begin{array}{cccccc}
 a^2 \ & \ ababw \ & \ aba^2 \ & \ ab^2 \ & \ ababp\ell \ & \ b  \\
 ababw \ & \ a^2 \ & \ aba^2 \ & \ ab^2 \ & \ ababp\ell \ & \ b
\end{array}        \right]                                 \cdot  $
$\left[ \begin{array}{cccc}
   a^2 \ & \ ab \ & \ ba \ & \ b^2    \\
   a^2 \ & \ aba \ & \ ab^2 \ & \ b
\end{array}        \right] $

 \makebox[2.9in]{}  $_{p > w, \ p\ell \not\geq w}$

\bigskip

$= \ \left[ \begin{array}{cccc}
 a^2 \ & \ aba \ & \ ab^2 \ & \ b        \\
 a^2 \ & \ ab  \ & \ ba   \ & \ b^2
\end{array}        \right]              \cdot $
$ \left[ \begin{array}{cccccc}
a^2   \ &\ ab^2w \ &\ aba\ &\ ab^2p\ell\ &  \ ba\ &\ bb \\
ababw \ & \ a^2 \ & \ aba^2\ &\ ababp\ell\ &\ ab^2\ &\ b
\end{array}        \right]       $

\medskip

$= \ \left[ \begin{array}{cccccc}
a^2 \ &\ abbw\ &\ aba\ &\ ab^2p\ell\ &\ ba\ &\ bb \\
abbw \ &\ a^2 \ &\ aba\ &\ ab^2p\ell\ &\ ba\ &\ bb
\end{array}        \right]       $

\medskip

$= \ (a^2|a b bw)$.

\bigskip

\noindent This proves Fact A2.2.  \ \ \ $\Box$

\newpage 


\section{Appendix A3}

In this appendix we give details of the proof of Theorem \ref{linearDist}. 

\bigskip

{\sc Case} $\alpha \alpha$: \ \ \ \ \
$\mu = \ \ \ \ \ \alpha^{h_m} \beta^{k_{m-1}} \ldots
        \beta^{k_i} \alpha^{h_i} \ldots \alpha^{h_2} \beta^{k_1} \alpha^{h_1}$,

\smallskip

{\sc Case} $\beta \alpha$: \ \ \ \ \
$\mu = \beta^{k_m} \alpha^{h_m} \beta^{k_{m-1}} \ldots
      \beta^{k_i} \alpha^{h_i} \ldots \alpha^{h_2} \beta^{k_1} \alpha^{h_1}$,

\bigskip

\noindent {\bf Claim ($*\alpha$).}  \ {\it
Let $\mu$ be as above, according to cases $\alpha \alpha$ or $\beta \alpha$.
Then,

\medskip

$\mu(a) =  (w_m v_m^{|k_m|-1}) \ t_m u_m^{|h_m|-1} w_{m-1} \ldots $
             $ w_i v_i^{|k_i|-1} t_i u_i^{|h_i|-1} \ldots $
                         $ w_1 v_1^{|k_1|-1} t_1 \ u_1^{|h_1|-1} w_0$

\medskip

\noindent where:

\medskip

$ w_0 = \left\{ \begin{array}{ll}
                ba  & \mbox{if \ \ $h_1> 0$} \\
                a^2 & \mbox{if \ \ $h_1 < 0$} \end{array}  \right. $

\medskip

\noindent For \ $m-1 \geq i \geq 1$: \ \ \ \
$ w_i = \left\{ \begin{array}{ll}
                ba^3  & \mbox{if \ \ $h_{i+1} > 0$, $k_i > 0$} \\
                ba^2b & \mbox{if \ \ $h_{i+1} > 0$, $k_i < 0$} \\
                a^4   & \mbox{if \ \ $h_{i+1} < 0$, $k_i > 0$} \\
                a^3b  & \mbox{if \ \ $h_{i+1} < 0$, $k_i < 0$}
                \end{array}  \right. $

\medskip

\noindent In case $\beta \alpha$: \ \ \ \
$ w_m = \left\{ \begin{array}{ll}
                a^3  & \mbox{if \ \ $k_m > 0$} \\
                a^2b & \mbox{if \ \ $k_m < 0$} \end{array}  \right. $

\medskip

\noindent For \ $m-1 \geq i \geq 1$,  or $i = m$ in case $\beta \alpha$:
 \ \ $ t_i = \left\{ \begin{array}{ll}
                bab^2  & \mbox{if \ \ $h_i > 0$, $k_i > 0$} \\
                a^2b^2 & \mbox{if \ \ $h_i > 0$, $k_i < 0$} \\
                baba   & \mbox{if \ \ $h_i < 0$, $k_i > 0$} \\
                a^2ba  & \mbox{if \ \ $h_i < 0$, $k_i < 0$}
                \end{array}  \right. $

\medskip

\noindent In case $\alpha \alpha$: \ \ \ \
The factor $(w_m v_m^{|k_m|-1})$ is absent, and \
 $t_m = \left\{ \begin{array}{ll}
                b^3  & \mbox{if \ \ $h_m > 0$} \\
                b^2a & \mbox{if \ \ $h_m < 0$} \end{array}  \right. $

\medskip

\noindent For $m \geq i \geq 1$: \ \ \ \ $ u_i = \left\{ \begin{array}{ll}
                a & \mbox{if \ \ $h_i > 0$} \\
                b & \mbox{if \ \ $h_i < 0$}
                \end{array}  \right. $   \ \ \ \ \  \ \ \
$ v_i = \left\{ \begin{array}{ll}
                a & \mbox{if \ \ $k_i > 0$} \\
                b & \mbox{if \ \ $k_i < 0$}
                \end{array}  \right. $

\bigskip

\noindent We also have:

\smallskip

$\mu(ba) = (w_m v_m^{|k_m|-1}) \ t_m u_m^{|h_m|-1} w_{m-1} \ldots $
             $ w_i v_i^{|k_i|-1} t_i u_i^{|h_i|-1} \ldots $
                         $ w_1 v_1^{|k_1|-1} t_1 u_1^{|h_1|-1} W_0$

\smallskip

\noindent where:

\medskip

$W_0 = \left\{ \begin{array}{ll}
                b^2  & \mbox{if \ \ $h_1> 0$} \\
                ab & \mbox{if \ \ $h_1 < 0$} \end{array}  \right. $

\medskip

\noindent and all other $w_i$, $v_i$, $t_i$, and $u_i$ are the same as
for $\mu(a)$.
} 

\medskip

\noindent {\bf Proof of Claim ($* \alpha$).} \ The proof goes by induction on
the number of exponents $(k_m,) \ h_m, \ldots, k_1, h_1$
$(\in \mathbb{Z} -\{0\})$. 

\medskip

\noindent {\sc Base of the Induction:} \  A straightforward induction on
$h_1$ shows that for all $h_1 \in \mathbb{Z} -\{0\}$,

\smallskip

$\alpha^{h_1}(a) = \left\{ \begin{array}{ll}
           b^3 a^{|h_1|-1} ba  & \mbox{if \ \ $h_1 > 0$} \\
           b^2a b^{|h_1|-1}a^2 & \mbox{if \ \ $h_1 < 0$} \end{array}  \right. $

\smallskip

$\alpha^{h_1}(ba) = \left\{ \begin{array}{ll}
           b^3 a^{|h_1|-1} b^2 & \mbox{if \ \ $h_1 > 0$} \\
           b^2a b^{|h_1|-1}ab  & \mbox{if \ \ $h_1 < 0$} \end{array}  \right. $

\smallskip

\noindent So Claim $(*\alpha)$ holds when we just have one non-zero
exponent.

\medskip

\noindent {\sc Inductive Step}: \ We assume $h_1 \neq 0$.

\smallskip

\noindent {\bf Case} $\alpha \alpha$: We consider the case where $\mu$ is of
the form $\alpha \alpha$ (with, $k_m = 0, \ h_m \neq 0$).

\smallskip

Assume $h_m, \ldots, k_1, h_1$ are non-zero, and $k_m = 0$. By induction,
$\mu(a)$ and $\mu(ba)$ are of the form

\smallskip

$\mu(a) = t_m u_m^{|h_m|-1} w_{m-1} \ldots $
                         $ w_1 v_1^{|k_1|-1} t_1 \ u_1^{|h_1|-1} w_0$

\smallskip

$\mu(ba) = t_m u_m^{|h_m|-1} w_{m-1} \ldots $
                         $ w_1 v_1^{|k_1|-1} t_1 u_1^{|h_1|-1} W_0$
\smallskip

\noindent with \  $t_m = \left\{ \begin{array}{ll}
                b^3  & \mbox{if \ \ $h_m > 0$} \\
                b^2a & \mbox{if \ \ $h_m < 0$} \end{array}  \right. $

\medskip

\noindent {\bf (1)} \ If $h_m > 0$ and another $\alpha$ is applied to $\mu(a)$
or $\mu(ba)$, then $t_m = b^3$; by looking at the entry $b^3$ in the table of
$\alpha$ we obtain \

\smallskip

 $\alpha \mu(a) = t_m u_m^{|h_m|+1-1} w_{m-1} \ldots $
     $w_1 v_1^{|k_1|-1} t_1 \ u_1^{|h_1|-1} w_0$ \ \ 
  (respectively $W_0$ instead of $w_0$ for $\alpha \mu(ba)$).

\smallskip

\noindent Also, $|h_m|+1-1 = |h_m +1| -1$ when $h_m > 0$, so we have verified
the induction hypothesis.

\smallskip

\noindent {\bf (2)} \ If $h_m < 0$ and another $\alpha^{-1}$ is applied to
$\mu(a)$ or to $\mu(ba)$, then $t_m = b^2a$; the entry $b^2a$ in the
range-row of the table of $\alpha$ yields then

\smallskip

 $\alpha^{-1} \mu(a) = t_m u_m^{|h_m|+1-1} w_{m-1} \ldots $
                         $ w_1 v_1^{|k_1|-1} t_1 \ u_1^{|h_1|-1} w_0$
(respectively $W_0$ instead of $w_0$ for $\alpha^{-1} \mu(ba)$).

\smallskip

\noindent When $h_m < 0$, $|h_m|+1-1 = |h_m -1|-1$, so we have verified
the induction hypothesis.

\smallskip

\noindent {\bf (3)} \ If $\beta$ is applied to $\mu(a)$ or to $\mu(ba)$,
and $h_m > 0$, we use the fact that $t_m$ $(= b^3)$ starts with $b$.
The entry $b$ in the table of $\beta$ implies that the leftmost $b$ in
$t_m$ is replaced by $a^3ba$.  So, $\mu(a)$ (or $\mu(ba)$) becomes

\smallskip

 $\beta(\mu(a))$ \ (or $\beta(\mu(ba)$)

\smallskip

 $= a^3 ba \, b^2 \, u_m^{|h_m|-1} w_{m-1} \ldots $
                         $ w_1 v_1^{|k_1|-1} t_1 \ u_1^{|h_1|-1} w_0$
       \ \ (respectively $W_0$).

\smallskip

 $= w'_m v_{m}^{|k_m|-1} t'_m u_m^{|h_m|-1} w_{m-1} \ldots $
                         $ w_1 v_1^{|k_1|-1} t_1 \ u_1^{|h_1|-1} w_0$
      \ \ (respectively $W_0$).

\smallskip

\noindent Indeed, now $k_m$ becomes $1$; when $k_m > 0$ we have for the new
factors $w'_m = a^3$ (and $v_m = b$, buy $|k_m|-1 = 0$ here); moreover,
when $k_m > 0$ and $h_m > 0$, the new value of ``$t_m$'' is  $t'_m = ab^2$.
So we have verified the induction hypothesis.

\smallskip

\noindent If $h_m < 0$, we obtain

\smallskip

 $\beta(\mu(a))$ \ (or $\beta(\mu(ba)$)

\smallskip

 $= a^3 ba \, ba \, u_m^{|h_m|-1} w_{m-1} \ldots $
                         $ w_1 v_1^{|k_1|-1} t_1 \ u_1^{|h_1|-1} w_0$
       \ \ (respectively $W_0$).

\smallskip

 $= w'_m v_m^{|k_m|-1} t'_m u_m^{|h_m|-1} w_{m-1} \ldots $
                         $ w_1 v_1^{|k_1|-1} t_1 \ u_1^{|h_1|-1} w_0$
       \ \ (respectively $W_0$).

\smallskip

\noindent Indeed, now $k_m = 1$; when $k_m > 0$ we have $w'_m = a^3$ (and
$v_m = b$, but $|k_m|-1 = 0$ here);
moreover, when $k_m > 0$ and $h_m < 0$, $t'_m = baba$.
So we have verified the induction hypothesis.

\smallskip

\noindent {\bf (4)} \ If $\beta^{-1}$  is applied to $\mu(a)$ or to $\mu(ba)$,
we again use the fact that $t_m$ starts with $b$. The entry $b$ in the
range-row of the table of $\beta$ implies that the leftmost $b$ in $t_m$
is replaced by $a^2ba^2$.

\smallskip

\noindent If $h_m > 0$, $\mu(a)$ (or $\mu(ba)$) becomes

\smallskip

 $\beta^{-1}(\mu(a))$ \ (or $\beta^{-1}(\mu(ba)$)

\smallskip

 $= a^2ba^2 \, b^2 \, u_m^{|h_m|-1} w_{m-1} \ldots $
                         $ w_1 v_1^{|k_1|-1} t_1 \ u_1^{|h_1|-1} w_0$
       \ \ (respectively $W_0$).

\smallskip

 $= w'_m v_m^{|k_m|-1} t'_m u_m^{|h_m|-1} w_{m-1} \ldots $
                         $ w_1 v_1^{|k_1|-1} t_1 \ u_1^{|h_1|-1} w_0$
       \ \ (respectively $W_0$).

\smallskip

\noindent Indeed, now, $k_m = -1$, and when $k_m > 0$ we have $w'_m = a^2b$
(and $v_m = a$, but $|k_m|-1 = 0$ here);
moreover, when $k_m < 0$  and $h_m > 0$, $t'_m = a^2b^2$.
So we have verified the induction hypothesis.

\smallskip

\noindent If $h_m < 0$, $\mu(a)$ (or $\mu(ba)$) becomes

\smallskip

 $\beta^{-1}(\mu(a))$ \ (or $\beta^{-1}(\mu(ba)$)

\smallskip

 $= a^2ba^2 \, ba \, u_m^{|h_m|-1} w_{m-1} \ldots $
                         $ w_1 v_1^{|k_1|-1} t_1 \ u_1^{|h_1|-1} w_0$
       (respectively $W_0$).

\smallskip

 $= w'_m v_m^{|k_m|-1} t'_m u_m^{|h_m|-1} w_{m-1} \ldots $
                         $ w_1 v_1^{|k_1|-1} t_1 \ u_1^{|h_1|-1} w_0$
       (respectively $W_0$).

\smallskip

\noindent Indeed, now, $k_m = -1$, and when $k_m > 0$ we have $w'_m = a^2b$
(and $v_m = a$, but $|k_m|-1 = 0$ here);
moreover, when $k_m < 0$  and $h_m < 0$, $t'_m = a^2ba$.

This completes the verification of the induction hypothesis in the case where
$\mu$ is of the form $\alpha \alpha$.

\medskip

\noindent {\bf Case} $\beta \alpha$: \ We consider the case where $\mu$ is of
the form $\beta \alpha$ (with $h_{m+1} = 0, k_m \neq 0$).
By induction, $\mu(a)$ and $\mu(ba)$ are of the form

\smallskip

$\mu(a) = w_m v_m^{|k_m|-1} t_m u_m^{|h_m|-1} w_{m-1} \ldots $
                         $ w_1 v_1^{|k_1|-1} t_1 \ u_1^{|h_1|-1} w_0$

\smallskip

$\mu(ba) = w_m v_m^{|k_m|-1} t_m u_m^{|h_m|-1} w_{m-1} \ldots $
                         $ w_1 v_1^{|k_1|-1} t_1 u_1^{|h_1|-1} W_0$

\smallskip

\noindent with $w_m = \left\{ \begin{array}{ll}
                a^3  & \mbox{if \ \ $k_m > 0$} \\
                a^2b & \mbox{if \ \ $k_m < 0$} \end{array}  \right.  $

\smallskip

\noindent {\bf (5)} \ If $\alpha$ is applied to the string $\mu(a)$
or $\mu(ba)$, we look at the entry $a$ in the table of $\alpha$
(since in all case, $w_m$ starts with $a$).

\smallskip

If $k_m > 0$, we use the entry $a$ in the table of $\alpha$ (since
$w_m =a^3$), and $\mu(a)$ (or $\mu(ba)$) becomes

\smallskip

$\alpha(\mu(a))$ (or $\alpha(\mu(ba)$)

\smallskip

$= b^4a \, a^2 \, v_m^{|k_m|-1} t_m u_m^{|h_m|-1} w_{m-1} \ldots $
                         $ w_1 v_1^{|k_1|-1} t_1 u_1^{|h_1|-1} w_0$
   \ \ (respectively $W_0$)

\smallskip

$=t'_{m+1} u_{m+1}^{|h_{m+1}|-1} w'_m v_m^{|k_m|-1} t_m u_m^{|h_m|-1} w_{m-1}$
  $ \ldots w_1 v_1^{|k_1|-1} t_1 u_1^{|h_1|-1} w_0$
   \ \ (respectively $W_0$)

\smallskip

\noindent Indeed, now $h_{m+1} = 1$; when $h_{m+1} > 0$ we have $t'_{m+1} = b^3$
(and $u_{m+1} = a$, but $|h_{m+1}|-1 = 0$ here anyway); moreover, when
$h_{m+1} > 0$ and $k_m > 0$, $w'_m = ba^3$.
So we have verified the induction hypothesis.

\smallskip

If $k_m < 0$, we use again the entry $a$ in the table of $\alpha$
(now, $w_m = a^2b$), and $\mu(a)$ (or $\mu(ba)$) becomes

\smallskip

$\alpha(\mu(a))$ (or $\alpha(\mu(ba)$)

\smallskip

$= b^4a \, ab \, v_m^{|k_m|-1} t_m u_m^{|h_m|-1} w_{m-1} \ldots $
                         $ w_1 v_1^{|k_1|-1} t_1 u_1^{|h_1|-1} w_0$
   \ \ (respectively $W_0$)

\smallskip

$= t'_{m+1} u_{m+1}^{|h_{m+1}|-1} w'_m v_m^{|k_m|-1} t_m u_m^{|h_m|-1} w_{m-1}$
  $ \ldots w_1 v_1^{|k_1|-1} t_1 u_1^{|h_1|-1} w_0$
   \ \ (respectively $W_0$)

\smallskip

\noindent Indeed, now $h_{m+1} = 1$; when $h_{m+1} > 0$ we have $t'_{m+1} = b^3$
(and $u_m = a$, but $|h_{m+1}|-1 = 0$ here anyway); moreover, when $h_{m+1} > 0$
and $k_m < 0$, $w'_m = ba^2b$.
So we have verified the induction hypothesis.

\smallskip

\noindent {\bf (6)} \ If $\alpha^{-1}$ is applied to the string $\mu(a)$
or $\mu(ba)$, we look at the entry $a$ in the table of $\alpha$
(since in all case, $w_m$ starts with $a$).

\smallskip

If $k_m > 0$, we use the entry $a$ in the range-row of the table of $\alpha$
(since $w_m =a^3$), and $\mu(a)$ (or $\mu(ba)$) becomes

\smallskip

$\alpha^{-1}(\mu(a))$ (or $\alpha^{-1}(\mu(ba)$)

\smallskip

$= b^2a^3 \, a^2 \, v_m^{|k_m|-1} t_m u_m^{|h_m|-1} w_{m-1} \ldots $
                         $ w_1 v_1^{|k_1|-1} t_1 u_1^{|h_1|-1} w_0$
   \ \ (respectively $W_0$)

\smallskip

$= t'_{m+1} u_{m+1}^{|h_{m+1}|-1} w'_m v_m^{|k_m|-1} t_m u_m^{|h_m|-1} w_{m-1}$
  $ \ldots w_1 v_1^{|k_1|-1} t_1 u_1^{|h_1|-1} w_0$
   \ \ (respectively $W_0$)

\smallskip

\noindent Indeed, now $h_{m+1} = -1$; when $h_{m+1} < 0$ we have 
$t'_{m+1} = b^2a$ (and $u_{m+1} = b$, but $|h_{m+1}|-1 = 0$ here); 
moreover, when $h_{m+1} < 0$ and $k_m > 0$, $w'_m = a^4$.
So we have verified the induction hypothesis.

\smallskip

If $k_m < 0$, we use again the entry $a$ in the table of $\alpha$
(now, $w_m = a^2b$), and $\mu(a)$ (or $\mu(ba)$) becomes

\smallskip

$\alpha^{-1}(\mu(a))$ (or $\alpha^{-1}(\mu(ba)$)

$= b^2a^3 \, ab \, v_m^{|k_m|-1} t_m u_m^{|h_m|-1} w_{m-1} \ldots $
                         $ w_1 v_1^{|k_1|-1} t_1 u_1^{|h_1|-1} w_0$
   \ \ (respectively $W_0$)

\smallskip

$= t'_{m+1} u_{m+1}^{|h_{m+1}|-1} w'_m v_m^{|k_m|-1} t_m u_m^{|h_m|-1} w_{m-1}$
  $ \ldots w_1 v_1^{|k_1|-1} t_1 u_1^{|h_1|-1} w_0$
   \ \ (respectively $W_0$)

\smallskip

\noindent Indeed, now $h_{m+1} = -1$; when $h_{m+1} < 0$ we have 
$t'_{m+1} = b^2a$ (and $u_m = b$, but $|h_{m+1}|-1 = 0$ here anyway); 
moreover, when $h_{m+1} < 0$ and $k_m < 0$, $w'_m = a^3b$.
So we have verified the induction hypothesis.

\smallskip

\noindent {\bf (7)} \ If $k_m > 0$ and another $\beta$ is applied to $\mu(a)$
or $\mu(ba)$, we look at the entry $a^3$ in the table of $\beta$
(since $w_m = a^3$). Then and $\mu(a)$ (or $\mu(ba)$) becomes

\smallskip

$\beta(\mu(a))$ (or $\beta(\mu(ba)$)

\smallskip

$= w_m v_m^{|k_m|+1 -1} t_m u_m^{|h_m|-1} w_{m-1}$
  $ \ldots w_1 v_1^{|k_1|-1} t_1 u_1^{|h_1|-1} w_0$
   \ \ (respectively $W_0$)

\smallskip

\noindent When $k_m > 0$, $v_m = a$ and $|k_m|+1 -1 = |k_m+1|-1$, so we have 
verified the induction hypothesis.

\smallskip

\noindent {\bf (8)} \ If $k_m < 0$ and another $\beta^{-1}$ is applied to 
$\mu(a)$ or $\mu(ba)$, we look at the entry $a^2b$ in the range-row of the 
table of $\beta$ (since $w_m = a^2b$). Then and $\mu(a)$ (or $\mu(ba)$) becomes

\smallskip

$\beta^{-1}(\mu(a))$ (or $beta^{-1}(\mu(ba)$)

\smallskip

$= w_m v_m^{|k_m|-1 -1} t_m u_m^{|h_m|-1} w_{m-1}$
  $ \ldots w_1 v_1^{|k_1|-1} t_1 u_1^{|h_1|-1} w_0$
   \ \ (respectively $W_0$)

\smallskip

\noindent When $k_m < 0$, $v_m = b$ and $|k_m|-1 -1 = |k_m+1|-1$, so we have 
verified the induction hypothesis.

\smallskip

This completes the proof of Claim  $(*\alpha)$.  
 \ \ \ $\Box$

\bigskip

\bigskip

{\sc Case} $\alpha \beta$: \ \ \ \ \
$\mu = \ \ \ \ \ \alpha^{h_m} \beta^{k_{m-1}} \ldots
      \beta^{k_i} \alpha^{h_i} \ldots \alpha^{h_2} \beta^{k_1}$,

\smallskip

{\sc Case} $\beta \beta$: \ \ \ \ \
$\mu = \beta^{k_m} \alpha^{h_m} \beta^{k_{m-1}} \ldots
       \beta^{k_i} \alpha^{h_i} \ldots \alpha^{h_2} \beta^{k_1}$.

\bigskip

\noindent {\bf Claim ($*\beta$).}  \ {\it
Let $\mu$ be as above, according to cases $\alpha \beta$ or $\beta \beta$.
Then,

\smallskip

$\mu(b) =  (w_m v_m^{|k_m|-1}) \ t_m u_m^{|h_m|-1} w_{m-1} \ldots $
             $ w_i v_i^{|k_i|-1} t_i u_i^{|h_i|-1} \ldots $
                         $ w_1 v_1^{|k_1|-1} t_1$

\smallskip

\noindent where:

\medskip

\noindent $t_1 = \left\{ \begin{array}{ll}
                  ba  & \mbox{if \ \ $k_1> 0$} \\
                  a^2 & \mbox{if \ \ $k_1 < 0$} \end{array}  \right. $

\medskip

\noindent and all other $w_i$, $v_i$, $t_i$, and $u_i$ are the same as
for $\mu(a)$ in Claim $(*\alpha)$.

\medskip

\noindent We also have:

\smallskip

$\mu(ab) = (w_m v_m^{|k_m|-1}) \ t_m u_m^{|h_m|-1} w_{m-1} \ldots $
             $ w_i v_i^{|k_i|-1} t_i u_i^{|h_i|-1} \ldots $
                         $ w_1 v_1^{|k_1|-1} T_1$

\smallskip

\noindent where:

\medskip

$T_1 = \left\{ \begin{array}{ll}
                  b^2  & \mbox{if \ \ $k_1> 0$} \\
                  ab   & \mbox{if \ \ $k_1 < 0$} \end{array}  \right. $

\medskip

\noindent and all other $w_i$, $v_i$, $t_i$, and $u_i$ are the same as
for $\mu(b)$.
} 

\medskip

\noindent {\bf Proof of Claim ($* \beta$).} \ The proof goes by induction on
the number of non-zero exponents $k_m, h_m, \ldots, k_1$ in $\mu$. By
assumption, in Claim $(* \beta)$ we have $k_1 \neq 0$.

\medskip

\noindent {\sc Base of the Induction:} \  A straightforward induction on
$k_1$ shows that for all $k_1 \in \mathbb{Z} -\{0\}$,

\smallskip

$\beta(b) = \left\{ \begin{array}{ll}
           a^3 a^{|k_1|-1} ba  & \mbox{if \ \ $k_1 > 0$} \\
           a^2b b^{|k_1|-1}a^2 & \mbox{if \ \ $k_1 < 0$} \end{array}  \right. $

\smallskip

$\beta(ab) = \left\{ \begin{array}{ll}
           a^3 a^{|k_1|-1} b^2 & \mbox{if \ \ $k_1 > 0$} \\
           a^2b b^{|k_1|-1}ab  & \mbox{if \ \ $k_1 < 0$} \end{array}  \right. $

\medskip

\noindent So Claim $(*\alpha)$ holds when we just have one non-zero
exponent.

\medskip

\noindent {\sc Inductive Step.}

\smallskip

The proof of the inductive step is very similar to the proof of the inductive
step of Claim $(*\alpha)$. The results are the same too, except for $t_1$ and
$T_1$, which are dealt with in the base of the induction.

The proof of Claim $(*\beta)$ completes the proof of Theorem \ref{linearDist}.
 \ \ \ $\Box$


\section{Appendix A4}

In this appendix we give the proofs of the Lemmas related to Theorem
\ref{thmThompson}, about the representation of Thompson groups in 
algebras.

\begin{lem}
\label{BclosedMult} \ \
${\cal B}_{\infty}$ is closed under multiplication.
\end{lem}
{\bf Proof.} \
The condition that $\{ (y_i, x_i) : i \in I\}$ is finite-to-finite will
guarantee that multiplication in ${\cal B}_{\infty}$ is well defined
(i.e., no infinite sums in $\mathbb{K}$ are used). Indeed, let \
$\sum_{j \in J} \alpha_j \, v_j u_j^{-1}$,
$\sum_{k \in K} \beta_k \, t_k s_k^{-1} \ \in \ {\cal B}_{\infty}$,
 \ and let

\medskip

$\sum_{i \in I} \kappa_i \, y_i x_i^{-1} \ =$ \
$\sum_{j \in J} \alpha_j \, v_j u_j^{-1} \ \cdot$ \
$\sum_{k \in K} \beta_k \, t_k s_k^{-1}$.

\medskip

\noindent Then, for any fixed $i \in I$ the coefficient of $y_i x_i^{-1}$
is

\smallskip

$\kappa_i \ = \ \sum \{ \alpha_j \, \beta_k \ : \ j \in J, k \in K$
 are such that \ $y_i x_i^{-1} = v_j u_j^{-1} \ t_k s_k^{-1} \}$.

\noindent If \ $y_i x_i^{-1} = v_j u_j^{-1} \ t_k s_k^{-1}$ \  then
we have the following two possibilities:

\smallskip

[case 1] \ \ \ $y_i x_i^{-1} = v_j r_k s_k^{-1}$, \ if \
                                 $t_k = u_j r_k \leq_{{\rm pref}} u_j$,
                 for some $r_k \in A^*$, \ \ \  or

\smallskip

[case 2] \ \ \ $y_i x_i^{-1} = v_j (s_k w_j)^{-1}$,  \ if \
                                 $t_k \geq_{{\rm pref}} u_j = t_k w_j$, \
                for some $w_j \in A^*$.

\smallskip

\noindent In both cases, there are only finitely many ways to choose $v_j$
for a given $y_i$ (since $v_j$ is a prefix of $y_i$). Hence, there are only
finitely many ways to choose $u_j$ (by the finite-to-finite property of the
sum $\sum \alpha_j \, v_ju_j^{-1}$).

In case 1, $r_k$ can be chosen in a finite number of ways (being a suffix
of $y_i$), hence there are only finitely many choices for $t_k$. Hence,
there are only finitely many choices for $s_k$ (by the finite-to-finite
property of the sum $\sum \beta_k \, t_ks_k^{-1}$).

In case 2, $t_k$ is a prefix of $u_j$, and there were only finitely many
possible choices for $t_k$ (since there are only finitely many choices for
$u_j$ in both cases, as we saw). By the the finite-to-finite condition
of the sum $\sum \beta_k \, t_ks_k^{-1}$,
there will only be finitely many choices for $s_k$.

In summary, if we fix just $y_i$ (irrespective of what $x_i$ might be),
$y_i x_i^{-1}$ has only a finite number of factorizations of
the form \ $y_i x_i^{-1} = v_j u_j^{-1} \ t_k s_k^{-1}$ \
(for a fixed \ $\sum_{j \in J} \alpha_j \, v_j u_j^{-1}$ and
$\sum_{k \in K} \beta_k \, t_k s_k^{-1} \ \in \ {\cal B}_{\infty}$).
Therefore, \ $\kappa_i = \sum \{ \alpha_j \, \beta_k \ : \ j, k$ etc. $\}$ \
is a finite sum, hence $\kappa_i$ is well defined.

\smallskip

The above also implies that  $y_i$ determines a finite number of
possibilities for $x$ such that \
$y_ix^{-1} \in \{y_i x_i^{-1} : i \in I\}$ \
(for a fixed \ $\sum_j \alpha_j \, v_j u_j^{-1}$,
$\sum_k \beta_k \, t_k s_k^{-1} \ \in \ {\cal B}_{\infty}$).
Hence, $y_i$ determines a finite number of possibilities for the value
of $x_i$. In a similar way one proves that, given $x_i$, there are only
finitely many choices for $y$ such that \
$yx_i^{-1} \in \{ y_i x_i^{-1} : i \in I \}$.
Hence the finite-to-finite property is preserved under multiplication.
 \ \ \ $\Box$

\begin{lem}
\label{charactA} \
The algebra ${\cal A}_{\infty}$ consists of the elements
$\sum_{i \in I} \kappa_i \,y_i x_i^{-1}$ of ${\cal B}_{\infty}$ that have
the following three properties:

\smallskip

\noindent (1) The relation $S = \{(x_i, y_i) : i \in I\}$ is
{\em bounded finite-to-finite} (i.e., there is a bound on the cardinalities
of all the sets $S(x_i)$ and $S^{-1}(y_i)$ as $i$ ranges over $I$).

\smallskip

\noindent (2) In $\{x_i : i \in I\}$ and in $\{y_i : i \in I\}$, all
$>_{\rm pref}$-chains have bounded length.

\smallskip

\noindent (3) The set $\{ \kappa_i : i \in I \}$ is finite (i.e., only
finitely many different coefficients occur).
\end{lem}
{\bf Proof.} \ Properties (1), (2) and (3) are straightforward consequences of
the definition of ${\cal A}_{\infty}$.

Conversely, suppose
$\sum_{i \in I} \kappa_i \,y_i x_i^{-1} \ \in {\cal B}_{\infty}$ satisfies
(1), (2) and (3). Let $n_{1,x}$, $n_{1,y}$, $n_{2,x}$ and $n_{2,y}$ be the
bounds that occur in properties (1) and (2); $n_{2,x}$ is the maximum length
of a $>_{\rm pref}$-chain in $\{x_i : i \in I\}$, and similarly for $n_{2,y}$;
$n_{1,x} = {\sf max}\{ |S(x_i)| : i \in I \}$, and
$n_{1,y} = {\sf max}\{ |S^{-1}(y_j)| : j \in I \}$.

We want to prove by induction on $n_{1,x} + n_{1,y} + n_{2,x} + n_{2,y}$ that
$\sum_{i \in I} \kappa_i \,y_i x_i^{-1}$ is a finite linear combination of
elements of $U_{\infty}^{\rm part}$.

Base of the induction: If $n_{2,x} = n_{2,y} = 1$ then $\{ x_i : i \in I \}$
and $\{ y_i : i \in I \}$ are prefix codes.
If in addition $n_{1,x} = n_{1,y} = 1$ then the relation $S$ is injective,
hence $\sum_{i \in I} y_i x_i^{-1} \in U_{\infty}^{\rm part}$.
It follows that $\sum_{i \in I} \kappa_i \,y_i x_i^{-1}$ is a
finite linear combination of elements of $U_{\infty}^{\rm part}$, since
$\{ \kappa_i : i \in I \}$ is finite.

Inductive steps:
If $n_{2,x} > 1$, let $P = \{x_i : i\in J\} \subseteq \{x_i : i \in I\}$
($J \subseteq I$) be a prefix code which is $\subseteq$-maximal in the set
$\{x_i : i \in I\}$. Then \ $\sum_{i \in I} \kappa_i \,y_i x_i^{-1} $
$= \sum_{i \in J} \kappa_i \,y_i x_i^{-1}$
$ + \sum_{i \in I-J} \kappa_i \,y_i x_i^{-1}$.
Since $P$ is maximal within $\{x_i : i \in I\}$, the longest
$>_{\rm pref}$-chain in $\{x_i : i \in I-J\}$ is strictly shorter than the
longest $>_{\rm pref}$-chain in $\{x_i : i \in I\}$.
So, the number $n_{2,x}$ has strictly decreased for
$\sum_{i \in I-J} \kappa_i \,y_i x_i^{-1}$, while the other three numbers did
not increase.  Hence, by induction,
$\sum_{i \in I-J} \kappa_i \,y_i x_i^{-1} \in {\cal A}_{\infty}$.

For the sum $\sum_{i \in J} \kappa_i \,y_i x_i^{-1}$ the number $n_{2,x} =1$;
if for this sum the number $n_{2,y} > 1$, we take a prefix code
$Q = \{y_i : i\in H \} \subseteq \{y_i : i \in I\}$ which is
$\subseteq$-maximal in $\{y_i : i \in I\}$. Then \
$\sum_{i \in J} \kappa_i \,y_i x_i^{-1} = $
$\sum_{i \in H} \kappa_i \,y_i x_i^{-1} + $
$\sum_{i \in J-H} \kappa_i \,y_i x_i^{-1}$.
Then for $\sum_{i \in J-H} \kappa_i \,y_i x_i^{-1}$, the number $n_{2,y}$
(i.e., the maximum length of a $>_{\rm pref}$-chain in $\{y_i : i\in J-H \}$)
has strictly decreased, while the other three numbers did not increase. Hence,
by induction,
$\sum_{i \in J-H} \kappa_i \,y_i x_i^{-1} \in {\cal A}_{\infty}$.

For the sum $\sum_{i \in H} \kappa_i \,y_i x_i^{-1}$, the numbers
$n_{2,x} = n_{2,y} = 1$ (since both $\{ x_i : i\in H\}$ and
$\{ y_i : i \in H\}$ are prefix codes). If for this sum $n_{1,x} > 1$ or
$n_{1,y} > 1$, we write $\sum_{i \in H} \kappa_i \,y_i x_i^{-1}$ as a finite
addition of sums for which the numbers $n_{1,x}$ and $n_{1,y}$ are 1
(while $n_{2,x}$ and $n_{2,y}$ remain 1). Now we are back in the base case.
 \ \ \ $\Box$

\begin{lem}
\label{corresp} \
There is a one-to-one correspondence between (1) the set of all isomorphisms
between (essential) right ideals of $\{a,b\}^*$, and (2) the set
${\cal U}_{\infty}^{\rm part}$ (respectively ${\cal U}_{\infty}$).

Similarly, there is a one-to-one correspondence between (1) the set of all
isomorphisms between finitely generated (essential) right ideals of
$\{a,b\}^*$, and (2) the set ${\cal U}_V^{\rm part}$
(respectively ${\cal U}_V$).
\end{lem}
{\bf Proof.} \ We give the proof for ${\cal U}_{\infty}^{\rm part}$ (and
for ${\cal U}_{\infty}$); for ${\cal U}_V^{\rm part}$ and
${\cal U}_V$ it is similar.  The correspondence map is 

\smallskip

\hspace{1in}
$\Sigma: \ \ \varphi \ \longmapsto \ \sum_{x \in P_1} \varphi(x) \, x^{-1}$

\smallskip

\noindent where $P_1$ (the domain code of $\varphi$) is a prefix code.
This map is clearly onto ${\cal U}_{\infty}^{\rm part}$.
Injectiveness of $\Sigma$ follows from the fact that the map defined
below, is the inverse of $\Sigma$:

\smallskip

\hspace{1in}
$\Phi: \ \ \sum_{i\in I}y_ix_i^{-1}\in {\cal U}_{\infty}^{\rm part}
  \ \longmapsto \
(\varphi : \{x_i:i\in I\} A^* \to \{y_i:i \in I\} A^*)$

\smallskip

\noindent where $\varphi$ is defined by \
$\varphi(x) = \sum_{i\in I} y_i x_i^{-1} \, x  \ \cap \ A^*$. Equivalently,
$\varphi(x) = y_j w$ if $x = x_j w$ for some $j \in I$, $w \in A^*$;
$\varphi(x)$ is undefined otherwise. Then,
since $\{x_i : i \in I\}$ and $\{y_i : i \in I\}$ are prefix codes,
$\varphi$ is an isomorphism between right ideals
(which are essential if $\sum_{i\in I}y_ix_i^{-1}\in {\cal U}_{\infty}$).
  \ \ \ $\Box$

\begin{lem}
\label{Umult} \
The sets ${\cal U}_{\infty}$, ${\cal U}_V$,
${\cal U}_{\infty}^{\rm part}$ and ${\cal U}_V^{\rm part}$ are
closed under multiplication.
\end{lem}
{\bf Proof.} \ Let $\sigma_2 = \sum_{j \in J} y_j \, x_j^{-1}$,
$\sigma_1 = \sum_{i \in I} v_i \, u_i^{-1}$
$\in {\cal U}_{\infty}^{\rm part}$. In the product \ $\sigma_2 \sigma_1$,
each term $y_j \, x_j^{-1} \cdot v_i \, u_i^{-1}$ falls into one of the
following four cases:

\smallskip

\noindent {\bf (0)} \ $x_j$ and $v_i$ are not prefix-comparable: \ Then \
$y_j \, x_j^{-1} \cdot v_i \, u_i^{-1} = {\bf 0}$.

\smallskip

\noindent {\bf (1)} \ $x_j <_{\rm pref} v_i$ \ : \ \ 
Since $\{v_i : i \in I\}$ is a prefix code there will be at most one
$v_i$ for a given $x_j$, such that we are in this case; so, there is a
partial function \ $f: j \in J \mapsto f(j) \in I$ \
such that $x_j <_{\rm pref} v_{f(j)}$. The domain of $f$ is \
dom$f = \{ j \in J : x_j \leq_{\rm pref} v_i \ {\rm for \ some} \ i \in I\}$.
For every $j \in {\rm dom}f$ there is (a unique) $z_j \in A^*$
such that $x_j = v_{f(j)} z_j$. Now, \
$y_j \, x_j^{-1} \cdot v_i \, u_i^{-1} = y_j \, z_j^{-1} \, u_{f(j)}^{-1}$.

\smallskip

\noindent {\bf (2)} \ $x_j >_{\rm pref} v_i$ \ : \ \  
Since
$\{ x_j : j \in J\}$ is a prefix code there will be at most one $x_j$ for a
given $v_i$, such that we are in this case; so, there is a partial function \
$g: i \in I \mapsto g(i) \in J$ \ such that $x_j >_{\rm pref} v_i$.
The domain of $g$ is \
dom$g = \{ i \in I : x_j \geq_{\rm pref} v_i \ {\rm for \ some} \ j \in J\}$.
Then there
is (a unique) $t_i \in A^*$ such that $v_i = x_{g(i)} t_i$. Now, \
$y_j \, x_j^{-1} \cdot v_i \, u_i^{-1} = y_{g(i)} \, t_i \, u_i^{-1}$.

\smallskip

\noindent {\bf (3)} \ $x_j = v_i$ \ : \ \   
 Then \ $y_j \, x_j^{-1} \cdot v_i \, u_i^{-1} = y_j \, u_i^{-1}$
Moreover, for a given $v_i$ there is at most one $x_j$ such that $x_j = v_i$.
The same reasoning as in the two previous cases applies here; we can write
$i = f(j)$ and $j = g(i)$, i.e., this case corresponds to
$j \in {\rm dom}f \cap {\rm im}g$, or equivalently,
$i \in {\rm dom}g \cap {\rm im}f$. Now,
 \ $y_j \, x_j^{-1} \cdot v_i \, u_i^{-1} = y_j \, u_{f(j)}^{-1}$.

\medskip

\noindent This yields the following formula for the multiplication in
${\cal U}_{\infty}^{\rm part}$:

\[ \sum_{j \in J} y_j \, x_j^{-1} \cdot \sum_{i \in I} v_i \, u_i^{-1}
 \ = \ \]
\[ \sum_{j \ \in \ {\rm dom}f \ \cap \ {\rm im}g} \hspace{-.3in}
                    y_j \, u_{f(j)}^{-1}
 \ +  \
\sum_{j \ \in \ {\rm dom}f \ - \ {\rm im}g} \hspace{-.3in}
                     y_j \, z_j^{-1} \, u_{f(j)}^{-1}
 \ + \
\sum_{i \ \in \ {\rm dom}g \ - \ {\rm im}f} \hspace{-.3in}
                       y_{g(i)} \, t_i \, u_i^{-1} \]

Let us check that the domain code \
$\{u_{f(j)}z_j : j \in {\rm dom} f\} \ \cup \ \{u_i : i \not\in {\rm im} f\}$
 \ of $\sigma_2 \sigma_1$ is a prefix code, which is maximal if
$\sigma_2, \sigma_1 \in  {\cal U}_{\infty}$.
(For the image code \
$\{y_{g(i)}t_i : i \in {\rm dom} g\} \ \cup \ \{y_j : j \not\in {\rm im} g\}$
 \ of $\sigma_2 \sigma_1$ the proof is similar.)
Indeed, for each $i \in {\rm im}f$ (this corresponds to cases (1) and (3)),
the set \

\smallskip

 \ \ \ \ \ \ \ \ \ \ \
 $Z_i \ = \ \{ z_j : f(j) = i\} \ = \ \overline{v_i} \, \{x_j : j \in J\}$

\smallskip

\noindent is a prefix code of $A^*$, which is maximal if $\{x_j : j \in J\}$
is a maximal prefix code; see Lemma \ref{Q} for the notation
$\overline{v_i} \, \{ \ldots \}$.
Similarly, for each $j \in {\rm im}g$ (corresponding to cases (2) and (3)),
the set \ $T_j = \ \{t_i : g(i)=j\} \ = \ \overline{x_j} \, \{v_i : i \in I\}$
 \ is a prefix code of $A^*$, which is maximal if \ $\{v_i : i \in I\}$ is a
maximal prefix code.
This follows from Lemma \ref{Q}, since $\{x_j : j \in J\}$ and
$\{v_i : i \in I\}$ are prefix codes.

Now, observe that \
$\{u_{f(j)}z_j : j \in {\rm dom} f\} \ \cup \ \{u_i : i \not\in {\rm im} f\}$
$ \ = \ \bigcup_{i \in I} u_i Z_i$.
By construction (2) in Example \ref{exPrefCod} this is a prefix code, which is
maximal if $\{ u_i : i \in I\}$ and each $Z_i$ are maximal (moreover, each
$Z_i$ will be maximal if $\{x_j : j \in J\}$ is maximal). Similarly, \
$\{y_{g(i)}t_i : i \in {\rm dom} g\} \ \cup \ \{y_j : j \not\in {\rm im} g\}$
$ \ = \ \bigcup_{j \in J} y_jT_j$ \ is a prefix code, which is maximal if
$\{y_j : j \in J\}$ and each $T_j$ are maximal (moreover, each $T_j$ will
be maximal if $\{v_i : i \in I\}$ is maximal).

Now, from the multiplication formula and the fact that the domain and image 
codes in that formula are indeed prefix codes, it follows that
${\cal U}_{\infty}^{\rm part}$ is closed under multiplication. Since the
domain and image codes in the formula are maximal prefix codes if the sums
are in ${\cal U}_{\infty}$, it follows that ${\cal U}_{\infty}$ is also
closed under multiplication.
For ${\cal U}_V^{\rm part}$ and ${\cal U}_V$ the proofs are
similar. \ \ \ $\Box$

\begin{lem}
\label{homU} \
The one-to-one correspondence $\Phi$ of Lemma \ref{corresp} is a
{\em homomorphism}, i.e., for all $\sigma_2$, $\sigma_1$
$ \in {\cal U}_{\infty}^{\rm part} :$ \ \
$\Phi(\sigma_2 \cdot \sigma_1) = \Phi(\sigma_2) \circ \Phi(\sigma_1)$.
\end{lem}
{\bf Proof.} Let $\sigma_2 = \sum_{j \in J} y_j x_j^{-1}$,
$\sigma_1 = \sum_{i \in I} v_i u_i^{-1}$,  and
$\pi = \sum y_j u_{f(j)}^{-1} + \sum y_j z_j^{-1} u_{f(j)}^{-1} + $
$\sum y_{g(i)} t_i u_i^{-1}$, as in the multiplication formula of Lemma 
\ref{Umult}.  We also saw in Lemma \ref{Umult} that \ $D = $
$\{u_{f(j)}z_j : j \in {\rm dom} f\} \ \cup \ \{u_i : i \not\in {\rm im} f\}$
 \ is a prefix code of $A^*$ (which is
maximal if $\sigma_2, \sigma_1 \in {\cal U}_{\infty}$).

It is straightforward to check that for all $x \in DA^*$, \
$\Phi(\pi)(x) = (\Phi(\sigma_2) \circ \Phi(\sigma_1))(x)$,
and both sides are defined.

For all $x \not\in DA^*$, \ $\Phi(\pi)(x)$ \ is undefined;
so, to complete the proof we still must show that
$(\Phi(\sigma_2) \circ \Phi(\sigma_1))(x)$  is undefined when
$x \not\in DA^*$. Note that $x \not\in DA^*$ iff $x >_{\rm pref} d$ for some
$d \in D$, or $x$ is not prefix-comparable with any word in $D$.
More technically, $\Phi(\sigma_2)(\Phi(\sigma_1)(x))$ is undefined iff
for all $i \in I$, $j \in J$:
either we have \ $v_i u_i^{-1} \, x \not\in A^*$, or we have \
$v_i u_i^{-1} \, x = s \in A^*$ \ but \ $y_j x_j^{-1} \, s \not\in A^*$.

\medskip

\noindent {\sf Case 1:} \ $x >_{\rm pref} d$ for some $d \in D$.

Let $\{ d_k : k \in K \} \subset D$ be the elements of $D$ that have $x$ as a
prefix; all other elements of $D$ are prefix-incomparable with $x$
(since $D$ is a prefix code).

\noindent $\bullet$  If $x >_{\rm pref} d_k$ and $d_k$ is of the form $u_i$
(with $v_i = x_j$ or $v_i <_{\rm pref} x_j$), or if $d_k$ is of the form
$u_i z_j$ and $x >_{\rm pref} u_i = xs_i$ (for some non-empty word
$s_i \in A^*$), then  \ $v_i u_i^{-1} \cdot x = v_i s_i^{-1} \not\in A^*$.

\noindent $\bullet$  If $d_k$ is of the form $u_i z_j$ (with
$v_i >_{\rm pref} x_j = v_i z_j$), and
$u_i >_{\rm pref} x = u_i s_i  >_{\rm pref} u_i z_j$, then \
$v_i u_i^{-1} \cdot x = v_i s_i$, and
$y_j x_j^{-1} \cdot v_i s_i = y_j z_j^{-1} s_i = {\bf 0}$ (since $s_i$ and
$z_j$ are not prefix-comparable).

So, in all these cases, \ $\Phi(\sigma_2)(\Phi(\sigma_1)(x))$ is undefined.

\medskip

\noindent {\sf Case 2:} \ $x$ is not prefix-comparable with any word in $D$.

If $x$ is not prefix-comparable with any word $u_i$ then $\Phi(\sigma_1)(x)$
is undefined ($v_i u_i^{-1} \cdot x = {\bf 0}$ for all $i$).

If $x \geq_{\rm pref} u_i$ for some $u_i$ then $x$ is prefix-comparable with
some word in $D$, but this cannot happen in Case 2.

If $x <_{\rm pref} u_i$ for some $u_i$ but $x$ is not prefix-comparable with
any word $u_iz_j \in D$, then for those words $u_i$ we have: $x = u_is_i$
where $s_i$ is not prefix-comparable with $z_j$. Then \
$v_i u_i^{-1} \, x = v_i s_i \in A^*$, but \ $y_j x_j^{-1} \, v_i s_i =$
$y_j (v_iz_j)^{-1} v_i s_i = y_j z_j^{-1} s_i = {\bf 0}$  \ (since $s_i$ is not
prefix-comparable with $z_j$).  \ \ \ $\Box$

\bigskip

It will be useful to extend the definition of the map $\Phi$ of Lemma
\ref{corresp} to all of ${\cal A}_{\infty}$. First some notation:
For a set $S$ and a field $\mathbb{K}$ we denote the set of all
{\it finite $\mathbb{K}$-multisets over} $S$ by $\mathbb{K}[S]$;
such a multiset has the form \
$\{ \kappa_j s_j : j = 1, \ldots , n\}$, with $n \in \mathbb{N}$,
$\kappa_j \in \mathbb{K}$ and $s_j \in S$.
For any \  $\sum_{i\in I} \kappa_i y_i x_i^{-1} \in {\cal A}_{\infty}$ \
we define \

\smallskip

$\Phi(\sum_{i\in I} \kappa_i y_i x_i^{-1}) \ = \ $
$(\varphi: \{x_i:i\in I\} A^* \to \mathbb{K}[\{y_i:i \in I\} A^*])$ \

\smallskip

\noindent where $\varphi$ is defined by \

\smallskip

$\varphi(x) \ = \ $
$\{ \kappa_i y_i x_i^{-1} \, x : i\in I \} \ \cap \ \mathbb{K}[A^*]$
$ \ = \ $
$\{ \kappa_i \, y_i \, w_i : i \in I, x_i \geq_{\rm pref} x = x_i w_i \}$.

\begin{lem}
\label{respectEquiv} \
The one-to-one correspondence $\Sigma$ of Lemma \ref{corresp}
respects the congruence relations on the set ${\cal U}_{\infty}$ (induced
by {\bf I}$_{\infty}$) and on the set of all isomorphisms between essential
right ideals
(i.e., two isomorphisms $\varphi_1$ and $\varphi_2$ between essential right 
ideals are congruent \ iff \   
$\Sigma(\varphi_1)$ and $\Sigma(\varphi_2)$ are congruent).
A similar fact holds for ${\cal U}_V$.
\end{lem}
{\bf Proof.} \ We only prove this for for ${\cal U}_{\infty}$;
for ${\cal U}_V$ it is similar.
Suppose $\varphi_1$ and $\varphi_2$ are congruent, i.e., $\varphi_2$ can be
obtained from $\varphi_1$ by extensions and restrictions.
Then, by Lemmas \ref{ExtCritInf} and \ref{findMaxExt}, any extension or
restriction of $\varphi_1$ can be carried out by repeatedly
(perhaps infinitely often) using maximal prefix codes $Q$, as in
Lemma \ref{findMaxExt}. This corresponds precisely to applying (perhaps
infinitely many) relations of the form \
$\sum_{q \in Q} q q^{-1} \to {\bf 1}$ \ for an extension, or with ``$\to$''
replaced by ``$\leftarrow$'' for a restriction. Hence $\Sigma(\varphi_1)$
and $\Sigma(\varphi_2)$ are congruent modulo ${\bf I}_{\infty}$.

Conversely, if $\Sigma(\varphi_1)$ and $\Sigma(\varphi_2)$ are congruent then
their difference is an element of ${\bf I}_{\infty}$, i.e., \
$\Sigma(\varphi_1) - \Sigma(\varphi_2)$ \ is a linear combination of elements
of the form \

\medskip

 \ \ \ \ \  \ \ \ \ \
$\sum_{j \in J} v_j u_j^{-1}$ $\cdot$
$\sum_{i \in I} y_i (\sum_{q \in Q_i} q q^{-1} - {\bf 1}) x_i^{-1}$
$\cdot$ $\sum_{k \in K} t_k s_k^{-1}$,

\medskip

\noindent
where $\{x_i, i \in I\}, \{y_i, i \in I\}$, and every $Q_i$ are maximal prefix
codes, and where \ $\sum_{j \in J} v_j u_j^{-1}$, $\sum_{k \in K} t_k s_k^{-1}$
$ \in $ ${\cal U}_{\infty}^{\rm part}$.

Since $\{x_i, i \in I\}, \{y_i, i \in I\}$, and every $Q_i$ are maximal prefix
codes, it follows that \  $\bigcup_{i \in I} x_iQ_i$ \ and \
$\bigcup_{i \in I} y_iQ_i$ \  are maximal prefix codes, by construction (2)
of Example \ref{exPrefCod}.
It is straightforward to check that \
$\sum_{i \in I} y_i (\sum_{q \in Q_i} q q^{-1} - {\bf 1}) x_i^{-1}$ \
acts as the empty map on the essential right ideal \
$\bigcup_{i \in I} x_iQ_i \, A^*$,
i.e., for every $x = x_{i_0} q_0 w$ (with $i_0 \in I, q_0 \in Q_{i_0}$ and
$w \in A^*$), \
$\sum_{i \in I} y_i (\sum_{q \in Q_i} q q^{-1} - {\bf 1}) x_i^{-1} \ \cdot x$
$ \ = \ {\bf 0}$. \ Therefore,
 \ $\sum_{i \in I} y_i (\sum_{q \in Q_i} q q^{-1} - {\bf 1}) x_i^{-1}$
$\cdot$ $\sum_{k \in K} t_k s_k^{-1}$ \ acts as the empty map on the set

\medskip

 \ \ \ \ \  \ \ \ \ \
$\Phi(\sum_{k \in K} t_k s_k^{-1})^{-1}(\bigcup_{i \in I} x_iQ_iA^*)$
$ \ = \ \Phi(\sum_{k \in K} s_k t_k^{-1})(\bigcup_{i \in I} x_iQ_iA^*)$.

\medskip

\noindent
This set is a right ideal (being the image of the right-ideal isomorphism
 $\Phi(\sum_{k \in K} s_k t_k^{-1}$)), but not necessarily an essential right
ideal. Let $RA^*$ be any essential right ideal containing  \
$\Phi(\sum_{k \in K} s_k t_k^{-1})(\bigcup_{i \in I} x_iQ_i A^*)$. Then \
$\sum_{i \in I} y_i (\sum_{q \in Q_i} q q^{-1} - {\bf 1}) x_i^{-1}$
$\cdot$ $\sum_{k \in K} t_k s_k^{-1}$ \
acts as the empty map on $RA^*$, and so does
 \ $\sum_{j \in J} v_j u_j^{-1}$ $\cdot$
$\sum_{i \in I} y_i (\sum_{q \in Q_i} q q^{-1} - {\bf 1}) x_i^{-1}$
$\cdot$ $\sum_{k \in K} t_k s_k^{-1}$.

Thus, $\varphi_1$ and $\varphi_2$ agree on the essential right ideal $RA^*$.
Therefore \
${\sf max} \, \varphi_1 = {\sf max} \, \varphi_2$ \ by uniqueness of the
maximum extension (Lemma \ref{uniqueMaxExt}).          \ \ \ $\Box$

\bigskip

\bigskip

\bigskip

\noindent {\bf Acknowledgements.} \ I first learned about Cuntz algebras from
{\it John Meakin} who noticed that the algebras I used for representing the 
Thompson groups are the Cuntz algebras (when completed). I also would like to 
thank a referee for very conscientious and insightful work.  


\bigskip

\bigskip

\noindent
{\bf Jean-Camille Birget} \\ 
Dept.\ of Computer Science \\
Rutgers University at Camden \\
Camden, NJ 08102, USA \\
{\tt birget@camden.rutgers.edu }

\end{document}